\documentclass[12pt]{article}
\usepackage{amsmath,amscd,amsbsy,amssymb,latexsym,url,bm,amsthm}
\usepackage[vlined,boxed,commentsnumbered,linesnumbered,ruled]{algorithm2e}
\usepackage{epsfig,graphicx,subfigure}
\usepackage{enumitem,balance,mathtools}
\usepackage{wrapfig}
\usepackage{mathrsfs, euscript}
\usepackage{xcolor}
\usepackage{hyperref}
\usepackage{epstopdf}
\usepackage{float}
\usepackage{comment}
\usepackage{abstract}

\newcommand{\R}{\mathbb R}
\newcommand{\Z}{\mathbb Z}
\newcommand{\N}{\mathbb N}

\author{}
\title{Well-Posedness of The Compressible Boundary Layer Equations with Analytic Initial Data}
\date{}
\begin{document}
	\maketitle
        \vspace{-2cm}
\newtheorem{theorem}{Theorem}[section]
\newtheorem{lemma}[theorem]{Lemma}
\newtheorem{proposition}[theorem]{Proposition}
\newtheorem{assumption}[theorem]{Assumption}
\newtheorem{corollary}[theorem]{Corollary}
\newtheorem{definition}[theorem]{Definition}
\newtheorem{remark}{Remark}[section]
\newtheorem{example}[theorem]{Example}
\newtheorem{exercise}{Exercise}
\newenvironment{solution}{\begin{proof}[Solution]}{\end{proof}}
	\setlength{\lineskiplimit}{2.625bp}
	\setlength{\lineskip}{2.625bp}
	\numberwithin{equation}{section}
	\newenvironment{partlist}[1][]
	{\begin{enumerate}[itemsep=0pt, label=(\arabic*), wide, labelindent=\parindent, listparindent=\parindent, #1]}
		{\end{enumerate}}
	\setcounter{page}{1}

\begin{center}

Ya-Guang Wang\footnote{email address: ygwang@sjtu.edu.cn}\\
School of Mathematical Sciences, Center for Applied Mathematics, MOE-LSC and SHL-MAC, Shanghai Jiao Tong University, 200240 Shanghai, China\\[3mm]
Yi-Lei  Zhao\footnote{corresponding author, email address:  zhaoyilei@outlook.com}\\
School of Mathematical Sciences, Shanghai Jiao Tong University, 200240 Shanghai, China

\end{center}

\vspace{.1in}

\begin{abstract}
 We study the well-posedness of the compressible  boundary layer equations with data being analytic in  the tangential variable of the boundary. The compressible boundary layer equations, a nonlinear coupled system of degenerate parabolic equations and an elliptic equation, describe the behavior of thermal layer and viscous layer in the small viscosity and heat conductivity limit, for the two-dimensional compressible viscous flow with heat conduction with nonslip and zero heat flux  boundary conditions.  We use the Littlewood-Paley theory to establish the a priori estimates for solutions of this compressible boundary layer problem, and obtain the local existence and uniqueness of the solution in the space of analytic in the tangential variable and Sobolev in the normal variable. 
\end{abstract}

 \tableofcontents

\section{Introduction}
The behavior of the boundary layer in a two-dimensional compressible non-isentropic flow can be governed by the following initial-boundary value problem for a nonlinear coupled system of degenerate parabolic equations and an elliptic equation in $\{(t,x,y)| t>0, x\in  \mathbb{R}, y>0\}$,
\begin{equation}\label{78}
    \begin{cases}\partial_tu+u\partial_xu+v\partial_yu+\frac{RP_x}{P}\theta=\frac{\gamma R}{P}\theta\partial_y^2u,\\\partial_t\theta+u\partial_x\theta+v\partial_y\theta=\frac{R}{(R+c_V)P}\theta\Big[\kappa\partial_y^2\theta+\gamma(\partial_yu)^2+P_x\cdot u+P_t\Big],\\\partial_xu+\partial_yv=\frac{R}{(R+c_V)P}\Big[\kappa\partial_y^2\theta+\gamma(\partial_yu)^2\Big]-\frac{c_V}{(R+c_V)P}\Big(P_x\cdot u+P_t\Big),\\
(u,v)|_{y=0}=0,\quad\partial_y\theta|_{y=0}=0, \quad  \lim\limits_{y\to+\infty}(u,\theta)(t,x,y)=(u^e,\theta^e)(t,x),\\(u,\theta)|_{t=0}=(u_0,\theta_0)(x,y),\end{cases}
\end{equation}
where $u,v,\theta$ represent the tangential velocity, normal velocity, temperature of fluid respectively, and $P(t,x)=R\rho\theta,u^e(t,x),\theta^e(t,x)$ are the pressure, velocity, temperature of the outflow respectively, which satisfy Bernoulli's law:
$$\begin{cases}u^e_t+u^e u^e_x+\frac{R\theta^e}{P}P_x=0,\\
\theta^e_t+u^e\theta^e_x-\frac{R\theta^e}{(R+c_V)P}\cdot(P_t+u^e P_x)=0.
\end{cases}$$
with $R$ being a positive constant related to the ideal gas law, and $c_V$ the specific heat capacity. The viscosity and heat conductivity satisfy $\gamma>0$ and $\kappa>0$. The boundary layer problem \eqref{78} is derived when one considers the small viscosity and heat conductivity limit for the two-dimensional compressible Navier-Stokes equations with heat conduction and non-slip boundary condition for the velocity in the case that the viscosity and heat conductivity have the same scale, see \cite{liu2021study} for the detail. It describes the interaction behavior of  strong viscous layer and thermal layer.

The study of the invisid limit of viscous flow is a classical problem in fluid dynamics. Prandtl \cite{prandtl1904uber} introduced the boundary layer theory for incompressible viscous flows with nonslip boundary conditions, derived what is now known as the Prandtl equation, describing the rapid change of the flow in the thin layer of the boundary. The first local well-posedness result of classical solutions to the two-dimensional Prandtl equation was obtained by Oleinik and her collaborators (\cite{oleinik1963, oleinik1999mathematical}), under the monotonicity assumption of tangential velocity with respect to the normal variable, by using the Crocco transformation. This result was re-obtained by Alexandre, Wang, Xu and Yang in \cite{alexandre2015well}, and Masmoudi and Wong in \cite{masmoudi2015local} in the Sobolev spaces by introducing a direct energy method. With an additional favorable condition on the pressure, Xin and Zhang \cite{xin2004global} obtained the global existence of weak solutions to this problem  by using the Crocco transformation, and the regularity of solutions was studied recently by Xin, Zhang and Zhao in \cite{xin2024global}.  
Without the monotonicity assumption, the unsteady Prandtl equation is ill-posed usually in finite order Sobolev spaces.
E and Enquist \cite{weinan1997blowup} showed that smooth solutions of the unsteady Prandtl equations blow up in finite time when the monotonicity assumption fails, see also Kukavica, Vicol and Wang's result on formulation of the van Dommelen and Shen singularity in \cite{Fei2017}. The linearized instability in Sobolev spaces for the Prandtl equations at a shear flow with tangential velocity profile having a non-degenerate critical point was obtained by Gérard-Varet and Dormy in \cite{gerard2010ill}. The characteristics of degeneracy and loss of derivative in the Prandtl equations lead to study this problem alternatively in analytic setting or Gevrey class spaces. The first local well-posedness in the analytic setting was obtained by Sammartino and Caflisch in \cite{sammartino1998zero}, and later improved in \cite{lombardo2003well} by requiring the analyticity only in the tangential variable by applying the abstract Cauchy–Kowalewskaya (CK) theorem. The almost global existence of analytic solutions with small data was obtained in Zhang and Zhang \cite{zhang2016long}, and Ignatova and Vicol \cite{ignatova2016almost}. Recently, Paicu and Zhang \cite{paicu2021global} proved the global existence of small analytic solutions of the Prandtl equations by means of the Littlewood-Paley theory. Gérard-Varet and Masmoudi \cite{gerard2015well} obtained the well-posedness result for a class of data with $G^{\frac 74}$ regularity, which was improved to hold in $G^2$ by Dietert and Gérard-Varet in \cite{dietert2019well}, Li and Yang in \cite{li-yang}, and Li, Masmoudi and Yang in \cite{li-mas-yang} for the three dimensional problem respectively.  

The results mentioned above pertain to the incompressible flow. However, in many physical scenarios, such as air flow, density and temperature significantly influence the behavior of fluids. Consequently, it is important to understand the role that the change of density and temperature plays in the development of the  boundary layer. However, there have been few results on the problems of the boundary layer in compressible flow. The small viscosity limit for isentropic compressible flow with Navier slip boundary condition was rigorously studied by Wang and Williams in \cite{wang2012inviscid}, in which the boundary layer is weak.  The local existence of classical solutions to the two-dimensional boundary layer equation from the isentropic compressible flow with non-slip boundary condition was obtained in \cite{wang2015local} under the Oleinik monotocity assumption by using the energy method, which was also studied by Gong, Guo and Wang in \cite{gong2016boundary} by using  the Crocco transformation. Recently, Chen, Huang and Li in \cite{chen2024global} obtained a global existence of solutions of the compressible Prandtl equations with small analytic data.

All the above results were limited to the isentropic flow, in which there does not exist thermal layer. However, when considering the effect of heat conduction, in addition to the velocity boundary layer, the heat produced by both heat conduction and viscous friction within the fluid leads to a rapid change of temperature near the boundary, resulting in the formation of a thermal boundary layer.  The stability of boundary layers in the two-dimensional nonisentropic compressible circularly symmetric  flow was studied by Liu and Wang in \cite{liu-wang2014}.
Recently, Liu, Wang and Yang in \cite{liu2021study} have formally derived problems of boundary layers for the small viscosity and heat conductivity limit in the two-dimensional compressible non-isentropic flows with non-slip condition for the velocity and heat flux condition for the temperature on the boundary, for different scales of viscosity and heat conductivity, and they also obtained the local existence of classical solutions to the boundary problem \eqref{78}  by using an energy approach when the tangential velocity is monotonically increasing in the normal variable, when the viscosity and heat conductivity have the same scale, and both of viscous layer and thermal layer are strong. 

The aim of this work is to study the local well-posedness of the problem \eqref{78} in the analytic setting without the monotonicity assumption. Comparing with the classical Prandtl boundary layer equations, the additional complexity and challenging issues of the problem   
\eqref{78} include the complicated interaction between viscous layer and thermal layer, and the divergence of velocity fields depending on the second order derivative of the temperature and the quadratic term of the derivative of velocity.

To simplify the calculation and notations, assume that the outflow given in (\ref{78}) is an uniform one:
$$(P,u^e,\theta^e)(t,x)\equiv(1,0,\theta^E),$$
for a positive constant $\theta^E$, and also let all $\gamma R,\kappa R,R+c_V$ be one.
In order to homogenize the boundary condition at infinity, by introducing  $\hat{\theta}=\theta-\theta^E$, the problem \eqref{78} is simplified as the following one in 
$\{(t,x,y)| t>0, x\in  \mathbb{R}, y>0\}$,
\begin{equation}\label{1}\left\{
    \begin{aligned}
        &\partial_t u+u\partial_x u+v\partial_y u=(\theta+\theta^E)\partial_y^2 u,\\
        &\partial_t \theta+u\partial_x \theta+v\partial_y \theta=(\theta+\theta^E)\partial_y^2 \theta+(\theta+\theta^E)(\partial_y u)^2,\\
        &\partial_x u+\partial_y v=\partial_y^2\theta+(\partial_y u)^2,\\
        &u|_{y=0}=v|_{y=0}=\partial_y\theta|_{y=0}=0,\quad \lim\limits_{y\to\infty}u=\lim\limits_{y\to\infty}\theta=0,\\
        &(u,\theta)|_{t=0}=(u_0,\theta_0)(x,y)
    \end{aligned}
    \right.
\end{equation}
 where we omit the hat of $\hat{\theta}$ for simplicity.
 
The main result of this paper is the following local well-posedness of the problem \eqref{1} in the space of functions analytic in the $x-$variable, when the initial temperature $\theta_0$ is proper small.

\begin{theorem}\label{50}
For a fixed positive constant $\delta>0$, if the initial data $(u_0,\theta_0)$ satisfy 
\begin{equation}\label{1.3}
(e^{\delta|D|}u_0, e^{\delta|D|}\theta_0)\in B^{1,1}_{\Psi_0},\quad (e^{\delta|D|}\partial_y^2u_0, e^{\delta|D|}\partial_y^2\theta_0)\in B^{\frac{1}{2},0}_{\Psi_0},
\end{equation}
and
\begin{equation}\label{1.4}
\|e^{\delta|D|}\theta_0\|_{B^{1,1}_{\Psi_0}}\leq\epsilon
\end{equation}
for a fixed small constant $0<\epsilon<\theta^E$, with $\Psi_0(y)=\frac{y^2}{16\theta^E}$, 
$B^{s,j}_{\Psi}$ being the weighted Besov space defined in Definition \ref{58}, and $|D|$ denoting the Fourier multiplier in the $x-$variable with symbol $|\xi|$, 
then there exists a positive time $T>0$, so that (\ref{1}) has a unique solution $(u,\theta)$ in $[0, T]$ satisfying
\begin{equation}\label{1.5}
e^{\Phi(t,D)}u\in\widetilde{L}_t^\infty(B^{1,1}_{\Psi}), \quad e^{\Phi(t,D)}\partial_yu\in\widetilde{L}_t^2(B^{1,1}_{\Psi}),\quad 
e^{\Phi(t,D)}\partial_y^2u\in\widetilde{L}_t^\infty(B^{\frac{1}{2},0}_{\Psi}),\quad 
e^{\Phi(t,D)}\partial_y^3u\in\widetilde{L}_t^2(B^{\frac{1}{2},0}_{\Psi})
\end{equation}
and
\begin{equation}\label{1.6}
e^{\Phi(t,D)}\theta\in\widetilde{L}_t^\infty(B^{1,1}_{\Psi}),\quad
e^{\Phi(t,D)}\partial_y\theta\in\widetilde{L}_t^2(B^{1,1}_{\Psi}),
\quad
e^{\Phi(t,D)}\partial_y^2\theta\in\widetilde{L}_t^\infty(B^{\frac{1}{2},0}_{\Psi}),
\quad
e^{\Phi(t,D)}\partial_y^3\theta\in\widetilde{L}_t^2(B^{\frac{1}{2},0}_{\Psi}),
\end{equation}
for any $t\leq T$, 
where $\widetilde{L}_t^\infty (B^{s,j}_{\Psi})$ and $\widetilde{L}_t^2(B^{s,j}_{\Psi})$ are the weighted Besov spaces defined in Definition \ref{48}
with $\Psi(t,y)=\frac{y^2}{16\theta^E(1+t)}$, and $\Phi(t,\xi)=(\delta-\lambda\mu(t))|\xi|$ is given in (\ref{3}).
\end{theorem}

The remainder of this paper is organized as follows. In Section 2, we recall some basic facts of the Littlewood-Paley decomposition and define the function spaces which will be used later. Then, we present a series of a priori estimates of solutions to the problem \eqref{1} in Section \ref{76}. Finally, the proof of the main result, Theorem \ref{50} will be given in the last section.

Before the end of this introduction, we give some notations, which shall be used later.
For $a\lesssim b$, we mean there is a positive constant $C$, which may be diﬀerent from line to line but be independent of $\epsilon$, such that $a\leq Cb$, $(a|b):=\int_0^t\int_{\mathbb{R}_+^2}a(t',x,y)\bar{b}(t',x,y)dxdydt'$ 
stands for the inner product of $a$ and $b$ in $L^2((0, t)\times\mathbb{R}_+^2)$. Finally, $\{d_k\}_{k\in\mathbb{Z}}$ designates a generic element of $\ell^1(\mathbb{Z})$ satisfying $\sum\limits_{k\in\mathbb{Z}}|d_k|=1$.

\section{Littlewood-Paley decomposition and functional spaces}\label{75}
We shall frequently use the Littlewood-Paley decomposition in the tangential variable $x$. First, recall some notations from \cite{bahouri2011fourier}. Let 
 $\chi(\xi)$, $\varphi(\xi)$ be two fixed smooth functions such that
$$\begin{aligned}
    &\mathrm{Supp}\ \varphi\subset\left\{\xi\in\mathbb{R}\big|~\frac{3}{4}\leqslant|\xi|\leqslant\frac{8}{3}\right\}&\mathrm{and}&\quad\sum_{k\in\mathbb{Z}}\varphi(2^{-k}\xi)=1\quad  (\forall\xi>0),\\
    &\mathrm{Supp}\ \chi\subset\left\{\xi\in\mathbb{R}\big|~|\xi|\leqslant\frac{4}{3}\right\}&\mathrm{and}&\quad\chi(\xi)+\sum_{k\geq0}\varphi(2^{-k}\xi)=1\quad (\forall\xi\ge 0).
\end{aligned}$$
For any given Schwartz distribution $a\in {\mathcal S}'(\R_+^2)$, denote by  
$$\begin{aligned}
&\Delta_ka=0\quad\mathrm{if}\quad k\leq -2,\quad\Delta_{-1}a=\mathcal{F}^{-1}(\chi(|\xi|)\hat{a}),\quad\Delta_ka=\mathcal{F}^{-1}(\varphi(2^{-k}|\xi|)\hat{a})\quad\mathrm{if}\quad k\geqslant 0,\\&\mathrm{and}\quad S_ka=\mathcal{F}^{-1}(\chi(2^{-k}|\xi|)\hat{a}), \end{aligned}$$
where $\hat{a}$ and $\mathcal{F}^{-1}(a)$ denote the partial Fourier transform and Fourier inverse transform of $a$ with respect to $x$ variable, that is, $\hat{a}(\xi,y)=\mathcal{F}_{x\rightarrow\xi}(a)(\xi,y)$.
The following Bony's decomposition shall be used constantly,
\begin{equation}\label{8}
    fg=T_fg+T_gf+R(f,g),
\end{equation}
where 
$$T_fg=\sum_kS_{k-1}f\Delta_kg,\quad R(f,g)=\sum_k\sum\limits_{k'=k-1}^{k+1}\Delta_{k'}f\Delta_kg.$$

\begin{definition}\label{58}
    For any $s\in\R,j\in\N$, and $u\in\mathcal{S}'(\R_+^2)$ satisfying $\lim\limits_{k\rightarrow-\infty}\|S_ku\|_{L^\infty}=0$, we define
    $$B^{s,j}(\R_+^2):=\{u\in\mathcal{S}'(\R_+^2)|\|u\|_{B^{s,j}}<\infty\},$$
  with  
    $$\|u\|_{B^{s,j}}:=\|(2^{ks}\sum\limits_{i=0}^j\|\Delta_k\partial_y^iu\|_{L_+^2})_k\|_{\ell^1(\Z)},$$
   and the notation $\|a\|_{L_+^2}=\|a\|_{L^2(\{y>0\})}$, and
    $$B^{s,j}_\Psi(\R_+^2):=\{u\in\mathcal{S}'(\R_+^2)|\|u\|_{B^{s,j}_\Psi}<\infty\}$$
    with 
    $$\|u\|_{B^{s,j}_\Psi}:=\|(2^{ks}\sum\limits_{i=0}^j\|e^\Psi\Delta_k\partial_y^iu\|_{L_+^2})_k\|_{\ell^1(\Z)}$$
    for a given weighted function $\Psi(y)$.
\end{definition}

\begin{definition}\label{48}
(1)    For fixed $p\in[1,+\infty]$ and $t>0$, define $\widetilde{L}^p(0,t;B^{s,j}_\Psi(\R_+^2))$ as the completion of $C([0,t];\mathcal{S}(\R_+^2))$ by the norm for $1\le p<\infty$   $$\|u\|_{\widetilde{L}^p_t(B^{s,j}_\Psi)}:=\sum\limits_{k\in\mathbb{Z}}\sum\limits_{i=0}^j2^{ks}\left(\int_0^t\|e^\Psi\Delta_k\partial_y^iu(t',\cdot)\|^p_{L_+^2}dt'\right)^{\frac{1}{p}},$$
    with the usual change as $p=\infty$.

(2) For a given non-negative function $f(t)\in L^1_{\mathrm{loc}}(\mathbb{R}_+)$, define
    $$\|u\|_{\widetilde{L}^p_{t,f}(B^{s,j}_\Psi)}:=\sum\limits_{k\in\mathbb{Z}}\sum\limits_{i=0}^j2^{ks}\left(\int_0^tf(t')\|e^\Psi\Delta_k\partial_y^iu(t',\cdot)\|^p_{L_+^2}dt'\right)^{\frac{1}{p}},$$
     with the usual change as $p=\infty$.
 \end{definition}

In the following discussion, we shall fix the weighted function $\Psi(t,y)$ as
\begin{equation}\label{74}\Psi(t,y)=\frac{y^2}{16\theta^E\langle t\rangle},\end{equation}
with $\langle t\rangle=1+t$, and $\Psi_0(y):=\Psi(0,y)=\frac{y^2}{16\theta^E}.$
Obviously, $\Psi(t,y)$ satisfies 
\begin{equation}\label{43}
    \partial_t\Psi+4\theta^E(\partial_y\Psi)^2=0.
\end{equation}
Once $\Psi$ is fixed, we will omit the subscription $\Psi$ of $B^{s,j}_\Psi$ for simplicity. 

To control the terms consisting of $\partial_y\Psi$ resulting from integration by parts, we recall the following lemma, which is a special case of Lemma 2.5 given in \cite{wang2024global}:
\begin{lemma}\label{126}
    Let $u(x,y)$ be a smooth function on $\R_+^2$ decaying to zero sufficiently fast as $y\to +\infty$. Then, one has 
    \begin{equation}\label{127}
        \int_{\mathbb{R}_+^2}|\partial_y\Psi u(x,y)|^2e^{2\Psi}dxdy\leq C\int_{\mathbb{R}_+^2}|\partial_y u(x,y)|^2e^{2\Psi}dxdy,
    \end{equation}
    where $C$ is a constant depending only on $\theta^E$.
\end{lemma}

For any given locally bounded function $\Phi(t,\xi)$ on $\R^+\times\R$, define
\begin{equation}\label{2}
u_\Phi(t,x,y)=\mathcal{F}^{-1}_{\xi\rightarrow x}(e^{\Phi(t,\xi)}\hat{u}(t,\xi,y)).
\end{equation}

A key quantity $\mu(t)$ for describing the evolution of the analytic radius of $u,\theta$ is determined by
\begin{equation}\label{46}\left\{
    \begin{aligned}
        &\dot{\mu}=\langle t\rangle^{\frac{1}{4}}\|(\partial_yu_\Phi,\partial_y\theta_\Phi)\|_{B^{\frac{1}{2},2}}+\langle t\rangle^{\frac{1}{2}}(\|u_\Phi\|_{B^{1,0}}^2+\|\theta_\Phi\|_{B^{\frac{1}{2},0}}^2+\|(\partial_yu_\Phi,\partial_y\theta_\Phi)\|_{B^{\frac{1}{2},1}}^2)+\|\partial_yu_\Phi\|_{B^{\frac{1}{2},1}}^4\\
        &\  \
        \quad+\|\theta_\Phi\|_{B^{\frac{1}{2},1}}^4+\|\partial_y^3u_\Phi\|_{B^{\frac{1}{2},0}}(\|\partial_yu_\Phi\|_{B^{\frac{1}{2},1}}+\|\theta_\Phi\|_{B^{\frac{1}{2},1}})+\|\partial_y^3\theta_\Phi\|_{B^{\frac{1}{2},0}}\|\partial_yu_\Phi\|_{B^{\frac{1}{2},0}},\\
        &\mu|_{t=0}=0.
    \end{aligned}
\right.\end{equation}
To overcome the difficulty of loss of derivatives, inspired by the idea given in \cite{zhang2016long}, the phase function $\Phi$ is chosen as
\begin{equation}\label{3}\Phi(t,\xi)=(\delta-\lambda\mu(t))|\xi|,\end{equation}
where $\lambda$ is a sufficiently large constant.

In what follows, we shall always assume that $T^*$ is determined by
\begin{equation}\label{73}T^*:=\sup\{t>0|~\mu( t)<\delta/\lambda\}.\end{equation}
Obviously, when  $0<t<T^*$, there holds the following convex inequality
\begin{equation}\label{9}\Phi(t,\xi)\leqslant\Phi(t,\xi-\eta)+\Phi(t,\eta),\quad \ \forall\xi,\eta\in\R.\end{equation}

\section{A priori estimate}\label{76}

The main goal of this section is to prove the following a priori estimate for the solution of the problem (\ref{1}).

\begin{theorem}\label{42}
Suppose that $(u,\theta)$ is a solution of (\ref{1}) with norms appeared in \eqref{45} being finite. Then, there exists a positive time $T>0$ and a large constant $\lambda>0$, so that for any $0<t<T$ one has
\begin{equation}\label{45}
\begin{aligned}
    &\|(u_\Phi, \theta_\Phi)\|_{\widetilde{L}_t^\infty(B^{1,1})}
+\|(\partial_yu_\Phi,\partial_y\theta_\Phi)\|_{\widetilde{L}_t^2(B^{1,1})}
    +\|(\partial_y^2u_\Phi, \partial_y^2\theta_\Phi)\|_{\widetilde{L}_t^\infty(B^{\frac{1}{2},0})}\\  &+\|(\partial_y^3u_\Phi, \partial_y^3\theta_\Phi)\|_{\widetilde{L}_t^2(B^{\frac{1}{2},0})}+\sqrt{\lambda}(\|(u_\Phi, \theta_\Phi)\|_{\widetilde{L}_{t,\dot{\mu}}^2(B^{\frac{3}{2},1})}+\|(\partial_y^2u_\Phi,
    \partial_y^2\theta_\Phi)  \|_{\widetilde{L}_{t,\dot{\mu}}^2(B^{1,0})})\\
&\lesssim\|
(e^{\delta|D|}u_0,
e^{\delta|D|}\theta_0)
\|_{B^{1,1}_{\Psi_0}}+\|
(e^{\delta|D|}\partial_y^2u_0,e^{\delta|D|}\partial_y^2\theta_0)
\|_{B^{\frac{1}{2},0}_{\Psi_0}}.
\end{aligned}
\end{equation}
\end{theorem}

The key ingredients for proving Theorem \ref{42} are estimates stated in the following six propositions, their proofs shall be given in the following sections.

\begin{proposition}\label{15}
Suppose that $(u,\theta)$ is a smooth solution of (\ref{1}) with norms appeared in \eqref{17} making sense. Then for any $t<T^*$ and $0<\eta<1$, it holds that
\begin{equation}\label{17}
\begin{aligned}
&\|u_\Phi(t)\|_{B^{1,0}}+\sqrt{\lambda}\|u_\Phi\|_{\widetilde{L}_{t,\dot{\mu}}^2(B^{\frac{3}{2},0})}+(\frac{\sqrt{\theta^E}}{2}-\epsilon^{\frac{1}{2}}\langle t\rangle^{\frac{1}{8}})\|\partial_yu_\Phi\|_{\widetilde{L}_t^2(B^{1,0})}\\
\lesssim&\|e^{\delta|D|}u_0\|_{B^{1,0}_{\Psi_0}}+\eta(\|\partial_yu_\Phi\|_{\widetilde{L}_t^2(B^{1,0})}+\|\partial_y\theta_\Phi\|_{\widetilde{L}_t^2(B^{1,0})})+C_\eta\|u_\Phi\|_{\widetilde{L}_{t,\dot{\mu}}^2(B^{\frac{3}{2},0})}.
\end{aligned}\end{equation}
\end{proposition}

\begin{proposition}\label{21}
    Suppose that $(u,\theta)$ is a solution of (\ref{1}) with norms appeared in \eqref{23} making sense. Then for any $t<T^*$ and $0<\eta<1$, one has
    \begin{equation}\label{23}
    \begin{aligned}    &\|\theta_\Phi(t)\|_{B^{1,0}}+\sqrt{\lambda}\|\theta_\Phi\|_{\widetilde{L}_{t,\dot{\mu}}^2(B^{\frac{3}{2},0})}+(\frac{\sqrt{\theta^E}}{2}-\epsilon^{\frac{1}{2}}\langle t\rangle^{\frac{1}{8}})\|\partial_y\theta_\Phi\|_{\widetilde{L}_t^2(B^{1,0})}\\ \lesssim&\|e^{\delta|D|}\theta_0\|_{B^{1,0}_{\Psi_0}}+\eta(\|\partial_yu_\Phi\|_{\widetilde{L}_t^2(B^{1,1})}+\|\partial_y\theta_\Phi\|_{\widetilde{L}_t^2(B^{1,0})})+C_\eta\|\theta_\Phi\|_{\widetilde{L}_{t,\dot{\mu}}^2(B^{\frac{3}{2},0})}+\|u_\Phi\|_{\widetilde{L}_{t,\dot{\mu}}^2(B^{\frac{3}{2},0})}.\\
    \end{aligned}
    \end{equation}
\end{proposition}

\begin{proposition}\label{31}
Suppose that $(u,\theta)$ is a solution of (\ref{1}) with norms appeared in \eqref{33} making sense. Then for any $t<T^*$ and $0<\eta<1$, it holds that
\begin{equation}\label{33}
\begin{aligned}
 &\|\partial_yu_\Phi(t)\|_{B^{1,0}}+\sqrt{\lambda}\|\partial_yu_\Phi\|_{\widetilde{L}_{t,\dot{\mu}}^2(B^{\frac{3}{2},0})}+(\frac{\sqrt{\theta^E}}{2}-\epsilon^{\frac{1}{2}}\langle t\rangle^{\frac{1}{8}})\|\partial^2_y{u}_\Phi\|_{\widetilde{L}_t^2(B^{1,0})}\\
\lesssim&\|e^{\delta|D|}\partial_y{u}_0\|_{B^{1,0}_{\Psi_0}}+\eta(\|\partial^2_y{u}_\Phi\|_{\widetilde{L}_t^2(B^{1,0})}+\|\partial_y\theta_\Phi\|_{\widetilde{L}_t^2(B^{1,1})})+C_\eta\|\partial_yu_\Phi\|_{\widetilde{L}_{t,\dot{\mu}}^2(B^{\frac{3}{2},0})}+\|u_\Phi\|_{\widetilde{L}_{t,\dot{\mu}}^2(B^{\frac{3}{2},0})}.
\end{aligned}\end{equation}
\end{proposition}

\begin{proposition}\label{38}
Suppose that $(u,\theta)$ is a solution of (\ref{1}) with norms appeared in \eqref{41} making sense. Then for any $t<T^*$ and $0<\eta<1$, it holds that
\begin{equation}\label{41}
    \begin{aligned}    &\|\partial_y\theta_\Phi(t)\|_{B^{1,0}}+\sqrt{\lambda}\|\partial_y\theta_\Phi\|_{\widetilde{L}_{t,\dot{\mu}}^2(B^{\frac{3}{2},0})}+(\frac{\sqrt{\theta^E}}{2}-\epsilon^{\frac{1}{2}}\langle t\rangle^{\frac{1}{8}})\|\partial^2_y{\theta}_\Phi\|_{\widetilde{L}_t^2(B^{1,0})}\\ \lesssim&\|e^{\delta|D|}\partial_y \theta_0\|_{B^{1,0}_{\Psi_0}}+\eta(\|\partial^2_y{u}_\Phi\|_{{\widetilde{L}}^2(B^{1,0})}+\|\partial_y\theta_\Phi\|_{\widetilde{L}_t^2(B^{1,1})})+C_\eta\|\partial_y\theta_\Phi\|_{\widetilde{L}_{t,\dot{\mu}}^2(B^{\frac{3}{2},0})}\\
        &+\|u_\Phi\|_{\widetilde{L}_{t,\dot{\mu}}^2(B^{\frac{3}{2},1})}+\|\theta_\Phi\|_{\widetilde{L}_{t,\dot{\mu}}^2(B^{\frac{3}{2},0})}.\\
    \end{aligned}
    \end{equation}

\end{proposition}

\begin{proposition}\label{79}
Suppose that $(u,\theta)$ is a solution of (\ref{1}) with norms appeared in \eqref{80} making sense. Then for any $t<T^*$ and $0<\eta<1$, it holds that
\begin{equation}\label{80}
\begin{aligned}
&\|\partial_y^2u_\Phi(t)\|_{B^{\frac{1}{2},0}}+\sqrt{\lambda}\|\partial_y^2u_\Phi\|_{\widetilde{L}_{t,\dot{\mu}}^2(B^{1,0})}+(\frac{\sqrt{\theta^E}}{2}-\epsilon^{\frac{1}{2}}\langle t\rangle^{\frac{1}{8}})\|\partial^3_yu_\Phi\|_{\widetilde{L}_t^2(B^{\frac{1}{2},0})}\\
\lesssim&\|e^{\delta|D|}\partial_y^2u_0\|_{B^{\frac{1}{2},0}_{\Psi_0}}+\eta(\|\partial^3_yu_\Phi\|_{\widetilde{L}_t^2(B^{\frac{1}{2},0})}+\|\partial_y^2\theta_\Phi\|_{\widetilde{L}_t^2(B^{\frac{1}{2},1})})+C_\eta(\|\partial_y^2u_\Phi\|_{\widetilde{L}_{t,\dot{\mu}}^2(B^{1,0})}+\|\partial_y^2\theta_\Phi\|_{\widetilde{L}_{t,\dot{\mu}}^2(B^{1,0})})\\
&+\|u_\Phi\|_{\widetilde{L}_{t,\dot{\mu}}^2(B^{1,1})}+\|\theta_\Phi\|_{\widetilde{L}_{t,\dot{\mu}}^2(B^{1,1})}.
\end{aligned}\end{equation}
\end{proposition}

\begin{proposition}\label{93}
Suppose that $(u,\theta)$ is a solution of (\ref{1}) with norms appeared in \eqref{94} making sense. Then for any $t<T^*$ and $0<\eta<1$, it holds that
\begin{equation}\label{94}
\begin{aligned}
&\|\partial_y^2\theta_\Phi(t)\|_{B^{\frac{1}{2},0}}+\sqrt{\lambda}\|\partial_y^2\theta_\Phi\|_{\widetilde{L}_{t,\dot{\mu}}^2(B^{1,0})}+(\frac{\sqrt{\theta^E}}{2}-\epsilon^{\frac{1}{2}}\langle t\rangle^{\frac{1}{8}})\|\partial_y^3\theta_\Phi\|_{\widetilde{L}_t^2(B^{\frac{1}{2},0})}\\
\lesssim&\|e^{\delta|D|}\partial_y^2\theta_0\|_{B^{\frac{1}{2},0}_{\Psi_0}}+\eta(\|\partial_y^3\theta_\Phi\|_{\widetilde{L}_t^2(B^{\frac{1}{2},0})}+\|\partial^3_yu_\Phi\|_{\widetilde{L}_t^2(B^{\frac{1}{2},0})})+C_\eta\|\partial_y^2\theta_\Phi\|_{\widetilde{L}_{t,\dot{\mu}}^2(B^{1,0})}+\|u_\Phi\|_{\widetilde{L}_{t,\dot{\mu}}^2(B^{1,2})}\\
&+\|\theta_\Phi\|_{\widetilde{L}_{t,\dot{\mu}}^2(B^{1,1})}+C_\eta(\|\partial_yu_\Phi\|_{\widetilde{L}_{t,\dot{\mu}}^2(B^{\frac{3}{2},0})}+\|\partial_y\theta_\Phi\|_{\widetilde{L}_{t,\dot{\mu}}^2(B^{\frac{3}{2},0})}).
\end{aligned}\end{equation}

\end{proposition}


The proofs of the above propositions shall be given in the next subsections.
Now, we turn to prove the main estimate given in Theorem \ref{42} by using these propositions.

\begin{proof}[Proof of Theorem \ref{42}]
By summing up the estimates given in Propositions \ref{15}-\ref{93}, we obtain that for any $0<\eta<1$,
\begin{equation}\label{44}
\begin{aligned}
&\|u_\Phi(t)\|_{B^{1,1}}+\|\partial_y^2u_\Phi(t)\|_{B^{\frac{1}{2},0}}+\|\theta_\Phi(t)\|_{B^{1,1}}+\|\partial_y^2\theta_\Phi(t)\|_{B^{\frac{1}{2},0}}\\
&+\sqrt{\lambda}(\|u_\Phi\|_{\widetilde{L}_{t,\dot{\mu}}^2(B^{\frac{3}{2},1})}+\|\partial_y^2u_\Phi\|_{\widetilde{L}_{t,\dot{\mu}}^2(B^{1,0})}+\|\theta_\Phi\|_{\widetilde{L}_{t,\dot{\mu}}^2(B^{\frac{3}{2},1})}+\|\partial_y^2\theta_\Phi\|_{\widetilde{L}_{t,\dot{\mu}}^2(B^{1,0})})\\
&+(\frac{\sqrt{\theta^E}}{2}-\epsilon^{\frac{1}{2}}\langle t\rangle^{\frac{1}{8}})(\|\partial_yu_\Phi\|_{\widetilde{L}_t^2(B^{1,1})}+\|\partial_y^3u_\Phi\|_{\widetilde{L}_t^2(B^{\frac{1}{2},0})}+\|\partial_y\theta_\Phi\|_{\widetilde{L}_t^2(B^{1,1})}+\|\partial_y^3\theta_\Phi\|_{\widetilde{L}_t^2(B^{\frac{1}{2},0})})\\
\lesssim&\|e^{\delta|D|}u_0\|_{B^{1,1}_{\Psi_0}}+\|e^{\delta|D|}\partial_y^2u_0\|_{B^{\frac{1}{2},0}_{\Psi_0}}+\|e^{\delta|D|}\theta_0\|_{B^{1,1}_{\Psi_0}}+\|e^{\delta|D|}\partial_y^2\theta_0\|_{B^{\frac{1}{2},0}_{\Psi_0}}\\
&+\eta(\|\partial_yu_\Phi\|_{\widetilde{L}_t^2(B^{1,1})}+\|\partial_y^3u_\Phi\|_{\widetilde{L}_t^2(B^{\frac{1}{2},0})}+\|\partial_y\theta_\Phi\|_{\widetilde{L}_t^2(B^{1,1})}+\|\partial_y^3\theta_\Phi\|_{\widetilde{L}_t^2(B^{\frac{1}{2},0})})\\
&+C_\eta(\|u_\Phi\|_{\widetilde{L}_{t,\dot{\mu}}^2(B^{\frac{3}{2},1})}+\|\partial_y^2u_\Phi\|_{\widetilde{L}_{t,\dot{\mu}}^2(B^{1,0})}+\|\theta_\Phi\|_{\widetilde{L}_{t,\dot{\mu}}^2(B^{\frac{3}{2},1})}+\|\partial_y^2\theta_\Phi\|_{\widetilde{L}_{t,\dot{\mu}}^2(B^{1,0})}).
\end{aligned}\end{equation}

Obviously, we have $$\sup_{t\in[0,\tau_\epsilon^*)}\epsilon^{\frac{1}{2}}\langle t\rangle^{\frac{1}{8}}\le \frac{\sqrt{\theta^E}}{4}\qquad\text{when~}\tau_\epsilon^*\le\left(\frac{\theta^E}{16\epsilon}\right)^{4}.$$
Therefore, we can take a sufficiently small $\eta\leq\frac{\sqrt{\theta^E}}{4}$ so that the terms containing $\eta(\cdot)$ on the right hand side of (\ref{44}) can be absorbed by the left hand side, and by choosing a large $\lambda$ in (\ref{44}) satisfying $\sqrt{\lambda}\gtrsim C_\eta$, the last terms containing 
$C_\eta(\cdot)$ can be controlled by the left hand side too in \eqref{44}.

On the other hand, for the terms given on the right hand side of \eqref{46}, we have
$$\begin{aligned}
&\int_0^t\langle t'\rangle^{\frac{1}{4}}\|(\partial_yu_\Phi, \partial_y\theta_\Phi)\|_{B^{\frac{1}{2},2}}dt'\lesssim(\langle t\rangle^{\frac{3}{2}}-1)^{\frac{1}{2}}\|(\partial_yu_\Phi, \partial_y\theta_\Phi)\|_{\widetilde{L}_t^2(B^{\frac{1}{2},2})},\\
&\int_0^t\langle t'\rangle^{\frac{1}{2}}\left(\|u_\Phi\|_{B^{1,0}}^2+\|(\partial_yu_\Phi, \partial_y\theta_\Phi)\|_{B^{\frac{1}{2},1}}^2\right)dt'\lesssim(\langle t\rangle^{\frac{3}{2}}-1)(\|u_\Phi\|_{\widetilde{L}_t^\infty(B^{1,0})}^2+\|(\partial_yu_\Phi, \partial_y\theta_\Phi)\|_{\widetilde{L}_t^\infty(B^{\frac{1}{2},1})}^2),\\
&\int_0^t(\|\partial_yu_\Phi\|_{B^{\frac{1}{2},1}}^4+\|\theta_\Phi\|_{B^{\frac{1}{2},1}}^4)dt'
\lesssim t(\|\partial_y u_\Phi\|_{\widetilde{L}_t^\infty(B^{\frac{1}{2},1})}^4+\|\theta_\Phi\|_{\widetilde{L}_t^\infty(B^{\frac{1}{2},1})}^4),\\
&\int_0^t \|\partial_y^3u_\Phi\|_{B^{\frac{1}{2},0}}(\|\partial_yu_\Phi\|_{B^{\frac{1}{2},1}}+\|\theta_\Phi\|_{B^{\frac{1}{2},1}})dt'\lesssim t^{\frac{1}{2}}\|\partial_y^3 u_\Phi\|_{\widetilde{L}_t^2(B^{\frac{1}{2},0})}(\|\partial_y u_\Phi\|_{\widetilde{L}_t^\infty(B^{\frac{1}{2},1})}+\|\theta_\Phi\|_{\widetilde{L}_t^\infty(B^{\frac{1}{2},1})}),\\
&\int_0^t \|\partial_y^3\theta_\Phi\|_{B^{\frac{1}{2},0}}\|\partial_yu_\Phi\|_{B^{\frac{1}{2},0}}dt'\lesssim t^{\frac{1}{2}}\|\partial_y^3 \theta_\Phi\|_{\widetilde{L}_t^2(B^{\frac{1}{2},0})}\|\partial_y u_\Phi\|_{\widetilde{L}_t^\infty(B^{\frac{1}{2},0})}.\\
\end{aligned}$$

Hence in view of (\ref{46}) and (\ref{44}), there exists a small $T^\star>0$ such that
$$\sup_{t\in[0,T^\star)}\mu(t)\leq \frac{\delta}{\lambda}.$$
Therefore in view of (\ref{73}), this ensures that $T^*\geq T^\star$.
Thus (\ref{44}) implies that (\ref{45}) holds for any $0<t<T_{\epsilon}$, where $T_{\epsilon}=\min\{\tau_\epsilon^*,T^\star\}$. This completes the proof of Theorem \ref{42}.
\end{proof}

\begin{remark}
From Propositions \ref{21} and \ref{38}, we obtain
$$\|\theta_\Phi\|_{\widetilde{L}_t^\infty({B^{1,1}})}\lesssim\|e^{\delta|D|}\theta_0\|_{B^{1,1}_{\Psi_0}}+\eta\|\partial_yu_\Phi\|_{\widetilde{L}_t^2(B^{1,1})}+\|u_\Phi\|_{\widetilde{L}_{t,\dot{\mu}}^2(B^{\frac{3}{2},1})}.$$
By taking $\eta<\epsilon, \lambda>\epsilon^{-2}$ in Theorem \ref{42}, we get
$$\eta\|\partial_yu_\Phi\|_{\widetilde{L}_t^2(B^{1,1})}+\|u_\Phi\|_{\widetilde{L}_{t,\dot{\mu}}^2(B^{\frac{3}{2},1})}\lesssim\epsilon.$$
As a result, under the assumption \eqref{1.4}, for any $t<T_\epsilon$, 
\begin{equation}\label{51}\|\theta_\Phi\|_{\widetilde{L}_t^\infty({B^{1,1}})}\lesssim\epsilon.
\end{equation}
From the Sobolev embedding, we have
$$\|\theta\|_{L_+^\infty}\lesssim\|\theta\|_{L_+^2}+\|\partial_x\theta\|_{L_+^2}+\|\partial_y\theta\|_{L_+^2}+\|\partial_y\partial_x\theta\|_{L_+^2}.$$
Therefore, by taking $\epsilon\leq\frac{\theta^E}{2}$, we have $\|\theta\|_{\widetilde{L}_t^\infty(L_+^\infty)}\leq\frac{\theta^E}{2}$, which implies
\begin{equation}\label{5}
\theta+\theta^E\geq\frac{\theta^E}{2}>0.
\end{equation}
\end{remark}

\subsection{Estimate of \texorpdfstring{$u$}{u}}

The proposal of this section is to establish the estimate \eqref{17} given in Proposition \ref{15}.

From (\ref{1}) we known that $u_\Phi=
\mathcal{F}^{-1}_{\xi\rightarrow x}(e^{\Phi(t,\xi)}\hat{u}(t,\xi,y))$ with $\Phi(t,\xi)=(\delta-\lambda\mu(t))|\xi|$ satisfies the following equation,
\begin{equation}\label{4}\partial_t u_\Phi+\lambda\dot{\mu}|D| u_\Phi+[u\partial_x u]_\Phi+[v\partial_y u]_\Phi=[(\theta+\theta^E)\partial_y^2 u]_\Phi.\end{equation}

By applying the dyadic operator $\Delta_k$ to (\ref{4}) and taking the $L_{t,+}^2=L^2(0,t;\R_+^2)$ inner product of the resulting equation with $e^{2\Psi}\Delta_ku_\Phi$, we get 
\begin{equation}\label{16}
\begin{aligned}
    &(e^\Psi\partial_t \Delta_k u_\Phi|e^\Psi\Delta_k u_\Phi)+\lambda(\dot{\mu}|D| e^\Psi\Delta_k u_\Phi|e^\Psi\Delta_k u_\Phi)-(e^\Psi\Delta_k[(\theta+\theta^E)\partial_y^2 u]_\Phi|e^\Psi\Delta_k u_\Phi)\\
    =&-(e^\Psi\Delta_k[u\partial_x u]_\Phi|e^\Psi\Delta_k u_\Phi)-(e^\Psi\Delta_k[v\partial_y u]_\Phi|e^\Psi\Delta_k u_\Phi).
\end{aligned}\end{equation}
Denote the three terms on the left side by $\mathcal{L}_1^k,\mathcal{L}_2^k,\mathcal{L}_3^k$ and the two terms on the right side by $\mathcal{R}_1^k,\mathcal{R}_2^k$ respectively.

By using $u|_{y=0}=0$, one has
\begin{equation}\label{3.13}
\begin{aligned}
    \mathcal{L}_3^k=&(e^\Psi\Delta_k[(\theta+\theta^E)\partial_y u]_\Phi|e^\Psi\Delta_k \partial_y u_\Phi)\\
    &+(2e^\Psi\partial_y\Psi\Delta_k[(\theta+\theta^E)\partial_y u]_\Phi|e^\Psi\Delta_k u_\Phi)+(e^\Psi\Delta_k[\partial_y\theta\partial_y u]_\Phi|e^\Psi\Delta_k u_\Phi)\\
    =&:\mathcal{L}_{3,1}^k+\mathcal{L}_{3,2}^k+\mathcal{L}_{3,3}^k.
\end{aligned}\end{equation}

\begin{lemma}\label{6}
 For any $0< t<T^*$ with $T^*$ being given in \eqref{73}, the term $\mathcal{L}_{3,1}^k$ given in \eqref{3.13} satisfies the following estimate $$\sum\limits_{k\in\mathbb{Z}}2^k\sqrt{|\mathcal{L}_{3,1}^k|}\gtrsim (\sqrt{\theta^E}-\epsilon^{\frac{1}{2}}\langle t\rangle^{\frac{1}{8}})\|\partial_yu_\Phi\|_{\widetilde{L}_t^2(B^{1,0})}.$$
\end{lemma}

\begin{proof}
Decompose $\mathcal{L}_{3,1}^k$ into 
    $$\begin{aligned}
        \mathcal{L}_{3,1}^k&=\theta^E(e^\Psi\Delta_k\partial_y u_\Phi|e^\Psi\Delta_k \partial_y u_\Phi)+(e^\Psi\Delta_k[\theta\partial_y u]_\Phi|e^\Psi\Delta_k \partial_y u_\Phi)\\
        &=: I^k_1+I^k_2.
    \end{aligned}$$
Applying Bony's decomposition (\ref{8}) in the horizontal variable to $\theta\partial_y u$ gives
$$\begin{aligned}
    I^k_2=&(e^\Psi\Delta_k[\theta\partial_y u]_\Phi|e^\Psi\Delta_k \partial_y u_\Phi)\\
    =&(e^\Psi\Delta_k[T_{\theta}\partial_y u]_\Phi|e^\Psi\Delta_k \partial_y u_\Phi)+(e^\Psi\Delta_k[T_{\partial_y u}\theta]_\Phi|e^\Psi\Delta_k \partial_y u_\Phi)\\
    &+(e^\Psi\Delta_k[R(\theta,\partial_y u)]_\Phi|e^\Psi\Delta_k \partial_y u_\Phi)\\
    =:& I^k_{2,1}+I^k_{2,2}+I^k_{2,3}
\end{aligned}$$
Noting (\ref{9}) and using the support properties to the Fourier transform of the terms in $T_\theta \partial_y u$, we get
$$|I^k_{2,1}|\lesssim\int_0^t\int_0^\infty e^{2\Psi}\sum\limits_{|k'-k|\leq4}\sum\limits_{k''\leq k'-2}2^{\frac{k''}{2}}\|\Delta_{k''}\theta_\Phi\|_{L_x^2}\|\Delta_{k'}\partial_y u_\Phi\|_{L_x^2}\|\Delta_k \partial_y u_\Phi\|_{L_x^2}dydt'.$$

Due to that $\lim\limits_{y\to\infty}\theta=0$, by using (\ref{51}) and Hölder's inequality, it follows
$$\begin{aligned}
    |I^k_{2,1}|&\lesssim\int_0^t\sum\limits_{|k'-k|\leq4}\langle t\rangle^{\frac{1}{4}}
\|\partial_y\theta_\Phi\|_{B^{\frac{1}{2},0}}
\|e^\Psi\Delta_{k'}\partial_y u_\Phi\|_{L_+^2}\|e^\Psi\Delta_k\partial_y u_\Phi\|_{L_+^2}dt'\\
    &\leq\epsilon\langle t\rangle^{\frac{1}{4}}\sum\limits_{|k'-k|\leq4}\left(\int_0^t\|e^\Psi\Delta_k\partial_y u_\Phi\|_{L_+^2}^2dt'\right)^{\frac{1}{2}}\left(\int_0^t\|e^\Psi\Delta_{k'}\partial_y u_\Phi\|_{L_+^2}^2dt'\right)^{\frac{1}{2}},
\end{aligned}$$
which implies 
$$\sum\limits_{k\in\mathbb{Z}}2^k\sqrt{|I^k_{2,1}|}\lesssim\epsilon^{\frac{1}{2}}\langle t\rangle^{\frac{1}{8}}\|\partial_yu_\Phi\|_{\widetilde{L}_t^2(B^{1,0})}
$$
by using Young's inequality.
Along the same line, one can deduce that
$$\begin{aligned}
    |I_{2,2}^k|&\lesssim\int_0^t\int_0^\infty e^{2\Psi}\sum\limits_{|k'-k|\leq4}\sum\limits_{k''\leq k'-2}2^{\frac{k''}{2}}\|\Delta_{k''}\partial_y u_\Phi\|_{L_x^2}\|\Delta_{k'}\theta_\Phi\|_{L_x^2}\|\Delta_k \partial_y u_\Phi\|_{L_x^2}dydt'\\
    &\leq\int_0^t\sum\limits_{|k'-k|\leq 4}\sum\limits_{k''\leq k'-2}\langle t\rangle^{\frac{1}{4}}2^{\frac{k''}{2}}\|e^\Psi\Delta_{k'}\partial_y \theta_\Phi\|_{L_+^2}\|e^\Psi\Delta_{k''}\partial_y u_\Phi\|_{L_+^2}\|e^\Psi\Delta_k\partial_y u_\Phi\|_{L_+^2}dt'\\
    &\leq\int_0^t\sum\limits_{k''\leq k+2}\sum\limits_{|k'-k|\leq 4}\langle t\rangle^{\frac{1}{4}}2^{\frac{k''}{2}}\|e^\Psi\Delta_{k'}\partial_y \theta_\Phi\|_{L_+^2}\|e^\Psi\Delta_{k''}\partial_y u_\Phi\|_{L_+^2}\|e^\Psi\Delta_k\partial_y u_\Phi\|_{L_+^2}dt'\\
    &\lesssim\epsilon\langle t\rangle^{\frac{1}{4}}\sum\limits_{k''\leq k+2}2^{\frac{k''-2k}{2}}d_k\left(\int_0^t\|e^\Psi\Delta_k\partial_y u_\Phi\|_{L_+^2}^2dt'\right)^{\frac{1}{2}}\left(\int_0^t\|e^\Psi\Delta_{k''}\partial_y u_\Phi\|_{L_+^2}^2dt'\right)^{\frac{1}{2}},
\end{aligned}
$$
where 
\begin{equation}
    d_k=\frac{2^k\sum\limits_{|k'-k|\leq 4}\|e^\Psi\Delta_{k'}\partial_y \theta_\Phi\|_{L_t^\infty(L_+^2)}}{\|\partial_y \theta_\Phi\|_{\widetilde{L}_t^\infty(B^{1,0})}\sum\limits_{|j|\leq4}2^j}.
\end{equation}
It follows from Definitions \ref{48}, \eqref{51} that
$$\begin{aligned}
    |I_{2,2}^k|&\lesssim\epsilon\langle t\rangle^{\frac{1}{4}}\sum\limits_{k''\leq k+2}2^{-\frac{4k+k''}{2}}d_k\widetilde{d}_k\widetilde{d}_{k''}\|\partial_yu_\Phi\|_{\widetilde{L}_t^2(B^{1,0})}^2,
\end{aligned}
$$
where 
\begin{equation}
    \widetilde{d}_k=\frac{2^k\left(\int_0^t\|e^\Psi\Delta_k\partial_y u_\Phi\|_{L_+^2}^2dt'\right)^\frac{1}{2}}{\|\partial_y u_\Phi\|_{\widetilde{L}_t^2(B^{1,0})}}.
\end{equation}
Taking square root of the above inequality and multiplying the resulting inequality by $2^k$ and summing over $k\in\mathbb{Z}$, we ﬁnd
$$\begin{aligned}
    \sum\limits_{k\in\mathbb{Z}}2^k\sqrt{|I_{2,2}^k|}&\lesssim\epsilon^{\frac{1}{2}}\langle t\rangle^{\frac{1}{8}}\sum\limits_{k\in\mathbb{Z}}\left(d_k\widetilde{d}_k\sum\limits_{k''\leq k+2}2^{-\frac{k''}{2}}\widetilde{d}_{k''}\right)^\frac{1}{2}\|\partial_yu_\Phi\|_{\widetilde{L}_t^2(B^{1,0})}\\
    &\lesssim\epsilon^{\frac{1}{2}}\langle t\rangle^{\frac{1}{8}}\sum\limits_{k\in\mathbb{Z}}(d_k+\widetilde{d}_k)\|\partial_yu_\Phi\|_{\widetilde{L}_t^2(B^{1,0})}\\
    &\lesssim\epsilon^{\frac{1}{2}}\langle t\rangle^{\frac{1}{8}}\|\partial_yu_\Phi\|_{\widetilde{L}_t^2(B^{1,0})}.
\end{aligned}
$$
From the support properties to the Fourier transform of the terms in $R(\theta,\partial_y u)$, we infer that
$$\begin{aligned}
    |I_{2,3}^k|&\lesssim\int_0^t\int_0^\infty e^{2\Psi}\sum\limits_{k'\geq k-3}\sum\limits_{k''=k'-1}^{k'+1}2^{\frac{k''}{2}}\|\Delta_{k''} \theta_\Phi\|_{L_x^2}\|\Delta_{k'}\partial_y u_\Phi\|_{L_x^2}\|\Delta_k \partial_y u_\Phi\|_{L_x^2}dydt'\\
    &\leq\int_0^t\sum\limits_{k'\geq k-3}\langle t\rangle^{\frac{1}{4}}\|\partial_y\theta_\Phi\|_{B^{\frac{1}{2},0}}\|e^\Psi\Delta_{k'}\partial_y u_\Phi\|_{L_+^2}\|e^\Psi\Delta_k\partial_y u_\Phi\|_{L_+^2}dt'\\
    &\leq\epsilon\langle t\rangle^{\frac{1}{4}}\sum\limits_{k'\geq k-3}\left(\int_0^t\|e^\Psi\Delta_k\partial_y u_\Phi\|_{L_+^2}^2dt'\right)^{\frac{1}{2}}\left(\int_0^t\|e^\Psi\Delta_{k'}\partial_y u_\Phi\|_{L_+^2}^2dt'\right)^{\frac{1}{2}}.
\end{aligned}$$
Applying Young's inequality leads to 
$$\sum\limits_{k\in\mathbb{Z}}2^k\sqrt{|I_{2,3}^k|}\lesssim\epsilon^{\frac{1}{2}}\langle t\rangle^{\frac{1}{8}}\|\partial_yu_\Phi\|_{\widetilde{L}_t^2(B^{1,0})}.
$$
As a result, it turns out that
$$\sum\limits_{k\in\mathbb{Z}}2^k\sqrt{|\mathcal{L}_{3,1}^k|}\gtrsim(\sqrt{\theta^E}-\epsilon^{\frac{1}{2}}\langle t\rangle^{\frac{1}{8}})\|\partial_yu_\Phi\|_{\widetilde{L}_t^2(B^{1,0})}.$$
Thus we complete the proof of Lemma \ref{6}.
\end{proof}

\begin{lemma}\label{10}
For any $ t<T^*$, and $\eta>0$, the term $\mathcal{L}_{3,3}^k$ given in \eqref{3.13} satisfies the following estimate 
    \begin{equation}\label{11}\sum\limits_{k\in\mathbb{Z}}2^k\sqrt{|\mathcal{L}_{3,3}^k|}\lesssim\eta(\|\partial_yu_\Phi\|_{\widetilde{L}_t^2(B^{1,0})}+\|\partial_y\theta_\Phi\|_{\widetilde{L}_t^2(B^{1,0})})+C_\eta\|u_\Phi\|_{\widetilde{L}_{t,\dot{\mu}}^2(B^{\frac{3}{2},0})}.\end{equation}
\end{lemma}
\begin{proof}
Applying Bony's decomposition (\ref{8}) in the horizontal variable to $\partial_y\theta\partial_yu$ gives
    $$\begin{aligned}
    \mathcal{L}_{3,3}^k=&(e^\Psi\Delta_k[\partial_y\theta\partial_y u]_\Phi|e^\Psi\Delta_k u_\Phi)\\
    =&(e^\Psi\Delta_k[T_{\partial_y\theta}\partial_y u]_\Phi|e^\Psi\Delta_k u_\Phi)+(e^\Psi\Delta_k[T_{\partial_y u}\partial_y \theta]_\Phi|e^\Psi\Delta_k u_\Phi)\\
    &+(e^\Psi\Delta_k[R(\partial_y\theta,\partial_y u)]_\Phi|e^\Psi\Delta_k u_\Phi)\\
    =:& I_1^k+I_2^k+I_3^k.
\end{aligned}$$
Considering (\ref{9}) and the support properties to the Fourier transform of the terms in $T_{\partial_y\theta}\partial_y u$, we write
$$|I_1^k|\lesssim\int_0^t\int_0^\infty e^{2\Psi}\sum\limits_{|k'-k|\leq4}\sum\limits_{k''\leq k'-2}2^{\frac{k''}{2}}\|\Delta_{k''}\partial_y \theta_\Phi\|_{L_x^2}\|\Delta_{k'}\partial_y u_\Phi\|_{L_x^2}\|\Delta_k u_\Phi\|_{L_x^2}dydt'.$$
Due to $u|_{y=0}=0$, Definition \ref{48} and Hölder's inequality,
$$\begin{aligned}
    |I_1^k|&\lesssim\int_0^t\sum\limits_{|k'-k|\leq4}\langle t'\rangle^{\frac{1}{4}}\|\partial_y^2\theta_\Phi\|_{B^{\frac{1}{2},0}}\|e^\Psi\Delta_{k'}\partial_y u_\Phi\|_{L_+^2}\|e^\Psi\Delta_ku_\Phi\|_{L_+^2}dt'\\
    &\leq\sum\limits_{|k'-k|\leq4}\left(\int_0^t\dot{\mu}\|e^\Psi\Delta_ku_\Phi\|_{L_+^2}^2dt'\right)^{\frac{1}{2}}\left(\int_0^t\|e^\Psi\Delta_{k'}\partial_y u_\Phi\|_{L_+^2}^2dt'\right)^{\frac{1}{2}}.
\end{aligned}$$
Along the same line, one can deduce that
$$\begin{aligned}
    |I_2^k|&\lesssim\int_0^t\sum\limits_{|k'-k|\leq4}\langle t'\rangle^{\frac{1}{4}}\|\partial_y^2u_\Phi\|_{B^{\frac{1}{2},0}}\|e^\Psi\Delta_{k'}\partial_y \theta_\Phi\|_{L_+^2}\|e^\Psi\Delta_ku_\Phi\|_{L_+^2}dt'\\
    &\leq\sum\limits_{|k'-k|\leq4}\left(\int_0^t\dot{\mu}\|e^\Psi\Delta_ku_\Phi\|_{L_+^2}^2dt'\right)^{\frac{1}{2}}\left(\int_0^t\|e^\Psi\Delta_{k'}\partial_y \theta_\Phi\|_{L_+^2}^2dt'\right)^{\frac{1}{2}}.
\end{aligned}$$
From the support properties to the Fourier transform of the terms in $R(\partial_y\theta,\partial_y u)$, we infer that
$$\begin{aligned}
    |I_{3}^k|&\lesssim\int_0^t\int_0^\infty e^{2\Psi}\sum\limits_{k'\geq k-3}\sum\limits_{k''=k'-1}^{k'+1}2^{\frac{k''}{2}}\|\Delta_{k''}\partial_y \theta_\Phi\|_{L_x^2}\|\Delta_{k'}\partial_y u_\Phi\|_{L_x^2}\|\Delta_ku_\Phi\|_{L_x^2}dydt'\\
    &\leq\int_0^t\sum\limits_{k'\geq k-3}\langle t'\rangle^{\frac{1}{4}}\|\partial_y^2\theta_\Phi\|_{B^{\frac{1}{2},0}}\|e^\Psi\Delta_{k'}\partial_y u_\Phi\|_{L_+^2}\|e^\Psi\Delta_ku_\Phi\|_{L_+^2}dt'\\
    &\leq\sum\limits_{k'\geq k-3}\left(\int_0^t\dot{\mu}\|e^\Psi\Delta_ku_\Phi\|_{L_+^2}^2dt'\right)^{\frac{1}{2}}\left(\int_0^t\|e^\Psi\Delta_{k'}\partial_y u_\Phi\|_{L_+^2}^2dt'\right)^{\frac{1}{2}}.
\end{aligned}$$
  Adding up the above estimates of $I_1^k, I_2^k, I_3^k$, taking square root of the resulting inequality, multiplying by $2^k$ and summing over $k\in\mathbb{Z}$, we ﬁnd that
$$\sum\limits_{k\in\mathbb{Z}}2^k\sqrt{|\mathcal{L}_{3,3}^k|}\lesssim\eta(\|\partial_yu_\Phi\|_{\widetilde{L}_t^2(B^{1,0})}+\|\partial_y\theta_\Phi\|_{\widetilde{L}_t^2(B^{1,0})})+C_\eta\|u_\Phi\|_{\widetilde{L}_{t,\dot{\mu}}^2(B^{\frac{3}{2},0})}.
$$
This finishes the proof of Lemma \ref{10}.
\end{proof}

\begin{lemma}\label{7}
For any $ t<T^*$, and $\eta>0$, the following estimate holds for
the terms $\mathcal{L}_{1}^k$, $\mathcal{L}_{2}^k$ given in \eqref{16} and $\mathcal{L}_{3,2}^k$ given in \eqref{3.13} respectively, 
\begin{equation}\label{52}
\begin{aligned}\sum\limits_{k\in\mathbb{Z}}&2^k\sqrt{|\mathcal{L}_1^k+\mathcal{L}_2^k+\mathcal{L}_{3,2}^k|}\gtrsim\|u_\Phi(t)\|_{B^{1,0}}-\|e^{\delta|D|}u_0\|_{B^{1,0}_{\Psi_0}}+(\sqrt{\lambda}-C_\eta)\|u_\Phi\|_{\widetilde{L}_{t,\dot{\mu}}^2(B^{\frac{3}{2},0})}\\
    &+\|\sqrt{-(\partial_t\Psi+4\theta^E(\partial_y\Psi)^2)}u_\Phi\|_{\widetilde{L}_t^2(B^{1,0})}-(\eta+\frac{\sqrt{\theta^E}}{2})\|\partial_yu_\Phi\|_{\widetilde{L}_t^2(B^{1,0})}-\eta\|\partial_y\theta_\Phi\|_{\widetilde{L}_t^2(B^{1,0})}\end{aligned}.\end{equation}
\end{lemma}

\begin{proof}
First, it is easy to have 
\begin{equation}\label{3.18}
\mathcal{L}_1^k=\frac{1}{2}\|e^\Psi\Delta_ku_\Phi(t)\|_{L_+^2}^2-\frac{1}{2}\|e^{\frac{y^2}{16\theta^E}}\Delta_ke^{\delta|D|}u_0\|_{L_+^2}^2-(\partial_t\Psi e^\Psi\Delta_ku_\Phi|e^\Psi\Delta_ku_\Phi).
\end{equation}
From the support properties of the operator $\Delta_k$,
\begin{equation}\label{3.19}
\mathcal{L}_2^k\gtrsim\lambda2^k(\dot{\mu}e^\Psi\Delta_k u_\Phi|e^\Psi\Delta_k u_\Phi).
\end{equation}
We write
\begin{equation}\label{3.20}
\begin{aligned}
\mathcal{L}_{3,2}^k&=(2e^\Psi\partial_y\Psi\Delta_k[(\theta+\theta^E)\partial_y u]_\Phi|e^\Psi\Delta_k u_\Phi)\\
&=\theta^E(2e^\Psi\partial_y\Psi\Delta_k\partial_y u_\Phi|e^\Psi\Delta_k u_\Phi)+(2e^\Psi\partial_y\Psi\Delta_k[\theta\partial_y u]_\Phi|e^\Psi\Delta_k u_\Phi)\\
&=:I_1^k+I_2^k.
\end{aligned}
\end{equation}
By Young's inequality,
\begin{equation}\label{3.21}
I_1^k\geq-4\theta^E((\partial_y\Psi)^2e^\Psi\Delta_ku_\Phi|e^\Psi\Delta_ku_\Phi)-\frac{1}{4}\theta^E(e^\Psi\Delta_k\partial_y u_\Phi|e^\Psi\Delta_k\partial_y u_\Phi).
\end{equation}
It follows from Lemma \ref{126} and a similar proof of Lemma \ref{10} that
$$\begin{aligned}
    |I_2^k|\lesssim&\int_0^t\int_0^\infty e^{2\Psi}\sum\limits_{|k'-k|\leq4}\sum\limits_{k''\leq k'-2}2^{\frac{k''}{2}}\|\partial_y\Psi\Delta_{k''}\theta_\Phi\|_{L_x^2}\|\Delta_{k'}\partial_y u_\Phi\|_{L_x^2}\|\Delta_k u_\Phi\|_{L_x^2}dydt'\\
    &+\int_0^t\int_0^\infty e^{2\Psi}\sum\limits_{|k'-k|\leq4}\sum\limits_{k''\leq k'-2}2^{\frac{k''}{2}}\|\Delta_{k''}\partial_y u_\Phi\|_{L_x^2}\|\partial_y\Psi\Delta_{k'}\theta_\Phi\|_{L_x^2}\|\Delta_k u_\Phi\|_{L_x^2}dydt'\\
    &+\int_0^t\int_0^\infty e^{2\Psi}\sum\limits_{k'\geq k-3}\sum\limits_{k''=k'-1}^{k'+1}2^{\frac{k''}{2}}\|\Delta_{k''}\partial_y u_\Phi\|_{L_x^2}\|\partial_y\Psi\Delta_{k'}\theta_\Phi\|_{L_x^2}\|\Delta_k u_\Phi\|_{L_x^2}dydt'\\
    \lesssim&\int_0^t\sum\limits_{|k'-k|\leq4}\sum\limits_{k''\leq k'-2}2^{\frac{k''}{2}}\langle t\rangle^{\frac{1}{4}}(\|\partial_y\Psi e^\Psi\Delta_{k''}\partial_y\theta_\Phi\|_{L_+^2}+\langle t\rangle^{-1}\|e^\Psi\Delta_{k''}\theta_\Phi\|_{L_+^2})\\
    &\times\|e^\Psi\Delta_{k'}\partial_y u_\Phi\|_{L_+^2}\|e^\Psi\Delta_k u_\Phi\|_{L_+^2}dt'\\
    &+\int_0^t\sum\limits_{|k'-k|\leq4}\sum\limits_{k''\leq k'-2}2^{\frac{k''}{2}}\langle t\rangle^{\frac{1}{4}}\|e^\Psi\Delta_{k''}\partial_y^2 u_\Phi\|_{L_+^2}\|\partial_y\Psi e^\Psi\Delta_{k'}\theta_\Phi\|_{L_+^2}\|e^\Psi\Delta_k u_\Phi\|_{L_+^2}dt'\\
    &+\int_0^t\sum\limits_{k'\geq k-3}\sum\limits_{k''=k'-1}^{k'+1}2^{\frac{k''}{2}}\langle t\rangle^{\frac{1}{4}}\|e^\Psi\Delta_{k''}\partial_y^2 u_\Phi\|_{L_+^2}\|\partial_y\Psi e^\Psi\Delta_{k'}\theta_\Phi\|_{L_+^2}\|e^\Psi\Delta_k u_\Phi\|_{L_+^2}dt'\\
\end{aligned}$$
$$\begin{aligned}    
    \lesssim&\int_0^t\sum\limits_{|k'-k|\leq4}\sum\limits_{k''\leq k'-2}2^{\frac{k''}{2}}\langle t\rangle^{\frac{1}{4}}(\|e^\Psi\Delta_{k''}\partial_y^2\theta_\Phi\|_{L_+^2}+\|e^\Psi\Delta_{k''}\theta_\Phi\|_{L_+^2})\\
    &\times\|e^\Psi\Delta_{k'}\partial_y u_\Phi\|_{L_+^2}\|e^\Psi\Delta_k u_\Phi\|_{L_+^2}dt'\\
    &+\int_0^t\sum\limits_{|k'-k|\leq4}\sum\limits_{k''\leq k'-2}2^{\frac{k''}{2}}\langle t\rangle^{\frac{1}{4}}\|e^\Psi\Delta_{k''}\partial_y^2 u_\Phi\|_{L_+^2}\|e^\Psi\Delta_{k'}\partial_y\theta_\Phi\|_{L_+^2}\|e^\Psi\Delta_k u_\Phi\|_{L_+^2}dt'\\
    &+\int_0^t\sum\limits_{k'\geq k-3}\sum\limits_{k''=k'-1}^{k'+1}2^{\frac{k''}{2}}\langle t\rangle^{\frac{1}{4}}\|e^\Psi\Delta_{k''}\partial_y^2 u_\Phi\|_{L_+^2}\|e^\Psi\Delta_{k'}\partial_y\theta_\Phi\|_{L_+^2}\|e^\Psi\Delta_k u_\Phi\|_{L_+^2}dt'\\
    \lesssim&\sum\limits_{|k'-k|\leq4}\left(\int_0^t\dot{\mu}\|e^\Psi\Delta_k u_\Phi\|_{L_+^2}^2dt'\right)^{\frac{1}{2}}\left(\int_0^t\|e^\Psi\Delta_{k'}\partial_y u_\Phi\|_{L_+^2}^2dt'\right)^{\frac{1}{2}}\\
    &+\sum\limits_{|k'-k|\leq4}\left(\int_0^t\dot{\mu}\|e^\Psi\Delta_k u_\Phi\|_{L_+^2}^2dt'\right)^{\frac{1}{2}}\left(\int_0^t\|e^\Psi\Delta_{k'}\partial_y\theta_\Phi\|_{L_+^2}^2dt'\right)^{\frac{1}{2}}\\
    &+\sum\limits_{k'\geq k-3}\left(\int_0^t\dot{\mu}\|e^\Psi\Delta_k u_\Phi\|_{L_+^2}^2dt'\right)^{\frac{1}{2}}\left(\int_0^t\|e^\Psi\Delta_{k'}\partial_y\theta_\Phi\|_{L_+^2}^2dt'\right)^{\frac{1}{2}}.
\end{aligned}$$
Taking square root of the above inequality, multiplying by $2^k$ and summing over $k\in\mathbb{Z}$, we have
\begin{equation}\label{3.22}
\sum\limits_{k\in\mathbb{Z}}2^k\sqrt{|I_2^k|}\lesssim\eta(\|\partial_yu_\Phi\|_{\widetilde{L}_t^2(B^{1,0})}+\|\partial_y\theta_\Phi\|_{\widetilde{L}_t^2(B^{1,0})})+C_\eta\|u_\Phi\|_{\widetilde{L}_{t,\dot{\mu}}^2(B^{\frac{3}{2},0})}.
\end{equation}

Plugging estimates \eqref{3.21} and \eqref{3.22} into \eqref{3.20}, and combining with \eqref{3.18} and \eqref{3.19}, 
we obtain the estimate (\ref{52}). This thus finishes the proof of Lemma \ref{7}.
\end{proof}

\begin{lemma}\label{12}
For any $ t<T^*$, the first term $\mathcal{R}_1^k$ given on the right hand side of \eqref{16} satisfies the estimate,    \begin{equation}\label{118}\sum\limits_{k\in\mathbb{Z}}2^k\sqrt{|\mathcal{R}_1^k|}\lesssim\|u_\Phi\|_{\widetilde{L}_{t,\dot{\mu}}^2(B^{\frac{3}{2},0})}.\end{equation}

\end{lemma}
\begin{proof}
This estimate can be obtained in a way similar to that for the estimate (3.9) given in \cite[Section 3]{zhang2016long}. By applying Bony's decomposition to $u\partial_xu$, we have
$$\begin{aligned}
|\mathcal{R}_1^k|\lesssim&\int_0^t\sum\limits_{|k'-k|\leq4}\sum\limits_{k''\leq k'-2}2^{k'+\frac{k''}{2}}\langle t\rangle^{\frac{1}{4}}\|e^\Psi\Delta_{k''}\partial_yu_\Phi\|_{L_+^2}\|e^\Psi\Delta_{k'} u_\Phi\|_{L_+^2}\|e^\Psi\Delta_ku_\Phi\|_{L_+^2}dt'\\
&+\int_0^t\sum\limits_{|k'-k|\leq4}\sum\limits_{k''\leq k'-2}2^{\frac{3}{2}k''}\langle t\rangle^{\frac{1}{4}}\|e^\Psi\Delta_{k''}\partial_y u_\Phi\|_{L_+^2}\|e^\Psi\Delta_{k'} u_\Phi\|_{L_+^2}\|e^\Psi\Delta_k u_\Phi\|_{L_+^2}dt'\\
&+\int_0^t\sum\limits_{k'\geq k-3}2^{k'+\frac{k''}{2}}\langle t\rangle^{\frac{1}{4}}\|e^\Psi\Delta_{k''}\partial_y u_\Phi\|_{L_+^2}\|e^\Psi\Delta_{k'} u_\Phi\|_{L_+^2}\|e^\Psi\Delta_k u_\Phi\|_{L_+^2}dt'\\
\lesssim&\int_0^t\sum\limits_{|k'-k|\leq4}2^{k'}\langle t\rangle^{\frac{1}{4}}\|\partial_yu_\Phi\|_{B^{\frac{1}{2},0}}\|e^\Psi\Delta_{k'} u_\Phi\|_{L_+^2}\|e^\Psi\Delta_ku_\Phi\|_{L_+^2}dt'\\
&+\int_0^t\sum\limits_{k'\geq k-3}2^{k'}\langle t\rangle^{\frac{1}{4}}\|\partial_yu_\Phi\|_{B^{\frac{1}{2},0}}\|e^\Psi\Delta_{k'} u_\Phi\|_{L_+^2}\|e^\Psi\Delta_ku_\Phi\|_{L_+^2}dt'\\
\leq&\sum\limits_{|k'-k|\leq4}2^{k'}\left(\int_0^t\dot{\mu}\|e^\Psi\Delta_ku_\Phi\|_{L_+^2}^2dt'\right)^{\frac{1}{2}}\left(\int_0^t\dot{\mu}\|e^\Psi\Delta_{k'} u_\Phi\|_{L_+^2}^2dt'\right)^{\frac{1}{2}}\\
&+\sum\limits_{|k'-k|\leq4}2^{k'}\left(\int_0^t\dot{\mu}\|e^\Psi\Delta_ku_\Phi\|_{L_+^2}^2dt'\right)^{\frac{1}{2}}\left(\int_0^t\dot{\mu}\|e^\Psi\Delta_{k'} u_\Phi\|_{L_+^2}^2dt'\right)^{\frac{1}{2}}\\
&+\sum\limits_{k'\geq k-3}2^{k'}\left(\int_0^t\dot{\mu}\|e^\Psi\Delta_ku_\Phi\|_{L_+^2}^2dt'\right)^{\frac{1}{2}}\left(\int_0^t\dot{\mu}\|e^\Psi\Delta_{k'} u_\Phi\|_{L_+^2}^2dt'\right)^{\frac{1}{2}}.\\
\end{aligned}$$
Therefore the above computation implies the estimate (\ref{118}) immediately by applying Young's inequality.
\end{proof}

\begin{lemma}\label{13}
For any $ t<T^*$ and $\eta>0$, the second term $\mathcal{R}_2^k$ given on the right hand side of \eqref{16} satisfies the estimate,     \begin{equation}\label{14}\sum\limits_{k\in\mathbb{Z}}2^k\sqrt{|\mathcal{R}_2^k|}\lesssim\eta(\|\partial_yu_\Phi\|_{\widetilde{L}_t^2(B^{1,0})}+\|\partial_y\theta_\Phi\|_{\widetilde{L}_t^2(B^{1,0})})+C_\eta\|u_\Phi\|_{\widetilde{L}_{t,\dot{\mu}}^2(B^{\frac{3}{2},0})}.\end{equation}
\end{lemma}

\begin{proof}
In view of (\ref{1}), we infer that 
\begin{equation}\label{3.25}
\begin{aligned}
\mathcal{R}_2^k=&-(e^\Psi\Delta_k[v\partial_yu]_\Phi|e^\Psi\Delta_ku_\Phi)\\
=&(e^\Psi\Delta_k[\int_0^y\partial_x udy'\partial_yu]_\Phi|e^\Psi\Delta_ku_\Phi)-(e^\Psi\Delta_k[\partial_y\theta\partial_yu]_\Phi|e^\Psi\Delta_ku_\Phi)\\
&-(e^\Psi\Delta_k[\int_0^y(\partial_y u)^2dy'\partial_yu]_\Phi|e^\Psi\Delta_ku_\Phi)\\
=:&\mathcal{R}_{2,1}^k-\mathcal{R}_{2,2}^k-\mathcal{R}_{2,3}^k.
\end{aligned}
\end{equation}
The estimate for the first term $\mathcal{R}_{2,1}^k$
can be derived in a way similar to that for (3.10) given \cite[Section 3]{zhang2016long}. By applying Bony's decomposition to $\int_0^y\partial_x udy'\partial_yu$, we obtain that
$$\begin{aligned}
|\mathcal{R}_{2,1}^k|\lesssim&\int_0^t\sum\limits_{|k'-k|\leq4}2^{\frac{k'}{2}}\langle t\rangle^{\frac{1}{4}}\|u_\Phi\|_{B^{1,0}}\|e^\Psi\Delta_{k'} \partial_yu_\Phi\|_{L_+^2}\|e^\Psi\Delta_ku_\Phi\|_{L_+^2}dt'\\
&+\int_0^t\sum\limits_{|k'-k|\leq4}2^{k'}\langle t\rangle^{\frac{1}{4}}\|\partial_yu_\Phi\|_{B^{\frac{1}{2},0}}\|e^\Psi\Delta_{k'} u_\Phi\|_{L_+^2}\|e^\Psi\Delta_ku_\Phi\|_{L_+^2}dt'\\
&+\int_0^t\sum\limits_{k'\geq k-3}2^{k'}\langle t\rangle^{\frac{1}{4}}\|\partial_yu_\Phi\|_{B^{\frac{1}{2},0}}\|e^\Psi\Delta_{k'} u_\Phi\|_{L_+^2}\|e^\Psi\Delta_ku_\Phi\|_{L_+^2}dt'\\
\leq&\sum\limits_{|k'-k|\leq4}2^{\frac{k'}2}\left(\int_0^t\dot{\mu}\|e^\Psi\Delta_ku_\Phi\|_{L_+^2}^2dt'\right)^{\frac{1}{2}}\left(\int_0^t\|e^\Psi\Delta_{k'} \partial_yu_\Phi\|_{L_+^2}^2dt'\right)^{\frac{1}{2}}\\
&+\sum\limits_{|k'-k|\leq4}2^{k'}\left(\int_0^t\dot{\mu}\|e^\Psi\Delta_ku_\Phi\|_{L_+^2}^2dt'\right)^{\frac{1}{2}}\left(\int_0^t\dot{\mu}\|e^\Psi\Delta_{k'} u_\Phi\|_{L_+^2}^2dt'\right)^{\frac{1}{2}}\\
&+\sum\limits_{k'\geq k-3}2^{k'}\left(\int_0^t\dot{\mu}\|e^\Psi\Delta_ku_\Phi\|_{L_+^2}^2dt'\right)^{\frac{1}{2}}\left(\int_0^t\dot{\mu}\|e^\Psi\Delta_{k'} u_\Phi\|_{L_+^2}^2dt'\right)^{\frac{1}{2}},\\
\end{aligned}$$
which implies 
\begin{equation}\label{3.26}
\sum\limits_{k\in\mathbb{Z}}2^k\sqrt{|\mathcal{R}_{2,1}^k|}\lesssim\eta\|\partial_yu_\Phi\|_{\widetilde{L}_t^2(B^{1,0})}+C_\eta\|u_\Phi\|_{\widetilde{L}_{t,\dot{\mu}}^2(B^{\frac{3}{2},0})},
\end{equation}
in a way similar to that given for (\ref{118}). For the second term  $\mathcal{R}_{2,2}^k$ given in \eqref{3.25},  from Lemma \ref{10} we have that 
\begin{equation}\label{3.27}
\sum\limits_{k\in\mathbb{Z}}2^k\sqrt{|\mathcal{R}_{2,2}^k|}\lesssim\eta(\|\partial_yu_\Phi\|_{\widetilde{L}_t^2(B^{1,0})}+\|\partial_y\theta_\Phi\|_{\widetilde{L}_t^2(B^{1,0})})+C_\eta\|u_\Phi\|_{\widetilde{L}_{t,\dot{\mu}}^2(B^{\frac{3}{2},0})}.
\end{equation}
To estimate $\mathcal{R}_{2,3}^k$ given in \eqref{3.25}, by applying Bony's decomposition (\ref{8}) in the horizontal variable to $\int_0^y(\partial_y u)^2dy'\partial_yu$, we obtain that
\begin{equation}\label{3.28}
\begin{aligned}
    \mathcal{R}_{2,3}^k=&(e^\Psi\Delta_k[T_{\int_0^y(\partial_y u)^2dy'}\partial_yu]_\Phi|e^\Psi\Delta_ku_\Phi)+(e^\Psi\Delta_k[T_{\partial_yu}\int_0^y(\partial_y u)^2dy']_\Phi|e^\Psi\Delta_ku_\Phi)\\
    &+(e^\Psi\Delta_k[R(\int_0^y(\partial_y u)^2dy',\partial_yu)]_\Phi|e^\Psi\Delta_ku_\Phi)\\
    =:&I_1^k+I_2^k+I_3^k.
\end{aligned}
\end{equation}
By applying Bony's decomposition again to $\int_0^y(\partial_y u)^2dy'$, we deduce that
$$\begin{aligned}
    I_1^k=&2(e^\Psi\Delta_k[T_{\int_0^y T_{\partial_y u}\partial_y udy'}\partial_yu]_\Phi|e^\Psi\Delta_ku_\Phi)+(e^\Psi\Delta_k[T_{\int_0^yR(\partial_yu,\partial_y u)dy'}\partial_yu]_\Phi|e^\Psi\Delta_ku_\Phi)\\
    =:&2I_{1,1}^k+I_{1,2}^k.
\end{aligned}$$
Consider the support properties to the Fourier transform of the terms in $T_{\int_0^y T_{\partial_y u}\partial_y udy'}$, we write
$$|I_{1,1}^k|\leq\int_0^t\int_0^\infty e^{2\Psi}\sum\limits_{|k'-k|\leq4}\sum\limits_{k''\leq k'-2}2^{\frac{k''}{2}}\|\int_0^y\Delta_{k''}[T_{\partial_y u}\partial_y u]_\Phi dy'\|_{L_x^2}\|\Delta_{k'}\partial_y u_\Phi\|_{L_x^2}\|\Delta_k u_\Phi\|_{L_x^2}dydt'.$$
Due to Minkowski's integral inequality, the support properties to the Fourier transform of the terms in $T_{\partial_y u}\partial_y u$ and Hölder's inequality,
$$
\begin{aligned}
    \|
    {\int_0^y
    }
    \Delta_{k''}[T_{\partial_y u}\partial_y u]_\Phi dy'\|_{L_x^2}
    {
    \leq
    }
    &\int_0^y\|\Delta_{k''}[T_{\partial_y u}\partial_y u]_\Phi\|_{L_x^2}dy'\\
    \lesssim&\sum\limits_{|l-k''|\leq 4}\sum\limits_{l'\leq l-2}2^{\frac{l'}{2}}\|\Delta_{l'}\partial_y u_\Phi\|_{L_+^2}\|\Delta_{l}\partial_y u_\Phi\|_{L_+^2}\\
    \lesssim&\sum\limits_{|l-k''|\leq 4}\|\partial_yu_\Phi\|_{B^{\frac{1}{2},0}}\|\Delta_l\partial_yu_\Phi\|_{L_+^2}.
\end{aligned}
$$
By applying Hölder's inequality again, one obtains that
\begin{equation}\label{27}
\begin{aligned}
    |I_{1,1}^k|\lesssim&\int_0^t\int_0^\infty e^{2\Psi}\sum\limits_{|k'-k|\leq4}\|\partial_yu_\Phi\|_{B^{\frac{1}{2},0}}^2\|\Delta_{k'}\partial_y u_\Phi\|_{L_x^2}\|\Delta_k u_\Phi\|_{L_x^2}dydt'\\
    \leq&\sum\limits_{|k'-k|\leq4}\left(\int_0^t\dot{\mu}\|e^\Psi\Delta_{k} u_\Phi\|_{L_+^2}^2 dt'\right)^{\frac{1}{2}}\left(\int_0^t\|e^\Psi\Delta_{k'}\partial_y u_\Phi\|_{L_+^2}^2 dt'\right)^{\frac{1}{2}}.
\end{aligned}
\end{equation}
Similarly,
$$|I_{1,2}^k|\leq\int_0^t\int_0^\infty e^{2\Psi}\sum\limits_{|k'-k|\leq4}\sum\limits_{k''\leq k'-2}2^{\frac{k''}{2}}\|\int_0^y\Delta_{k''}[R(\partial_y u,\partial_y u)]_\Phi dy'\|_{L_x^2}\|\Delta_{k'}\partial_y u_\Phi\|_{L_x^2}\|\Delta_k u_\Phi\|_{L_x^2}dydt',$$
where
$$
\begin{aligned}
    \|
    {
    \int_0^y
    }
    \Delta_{k''}[R(\partial_y u,\partial_y u)]_\Phi dy'\|_{L_x^2}
    {
    \leq
    }
    &\int_0^y\|\Delta_{k''}[R(\partial_y u,\partial_y u)]_\Phi\|_{L_x^2}dy'\\
    \lesssim&\sum\limits_{l\geq k''-3}\sum\limits_{l'=l-1}^{l+1}2^{\frac{l'}{2}}\|\Delta_{l'}\partial_y u_\Phi\|_{L_+^2}\|\Delta_{l}\partial_y u_\Phi\|_{L_+^2}\\
    \lesssim&\sum\limits_{l\geq k''-3}\|\partial_yu_\Phi\|_{B^{\frac{1}{2},0}}\|\Delta_l\partial_yu_\Phi\|_{L_+^2}. 
\end{aligned}
$$
Therefore, we get that
\begin{equation}\label{28}
\begin{aligned}
    |I_{1,2}^k|\lesssim&\int_0^t\int_0^\infty e^{2\Psi}\sum\limits_{|k'-k|\leq4}\|\partial_yu_\Phi\|_{B^{\frac{1}{2},0}}^2\|\Delta_{k'}\partial_y u_\Phi\|_{L_x^2}\|\Delta_k u_\Phi\|_{L_x^2}dydt'\\
    \leq&\sum\limits_{|k'-k|\leq4}\left(\int_0^t\dot{\mu}\|e^\Psi\Delta_{k} u_\Phi\|_{L_+^2}^2 dt'\right)^{\frac{1}{2}}\left(\int_0^t\|e^\Psi\Delta_{k'}\partial_y u_\Phi\|_{L_+^2}^2 dt'\right)^{\frac{1}{2}}.
\end{aligned}
\end{equation}
By applying Bony's decomposition to $\int_0^y(\partial_y u)^2dy'$ in $I_2$ of \eqref{3.28}, we get
$$\begin{aligned}
    I_2^k=&2(e^\Psi\Delta_k[T_{\partial_yu}\int_0^y T_{\partial_y u}\partial_y udy']_\Phi|e^\Psi\Delta_ku_\Phi)+(e^\Psi\Delta_k[T_{\partial_yu}\int_0^yR(\partial_yu,\partial_y u)dy']_\Phi|e^\Psi\Delta_ku_\Phi)\\
    =:&2I_{2,1}^k+I_{2,2}^k.
\end{aligned}$$
In the same manner as \eqref{27}, we have
\begin{equation}\label{3.31}
\begin{aligned}
    |I_{2,1}^k|\leq&\int_0^t\int_0^\infty e^{2\Psi}\sum\limits_{|k'-k|\leq4}\sum\limits_{k''\leq k'-2}2^{\frac{k''}{2}}\|\Delta_{k''}\partial_y u_\Phi\|_{L_x^2}\|\int_0^y\Delta_{k'}[T_{\partial_y u}\partial_y u]_\Phi dy'\|_{L_x^2}\|\Delta_k u_\Phi\|_{L_x^2}dydt'\\
    \leq&\int_0^t \sum\limits_{|k'-k|\leq 4}\sum\limits_{|l-k'|\leq4}\|\partial_yu_\Phi\|_{B^{\frac{1}{2},0}}^2\|e^\Psi\Delta_{l}\partial_y u_\Phi\|_{L_+^2}\|e^\Psi\Delta_k u_\Phi\|_{L_+^2}dt'\\
    \lesssim&\sum\limits_{|k'-k|\leq 4}\sum\limits_{|l-k'|\leq4}\left(\int_0^t\dot{\mu}\|e^\Psi\Delta_k u_\Phi\|_{L_+^2}^2dt'\right)^{\frac{1}{2}}\left(\int_0^t\|e^\Psi\Delta_{l}\partial_y u_\Phi\|_{L_+^2}^2dt'\right)^{\frac{1}{2}}.
\end{aligned}
\end{equation}

Similarly to \eqref{28}, we have
\begin{equation}\label{3.32}
\begin{aligned}
    |I_{2,2}^k|\leq&\int_0^t\int_0^\infty e^{2\Psi}\sum\limits_{|k'-k|\leq4}\sum\limits_{k''\leq k'-2}2^{\frac{k''}{2}}\|\Delta_{k''}\partial_y u_\Phi\|_{L_x^2}\|\int_0^y\Delta_{k'}[R(\partial_y u,\partial_y u)]_\Phi dy'\|_{L_x^2}\|\Delta_k u_\Phi\|_{L_x^2}dydt'\\
    \leq&\int_0^t\sum\limits_{|k'-k|\leq4}\sum\limits_{l\geq k'-3}\|\partial_yu_\Phi\|_{B^{\frac{1}{2},0}}^2\|e^\Psi\Delta_{l}\partial_y u_\Phi\|_{L_+^2}\|e^\Psi\Delta_k u_\Phi\|_{L_+^2}dt'\\
    \lesssim&\sum\limits_{|k'-k|\leq4}\sum\limits_{l\geq k'-3}\left(\int_0^t\dot{\mu}\|e^\Psi\Delta_k u_\Phi\|_{L_+^2}^2dt'\right)^{\frac{1}{2}}\left(\int_0^t\|e^\Psi\Delta_{l}\partial_y u_\Phi\|_{L_+^2}^2dt'\right)^{\frac{1}{2}}.
\end{aligned}
\end{equation}

Along the same way, for $I_3^k$ given in \eqref{3.28} we infer that
$$\begin{aligned}
    I_3^k=&2(e^\Psi\Delta_k[R(\int_0^y T_{\partial_y u}\partial_y udy',\partial_yu)]_\Phi|e^\Psi\Delta_ku_\Phi)+(e^\Psi\Delta_k[R(\int_0^yR(\partial_yu,\partial_y u)dy',\partial_yu)]_\Phi|e^\Psi\Delta_ku_\Phi)\\
    =:&2I_{3,1}^k+I_{3,2}^k,
\end{aligned}$$
from which one can obtain that
\begin{equation}\label{3.33}
\begin{aligned}
    |I_{3,1}^k|\leq&\int_0^t\int_0^\infty e^{2\Psi}\sum\limits_{k'\geq k-3}\sum\limits_{k''= k'-1}^{k'+1}2^{\frac{k''}{2}}\|\int_0^y\Delta_{k''}[T_{\partial_y u}\partial_y u]_\Phi dy'\|_{L_x^2}\|\Delta_{k'}\partial_y u_\Phi\|_{L_x^2}\|\Delta_k u_\Phi\|_{L_x^2}dydt'\\
    \lesssim&\int_0^t\int_0^\infty e^{2\Psi}\sum\limits_{k'\geq k-3}\|\partial_yu_\Phi\|_{B^{\frac{1}{2},0}}^2\|\Delta_{k'}\partial_y u_\Phi\|_{L_x^2}\|\Delta_k u_\Phi\|_{L_x^2}dydt'\\
    \leq&\sum\limits_{k'\geq k-3}\left(\int_0^t\dot{\mu}\|e^\Psi\Delta_{k} u_\Phi\|_{L_+^2}^2 dt'\right)^{\frac{1}{2}}\left(\int_0^t\|e^\Psi\Delta_{k'}\partial_y u_\Phi\|_{L_+^2}^2 dt'\right)^{\frac{1}{2}}
\end{aligned}
\end{equation}
and
\begin{equation}\label{3.34}
\begin{aligned}
    |I_{3,2}^k|\leq&\int_0^t\int_0^\infty e^{2\Psi}\sum\limits_{k'\geq k-3}\sum\limits_{k''= k'-1}^{k'+1}2^{\frac{k''}{2}}\|\int_0^y\Delta_{k''}[R(\partial_y u,\partial_y u)]_\Phi dy'\|_{L_x^2}\|\Delta_{k'}\partial_y u_\Phi\|_{L_x^2}\|\Delta_k u_\Phi\|_{L_x^2}dydt'\\
    \lesssim&\int_0^t\int_0^\infty e^{2\Psi}\sum\limits_{k'\geq k-3}\|\partial_yu_\Phi\|_{B^{\frac{1}{2},0}}^2\|\Delta_{k'}\partial_y u_\Phi\|_{L_x^2}\|\Delta_k u_\Phi\|_{L_x^2}dydt'\\
    \leq&\sum\limits_{k'\geq k-3}\left(\int_0^t\dot{\mu}\|e^\Psi\Delta_{k} u_\Phi\|_{L_+^2}^2 dt'\right)^{\frac{1}{2}}\left(\int_0^t\|e^\Psi\Delta_{k'}\partial_y u_\Phi\|_{L_+^2}^2 dt'\right)^{\frac{1}{2}}.
\end{aligned}
\end{equation}

Finally, applying Young's inequality and summing up all estimates \eqref{3.26}-\eqref{3.34}, it follows the conclusion (\ref{14}). 
\end{proof}

\begin{proof}[Proof of Proposition \ref{15}]
From (\ref{16}), we deduce that
$$\sum\limits_{k\in\mathbb{Z}}2^k\left(\sqrt{|\mathcal{L}_1^k+\mathcal{L}_2^k+\mathcal{L}_{3,2}^k|}+\sqrt{|\mathcal{L}_{3,1}^k|}-\sqrt{|\mathcal{L}_{3,3}^k|}\right)\lesssim\sum\limits_{k\in\mathbb{Z}}2^k\left(\sqrt{|\mathcal{R}_1^k|}+\sqrt{|\mathcal{R}_2^k|}\right).$$
Combining estimates given in Lemmas \ref{6}-\ref{13}, we obtain the estimate (\ref{17}) immediately.
\end{proof}

\subsection{Estimate of \texorpdfstring{$\theta$}{theta}}
In view of (\ref{1}), by a similar derivation of (\ref{4}), we get
$$\partial_t \theta_\Phi+\lambda\dot{\mu}|D|\theta_\Phi+[u\partial_x \theta]_\Phi+[v\partial_y \theta]_\Phi=[(\theta+\theta^E)\partial_y^2 \theta]_\Phi+[(\theta+\theta^E)(\partial_y u)^2]_\Phi.$$
By applying the dyadic operator $\Delta_k$ to the above equation and taking the $L_{t,+}^2=L^2(0,t;\R_+^2)$ inner product of the resulting equation with $e^{2\Psi}\Delta_k\theta_\Phi$, we have
\begin{equation}\label{20}
\begin{aligned}
    &(e^\Psi\partial_t \Delta_k\theta_\Phi|e^\Psi\Delta_k\theta_\Phi)+\lambda(\dot{\mu}|D| e^\Psi\Delta_k\theta_\Phi|e^\Psi\Delta_k\theta_\Phi)-(e^\Psi\Delta_k[(\theta+\theta^E)\partial_y^2 \theta]_\Phi|e^\Psi\Delta_k\theta_\Phi)\\
    =&-(e^\Psi\Delta_k[u\partial_x \theta]_\Phi|e^\Psi\Delta_k\theta_\Phi)-(e^\Psi\Delta_k[v\partial_y \theta]_\Phi|e^\Psi\Delta_k\theta_\Phi)+(e^\Psi\Delta_k[(\theta+\theta^E)(\partial_y u)^2]_\Phi|e^\Psi\Delta_k\theta_\Phi).
\end{aligned}\end{equation}
Denote the three terms on the left side by $\mathcal{L}_1^k,\mathcal{L}_2^k,\mathcal{L}_3^k$ and the three terms on the right side by $\mathcal{R}_1^k,\mathcal{R}_2^k,\mathcal{R}_3^k$, respectively.

Due to $\partial_y\theta|_{y=0}=0$, 
$$\begin{aligned}
    \mathcal{L}_3^k=&(e^\Psi\Delta_k[(\theta+\theta^E)\partial_y\theta]_\Phi|e^\Psi\Delta_k \partial_y\theta_\Phi)\\
    &+(2e^\Psi\partial_y\Psi\Delta_k[(\theta+\theta^E)\partial_y\theta]_\Phi|e^\Psi\Delta_k \theta_\Phi)+(e^\Psi\Delta_k[(\partial_y\theta)^2]_\Phi|e^\Psi\Delta_k\theta_\Phi).
\end{aligned}$$

By similar computation as given in the proofs of Lemmas \ref{6}-\ref{13}, we can deduce the following lemmas.

\begin{lemma}\label{22}
    For any $ t<T^*$, and $\eta>0$, it holds that  
    $$\begin{aligned}    \sum\limits_{k\in\mathbb{Z}}&2^k\left(\sqrt{|\mathcal{L}_1^k+\mathcal{L}_2^k+\mathcal{L}_{3}^k|}\right)\gtrsim\|\theta_\Phi(t)\|_{B^{1,0}}-\|e^{\delta|D|}\theta_0\|_{B^{1,0}_{\Psi_0}}+(\sqrt{\lambda}-C_\eta)\|\theta_\Phi\|_{\widetilde{L}_{t,\dot{\mu}}^2(B^{\frac{3}{2},0})}\\
    &+\|\sqrt{-(\partial_t\Psi+4\theta^E(\partial_y\Psi)^2)}\theta_\Phi\|_{\widetilde{L}_t^2(B^{1,0})}+(\frac{\sqrt{\theta^E}}{2}-\epsilon^{\frac{1}{2}}\langle t\rangle^{\frac{1}{8}}-\eta)\|\partial_y\theta_\Phi\|_{\widetilde{L}_t^2(B^{1,0})},
    \end{aligned}$$
and $$\sum\limits_{k\in\mathbb{Z}}2^k\left(\sqrt{|\mathcal{R}_1^k|}+\sqrt{|\mathcal{R}_2^k|}\right)\lesssim\eta(\|\partial_yu_\Phi\|_{\widetilde{L}_t^2(B^{1,0})}+\|\partial_y\theta_\Phi\|_{\widetilde{L}_t^2(B^{1,0})})+C_\eta\|\theta_\Phi\|_{\widetilde{L}_{t,\dot{\mu}}^2(B^{\frac{3}{2},0})}+\|u_\Phi\|_{\widetilde{L}_{t,\dot{\mu}}^2(B^{\frac{3}{2},0})}.$$
\end{lemma}

It remains to study $\mathcal{R}_3^k$ given in \eqref{20}.

\begin{lemma}\label{19}
    For any $ t<T^*$, for any $\eta>0$, it holds that  
    \begin{equation}\label{18}\sum\limits_{k\in\mathbb{Z}}2^k\sqrt{|\mathcal{R}_3^k|}\lesssim\eta\|\partial_yu_\Phi\|_{\widetilde{L}_{t}^2(B^{1,1})}+C_\eta\|\theta_\Phi\|_{\widetilde{L}_{t,\dot{\mu}}^2(B^{\frac{3}{2},0})}.\end{equation}
\end{lemma}

\begin{proof}
We write 
$$\begin{aligned}
\mathcal{R}_3^k=&(e^\Psi\Delta_k[(\theta+\theta^E)(\partial_y u)^2]_\Phi|e^\Psi\Delta_k\theta_\Phi)\\
=&(e^\Psi\Delta_k[\theta(\partial_y u)^2]_\Phi|e^\Psi\Delta_k\theta_\Phi)+\theta^E(e^\Psi\Delta_k[(\partial_y u)^2]_\Phi|e^\Psi\Delta_k\theta_\Phi)\\
=&:I_1^k+I_2^k.
\end{aligned}$$
Due to $\lim\limits_{y\to\infty}\partial_y u=0$, by a similar argument as given in (\ref{27}), (\ref{28}) we have that
\begin{equation}\label{4.3}
\begin{aligned}
    &\|\Delta_{k}[(\partial_y u)^2]_\Phi\|_{L_x^2}=\|\int_y^\infty\Delta_{k}[2\partial_y u\partial_y^2 u]_\Phi dy'\|_{L_x^2}\\
    \lesssim&\sum\limits_{|l-k|\leq 4}(\|\partial_yu_\Phi\|_{B^{\frac{1}{2},0}}\|\Delta_l\partial_y^2u_\Phi\|_{L_+^2}+\|\partial_y^2u_\Phi\|_{B^{\frac{1}{2},0}}\|\Delta_l\partial_yu_\Phi\|_{L_+^2})+\sum\limits_{l\geq k-3}\|\partial_y^2u_\Phi\|_{B^{\frac{1}{2},0}}\|\Delta_l\partial_yu_\Phi\|_{L_+^2}.
\end{aligned}
\end{equation}
By Bony's decomposition it follows
$$\begin{aligned}
    I_1^k=(e^\Psi\Delta_k[T_{(\partial_y u)^2}\theta]_\Phi|e^\Psi\Delta_k\theta_\Phi)+(e^\Psi\Delta_k[T_{\theta}\partial_y^2 u]_\Phi|e^\Psi\Delta_k\theta_\Phi)+(e^\Psi\Delta_k[R((\partial_y u)^2,\theta)]_\Phi|e^\Psi\Delta_k\theta_\Phi).
\end{aligned}$$
From a similar derivation as given in the proof of Lemma \ref{13}, and using \eqref{4.3} we have
$$\begin{aligned}
    |I_1^k|\lesssim&\int_0^t\int_0^\infty e^{2\Psi}\sum\limits_{|k'-k|\leq4}\sum\limits_{k''\leq k'-2}2^{\frac{k''}{2}}\|\Delta_{k''}[(\partial_y u)^2]_\Phi\|_{L_x^2}\|\Delta_{k'}\theta_\Phi\|_{L_x^2}\|\Delta_k \theta_\Phi\|_{L_x^2}dydt'\\
    &+\int_0^t\int_0^\infty e^{2\Psi}\sum\limits_{|k'-k|\leq4}\sum\limits_{k''\leq k'-2}2^{\frac{k''}{2}}\|\Delta_{k''}\theta_\Phi\|_{L_x^2}\|\Delta_{k'}[(\partial_y u)^2]_\Phi\|_{L_x^2}\|\Delta_k \theta_\Phi\|_{L_x^2}dydt'\\
    &+\int_0^t\int_0^\infty e^{2\Psi}\sum\limits_{k'\geq k-3}\sum\limits_{k''=k'-1}^{k'+1}2^{\frac{k''}{2}}\|\Delta_{k''}[(\partial_y u)^2]_\Phi\|_{L_x^2}\|\Delta_{k'}\theta_\Phi\|_{L_x^2}\|\Delta_k \theta_\Phi\|_{L_x^2}dydt'\\
    \lesssim&\sum\limits_{|k'-k|\leq4}\left(\int_0^t\dot{\mu}\|e^\Psi\Delta_k\theta_\Phi\|_{L_+^2}^2 dt'\right)^{\frac{1}{2}}\left(\int_0^t\dot{\mu}\|e^\Psi\Delta_{k'}\theta_\Phi\|_{L_+^2}^2 dt'\right)^{\frac{1}{2}}\\
    &+\sum\limits_{|k'-k|\leq4}\sum\limits_{|l-k'|\leq4}\left(\int_0^t\dot{\mu}\|e^\Psi\Delta_k\theta_\Phi\|_{L_+^2}^2 dt'\right)^{\frac{1}{2}}\left(\int_0^t\|e^\Psi\Delta_l\partial_yu_\Phi\|_{L_+^2}^2 dt'\right)^{\frac{1}{2}}\\
    &+\sum\limits_{|k'-k|\leq4}\sum\limits_{|l-k'|\leq4}\left(\int_0^t\dot{\mu}\|e^\Psi\Delta_k\theta_\Phi\|_{L_+^2}^2 dt'\right)^{\frac{1}{2}}\left(\int_0^t\|e^\Psi\Delta_l\partial_y^2u_\Phi\|_{L_+^2}^2 dt'\right)^{\frac{1}{2}}\\
    &+\sum\limits_{|k'-k|\leq4}\sum\limits_{l\geq k'-3}\left(\int_0^t\dot{\mu}\|e^\Psi\Delta_k\theta_\Phi\|_{L_+^2}^2 dt'\right)^{\frac{1}{2}}\left(\int_0^t\|e^\Psi\Delta_l\partial_yu_\Phi\|_{L_+^2}^2 dt'\right)^{\frac{1}{2}}\\
    &+\sum\limits_{k'\geq k-3}\left(\int_0^t\dot{\mu}\|e^\Psi\Delta_k\theta_\Phi\|_{L_+^2}^2 dt'\right)^{\frac{1}{2}}\left(\int_0^t\dot{\mu}\|e^\Psi\Delta_{k'}\theta_\Phi\|_{L_+^2}^2 dt'\right)^{\frac{1}{2}}.
\end{aligned}$$
By applying Bony's decomposition to $(\partial_yu)^2$, we get that
$$\begin{aligned}|I_2^k|\lesssim&2\int_0^t\int_0^\infty e^{2\Psi}\|\Delta_k[T_{\partial_y u}\partial_y u]_\Phi\|_{L_x^2}\|\Delta_k\theta_\Phi\|_{L_x^2}dydt'\\
&+\int_0^t\int_0^\infty e^{2\Psi}\|\Delta_k[R(\partial_y u,\partial_y u)]_\Phi\|_{L_x^2}\|\Delta_k\theta_\Phi\|_{L_x^2}dydt'\\
=&:2I_{2,1}^k+I_{2,2}^k.
\end{aligned}$$
Considering the Fourier support properties of $T_{\partial_y u}\partial_y u$, we have
$$\begin{aligned}|I_{2,1}^k|\lesssim&\int_0^t\int_0^\infty e^{2\Psi}\sum\limits_{|k'-k|\leq4}\sum\limits_{k''\leq k'-2}2^{\frac{k''}{2}}\|\Delta_{k''}\partial_y u_\Phi\|_{L_x^2}\|\Delta_{k'}\partial_y u_\Phi\|_{L_x^2}\|\Delta_{k}\theta_\Phi\|_{L_x^2}dydt'\\
\leq&\int_0^t\int_0^\infty \sum\limits_{|k'-k|\leq4}\sum\limits_{k''\leq k'-2}2^{\frac{k''}{2}}\langle t\rangle^{\frac{1}{4}}\|e^\Psi\Delta_{k''}\partial_y^2 u_\Phi\|_{L_+^2}\|e^\Psi\Delta_{k'}\partial_y u_\Phi\|_{L_x^2}\|e^\Psi\Delta_k\theta_\Phi\|_{L_x^2}dydt'\\
\leq&\sum\limits_{|k'-k|\leq4}\left(\int_0^t\|e^\Psi\Delta_{k'}\partial_y u_\Phi\|^2_{L_+^2}dt'\right)^\frac{1}{2}\left(\int_0^t\dot{\mu}\|e^\Psi\Delta_{k}\theta_\Phi\|_{L_+^2}^2 dt'\right)^\frac{1}{2}.\\
\end{aligned}$$
Similarly, one has
$$\begin{aligned}|I_{2,2}^k|\lesssim&\int_0^t\int_0^\infty e^{2\Psi}\sum\limits_{k'\geq k-3}\sum\limits_{k''=k'-1}^{k'+1}2^{\frac{k''}{2}}\|\Delta_{k''}\partial_y u_\Phi\|_{L_x^2}\|\Delta_{k'}\partial_y u_\Phi\|_{L_x^2}\|\Delta_{k}\theta_\Phi\|_{L_x^2}dydt'\\
\leq&\int_0^t\int_0^\infty \sum\limits_{k'\geq k-3}\sum\limits_{k''=k'-1}^{k'+1}2^{\frac{k''}{2}}\langle t\rangle^{\frac{1}{4}}\|e^\Psi\Delta_{k''}\partial_y^2 u_\Phi\|_{L_+^2}\|e^\Psi\Delta_{k'}\partial_y u_\Phi\|_{L_x^2}\|e^\Psi\Delta_k\theta_\Phi\|_{L_x^2}dydt'\\
\leq&\sum\limits_{k'\geq k-3}\left(\int_0^t\|e^\Psi\Delta_{k'}\partial_y u_\Phi\|^2_{L_+^2}dt'\right)^\frac{1}{2}\left(\int_0^t\dot{\mu}\|e^\Psi\Delta_{k}\theta_\Phi\|_{L_+^2}^2 dt'\right)^\frac{1}{2}.\\
\end{aligned}$$
Taking square root of the above estimates of $I_1^k, I_2^k$ and multiplying the resulting inequalities by $2^k$ and
summing over $k\in\mathbb{Z}$ lead to the estimate (\ref{18}). This completes the proof of Lemma \ref{19}.
\end{proof}

\begin{proof}[Proof of Proposition \ref{21}]
Due to (\ref{20}), we infer that
$$\sum\limits_{k\in\mathbb{Z}}2^k\left(\sqrt{|\mathcal{L}_1^k+\mathcal{L}_2^k+\mathcal{L}_{3}^k|}\right)\lesssim\sum\limits_{k\in\mathbb{Z}}2^k\left(\sqrt{|\mathcal{R}_1^k|}+\sqrt{|\mathcal{R}_2^k|}+\sqrt{|\mathcal{R}_3^k|}\right).$$
By combining Lemmas \ref{22}-\ref{19}, we conclude (\ref{23}). We thus finish the proof of Proposition \ref{21}.
\end{proof}

\subsection{Estimate of \texorpdfstring{$\partial_yu$}{partial y u}}
From (\ref{1}), we know that
$\bar{u}=\partial_yu$ satisfies
$$\partial_t\bar{u}+u\partial_x\bar{u}+\bar{u}\partial_y^2\theta+\bar{u}^3-\int_0^y\partial_xudy'\partial_y\bar{u}+\int_0^y\bar{u}^2dy'\partial_y\bar{u}=(\theta+\theta^E)\partial_y^2\bar{u},$$
which implies that ${\bar u}_\Phi=
\mathcal{F}^{-1}_{\xi\rightarrow x}(e^{\Phi(t,\xi)}\hat{\bar u}(t,\xi,y))$ satisfies the following equation,
$$\partial_t\bar{u}_\Phi+\lambda\dot{\mu}|D|\bar{u}_\Phi+[u\partial_x\bar{u}]_\Phi+[\bar{u}\partial_y^2\theta]_\Phi+[\bar{u}^3]_\Phi-[\int_0^y\partial_xudy'\partial_y\bar{u}]_\Phi+[\int_0^y\bar{u}^2dy'\partial_y\bar{u}]_\Phi=[(\theta+\theta^E)\partial_y^2\bar{u}]_\Phi.$$
Applying the dyadic operator $\Delta_k$ to the above equation and taking the $L_{t,+}^2=L^2(0,t;\R_+^2)$ inner product of the resulting equation with $e^{2\Psi}\Delta_k{\bar u}_\Phi$, we get 
\begin{equation}\label{32}
\begin{aligned}
    &(e^\Psi\partial_t\Delta_k\bar{u}_\Phi|e^\Psi\Delta_k\bar{u}_\Phi)
    +\lambda(\dot{\mu}|D| e^\Psi\Delta_k\bar{u}_\Phi|e^\Psi\Delta_k\bar{u}_\Phi)-(e^\Psi\Delta_k[(\theta+\theta^E)\partial_y^2\bar{u}]_\Phi|e^\Psi\Delta_k\bar{u}_\Phi)\\
    =&-(e^\Psi\Delta_k[u\partial_x\bar{u}]_\Phi|e^\Psi\Delta_k\bar{u}_\Phi)+(e^\Psi\Delta_k[\int_0^y\partial_xudy'\partial_y\bar{u}]_\Phi
    |e^\Psi\Delta_k\bar{u}_\Phi)\\
    &-(e^\Psi\Delta_k[\int_0^y\bar{u}^2dy'\partial_y\bar{u}]_\Phi|e^\Psi\Delta_k\bar{u}_\Phi)-(e^\Psi\Delta_k[\bar{u}\partial_y^2\theta]_\Phi|e^\Psi\Delta_k\bar{u}_\Phi)-(e^\Psi\Delta_k[\bar{u}^3]_\Phi|e^\Psi\Delta_k\bar{u}_\Phi).
\end{aligned}
\end{equation} 
Denote the three terms on the left side by $\mathcal{L}_1^k,\mathcal{L}_2^k,\mathcal{L}_3^k$ and the five terms on the right side by $\mathcal{R}_1^k,\cdots,\mathcal{R}_5^k$ respectively.

By the fact that $\partial_y\bar{u}|_{y=0}=0$ from (\ref{1}), 
$$\begin{aligned}
    \mathcal{L}_3^k=&(e^\Psi\Delta_k[(\theta+\theta^E)\partial_y \bar{u}]_\Phi|e^\Psi\Delta_k \partial_y \bar{u}_\Phi)\\
    &+(2e^\Psi\partial_y\Psi\Delta_k[(\theta+\theta^E)\partial_y \bar{u}]_\Phi|e^\Psi\Delta_k \bar{u}_\Phi)+(e^\Psi\Delta_k[\partial_y\theta\partial_y \bar{u}]_\Phi|e^\Psi\Delta_k \bar{u}_\Phi)\\
    =&:\mathcal{L}_{3,1}^k+\mathcal{L}_{3,2}^k+\mathcal{L}_{3,3}^k.
\end{aligned}$$

In a way similar to proofs of Lemmas \ref{6}-\ref{13}, we deduce the following results.

\begin{lemma}\label{24}
    For any $ t<T^*$ and $\eta>0$, it holds that  
    $$\begin{aligned}    \sum\limits_{k\in\mathbb{Z}}&2^k\left(\sqrt{|\mathcal{L}_1^k+\mathcal{L}_2^k+\mathcal{L}_{3}^k|}\right)\gtrsim\|\bar{u}_\Phi(t)\|_{B^{1,0}}-\|e^{\delta|D|}\partial_y{u}_0\|_{B^{1,0}_{\Psi_0}}+(\sqrt{\lambda}-C_\eta)\|\bar{u}_\Phi\|_{\widetilde{L}_{t,\dot{\mu}}^2(B^{\frac{3}{2},0})}\\
    &+\|\sqrt{-(\partial_t\Psi+4\theta^E(\partial_y\Psi)^2)}\bar{u}_\Phi\|_{\widetilde{L}_t^2(B^{1,0})}+(\frac{\sqrt{\theta^E}}{2}-\epsilon^{\frac{1}{2}}\langle t\rangle^{\frac{1}{8}}-\eta)\|\partial_y\bar{u}_\Phi\|_{\widetilde{L}_t^2(B^{1,0})}-\eta\|\partial_y\theta_\Phi\|_{\widetilde{L}_t^2(B^{1,1})}.
    \end{aligned}$$
\end{lemma}
\begin{proof}
    By applying Bony's decomposition to $\partial_y\theta\partial_y\bar{u}$ and noting $\partial_y\theta|_{y=0}=0$, we have
    $$\begin{aligned}
        |\mathcal{L}_{3,3}^k|\lesssim&\int_0^t\sum\limits_{|k'-k|\leq4}\langle t\rangle^{\frac{1}{4}}\|\partial_y^2\theta_\Phi\|_{B^{\frac{1}{2},0}}\|e^\Psi\Delta_{k'}\partial_y\bar{u}_\Phi\|_{L_+^2}\|e^\Psi\Delta_k\bar{    u}_\Phi\|_{L_+^2}dt'\\
        &+\int_0^t\sum\limits_{|k'-k|\leq4}\langle t\rangle^{\frac{1}{4}}\|\partial_y^2u_\Phi\|_{B^{\frac{1}{2},0}}\|e^\Psi\Delta_{k'}\partial_y^2 \theta_\Phi\|_{L_+^2}\|e^\Psi\Delta_k\bar{u}_\Phi\|_{L_+^2}dt'\\
        &+\int_0^t\sum\limits_{k'\geq k-3}\langle t\rangle^{\frac{1}{4}}\|\partial_y^2\theta_\Phi\|_{B^{\frac{1}{2},0}}\|e^\Psi\Delta_{k'}\partial_y \bar{u}_\Phi\|_{L_+^2}\|e^\Psi\Delta_k\bar{u}_\Phi\|_{L_+^2}dt'\\
    \end{aligned}$$
    $$\begin{aligned}
        \leq&\sum\limits_{|k'-k|\leq4}\left(\int_0^t\dot{\mu}\|e^\Psi\Delta_k\bar{u}_\Phi\|_{L_+^2}^2dt'\right)^{\frac{1}{2}}\left(\int_0^t\|e^\Psi\Delta_{k'}\partial_y \bar{u}_\Phi\|_{L_+^2}^2dt'\right)^{\frac{1}{2}}\\
        &+\sum\limits_{|k'-k|\leq4}\left(\int_0^t\dot{\mu}\|e^\Psi\Delta_k\bar{u}_\Phi\|_{L_+^2}^2dt'\right)^{\frac{1}{2}}\left(\int_0^t\|e^\Psi\Delta_{k'}\partial_y \bar{\theta}_\Phi\|_{L_+^2}^2dt'\right)^{\frac{1}{2}}\\
        &+\sum\limits_{|k'-k|\leq4}\left(\int_0^t\dot{\mu}\|e^\Psi\Delta_k\bar{u}_\Phi\|_{L_+^2}^2dt'\right)^{\frac{1}{2}}\left(\int_0^t\|e^\Psi\Delta_{k'}\partial_y \bar{u}_\Phi\|_{L_+^2}^2dt'\right)^{\frac{1}{2}},
    \end{aligned}$$
which implies $$\sum\limits_{k\in\mathbb{Z}}2^k\sqrt{|\mathcal{L}_{3,3}^k|}\lesssim\eta(\|\partial_y\bar{u}_\Phi\|_{\widetilde{L}_t^2(B^{1,0})}+\|\partial_y\bar{\theta}_\Phi\|_{\widetilde{L}_t^2(B^{1,0})})+C_\eta\|\bar{u}_\Phi\|_{\widetilde{L}_{t,\dot{\mu}}^2(B^{\frac{3}{2},0})},$$
with $\bar \theta$ denoting $\partial_y\theta$.

The estimate of the remaining terms on the left side of \eqref{32} can be derived along the same way as given in Lemmas \ref{6}-\ref{7}, and then the estimate claimed in this lemma follows.
\end{proof}

\begin{lemma}
    For any $ t<T^*$ and $\eta>0$, the terms $\{{\cal R}_j^k\}_{1\le j\le 4}$ given on the right hand side of \eqref{32} satisfy the following estimate  $$\sum\limits_{j=1}^4\sum\limits_{k\in\mathbb{Z}}2^k\sqrt{|\mathcal{R}_j^k|}\lesssim\eta(\|\partial_y\bar{u}_\Phi\|_{\widetilde{L}_t^2(B^{1,0})}+\|\partial_y\bar{\theta}_\Phi\|_{\widetilde{L}_t^2(B^{1,0})})+C_\eta\|\bar{u}_\Phi\|_{\widetilde{L}_{t,\dot{\mu}}^2(B^{\frac{3}{2},0})}+\|u_\Phi\|_{\widetilde{L}_{t,\dot{\mu}}^2(B^{\frac{3}{2},0})}.$$
\end{lemma}
One can estimate $\{{\cal R}_j^k\}_{1\le j\le 4}$ in a way similar to that given in the  proofs of Lemmas \ref{12} and \ref{13} to conclude this lemma, so wo omit the detail. It remains to estimate $\mathcal{R}_5^k$ given in \eqref{32} .

\begin{lemma}\label{30}
For any $ t<T^*$ and $\eta>0$, it holds that
\begin{equation}\label{29}
\sum\limits_{k\in\mathbb{Z}}2^k\sqrt{|\mathcal{R}_5^k|}\lesssim\eta\|\partial_y\bar{u}_\Phi\|_{\widetilde{L}_{t,\dot{\mu}}^2(B^{1,0})}+C_\eta\|\bar{u}_\Phi\|_{\widetilde{L}_{t,\dot{\mu}}^2(B^{\frac{3}{2},0})}.
\end{equation}
\end{lemma}

\begin{proof}
We recall from Theorem 1.1 of Chapter 5 in \cite{chenbook} that 
$$\begin{aligned}
\bar{u}^3=&T_{\bar{u}^2}\bar{u}+T_{\bar{u}}\bar{u}^2+R(\bar{u},\bar{u}^2)\\
=&3T_{\bar{u}^2}\bar{u}+2(T_{\bar{u}}T_{\bar{u}}-T_{\bar{u}^2})\bar{u}+T_{\bar{u}}R(\bar{u},\bar{u})+R(\bar{u},\bar{u}^2).
\end{aligned}$$
Therefore we can decompose $\mathcal{R}_5^k$ into 
    $$\begin{aligned}
        \mathcal{R}_5^k=&-(e^\Psi\Delta_k[\bar{u}^3]_\Phi|e^\Psi\Delta_k\bar{u}_\Phi)\\
        =&-3(e^\Psi\Delta_k[T_{\bar{u}^2}\bar{u}]_\Phi|e^\Psi\Delta_k\bar{u}_\Phi)-2(e^\Psi\Delta_k[(T_{\bar{u}}T_{\bar{u}}-T_{\bar{u}^2})\bar{u}]_\Phi|e^\Psi\Delta_k\bar{u}_\Phi)\\
        &-(e^\Psi\Delta_k[T_{\bar{u}}R(\bar{u},\bar{u})]_\Phi|e^\Psi\Delta_k\bar{u}_\Phi)-(e^\Psi\Delta_k[R(\bar{u},\bar{u}^2)]_\Phi|e^\Psi\Delta_k\bar{u}_\Phi)\\
    \end{aligned}$$
Due to $\lim\limits_{y\to\infty}\bar{u}=0$, it follows from a similar proof of (\ref{27}), (\ref{28}) that
$$\begin{aligned}
    &\|\Delta_{k}[\bar{u}^2]_\Phi\|_{L_x^2}
=\|\int_y^\infty\Delta_{k}[2\bar{u}\partial_y\bar{u}]_\Phi dy'\|_{L_x^2}\\
    & \qquad\lesssim\sum\limits_{|l-k|\leq 4}(\|\bar{u}_\Phi\|_{B^{\frac{1}{2},0}}\|\Delta_l\partial_y\bar{u}_\Phi\|_{L_+^2}+\|\partial_y\bar{u}_\Phi\|_{B^{\frac{1}{2},0}}\|\Delta_l\bar{u}_\Phi\|_{L_+^2})
    +\sum\limits_{l\geq k-3}\|\partial_y\bar{u}_\Phi\|_{B^{\frac{1}{2},0}}\|\Delta_l\bar{u}_\Phi\|_{L_+^2}.
\end{aligned}$$
Then, by a similar derivation as to that of Lemma \ref{13}, one can have
$$\sum\limits_{k\in\mathbb{Z}}2^k\sqrt{|\mathcal{R}_5|}\lesssim\eta\|\partial_y\bar{u}_\Phi\|_{\widetilde{L}_{t,\dot{\mu}}^2(B^{1,0})}+C_\eta\|\bar{u}_\Phi\|_{\widetilde{L}_{t,\dot{\mu}}^2(B^{\frac{3}{2},0})}.$$ 

\end{proof}

\begin{proof}[Proof of Proposition \ref{31}]
We deduce from (\ref{32}) that 
$$\sum\limits_{k\in\mathbb{Z}}2^k\left(\sqrt{|\mathcal{L}_1^k+\mathcal{L}_2^k+\mathcal{L}_{3}^k|}\right)\lesssim\sum\limits_{j=1}^5\sum\limits_{k\in\mathbb{Z}}2^k\sqrt{|\mathcal{R}_j^k|}.$$
Combining results given in Lemmas \ref{24}-\ref{30}, the estimate (\ref{33}) follows immediately. We thus complete the proof of Proposition \ref{31}.
\end{proof}

\subsection{Estimate of \texorpdfstring{$\partial_y\theta$}{partial y theta}}
First, from (\ref{1}) we see that ${\bar\theta}=\partial_y\theta$ satisfies the following equation,
\begin{equation}\label{85}
\begin{aligned}
\partial_t\bar{\theta}+\partial_yu\partial_x\theta+u\partial_x\bar{\theta}-\bar{\theta}\partial_xu+\bar{\theta}\partial_y\bar{\theta}-\int_0^y\partial_xudy'\partial_y\bar{\theta}+\int_0^y(\partial_yu)^2dy'\partial_y\bar{\theta}&\\
=(\theta+\theta^E)\partial_y^2\bar{\theta}+2(\theta+\theta^E)\partial_yu\partial_y^2u&
\end{aligned}\end{equation}
which implies that 
 ${\bar \theta}_\Phi=
\mathcal{F}^{-1}_{\xi\rightarrow x}(e^{\Phi(t,\xi)}\hat{\bar \theta}(t,\xi,y))$ satisfies 
$$\begin{aligned}
\partial_t\bar{\theta}_\Phi+&\lambda\dot{\mu}|D|\bar{\theta}_\Phi+[\partial_yu\partial_x\theta]_\Phi+[u\partial_x\bar{\theta}]_\Phi-[\bar{\theta}\partial_xu]_\Phi+[\bar{\theta}\partial_y\bar{\theta}]_\Phi\\
-&[\int_0^y\partial_xudy'\partial_y\bar{\theta}]_\Phi+[\int_0^y(\partial_yu)^2dy'\partial_y\bar{\theta}]_\Phi=[(\theta+\theta^E)\partial_y^2\bar{\theta}]_\Phi+2[(\theta+\theta^E)\partial_yu\partial_y^2u]_\Phi.
\end{aligned}$$

By applying the dyadic operator $\Delta_k$ to the above equation and then take the $L_{t,+}^2$ inner product of the resulting equation with $e^{2\Psi}\Delta_k\bar{\theta}_\Phi$, it follows
\begin{equation}\label{39}
\begin{aligned}
&(e^\Psi\partial_t\Delta_k\bar{\theta}_\Phi|e^\Psi\Delta_k\bar{\theta}_\Phi)+\lambda(\dot{\mu}|D|e^\Psi\Delta_k\bar{\theta}_\Phi|e^\Psi\Delta_k\bar{\theta}_\Phi)-(e^\Psi[(\theta+\theta^E)\partial_y^2\bar{\theta}]_\Phi|e^\Psi\Delta_k\bar{\theta}_\Phi)\\
=&-(e^\Psi[u\partial_x\bar{\theta}]_\Phi|e^\Psi\Delta_k\bar{\theta}_\Phi)+(e^\Psi[\int_0^y\partial_xudy'\partial_y\bar{\theta}]_\Phi|e^\Psi\Delta_k\bar{\theta}_\Phi)-(e^\Psi[\int_0^y(\partial_yu)^2dy'\partial_y\bar{\theta}]_\Phi|e^\Psi\Delta_k\bar{\theta}_\Phi)\\
&+(e^\Psi[\bar{\theta}\partial_xu]_\Phi|e^\Psi\Delta_k\bar{\theta}_\Phi)-(e^\Psi[\bar{\theta}\partial_y\bar{\theta}]_\Phi|e^\Psi\Delta_k\bar{\theta}_\Phi)-(e^\Psi[\partial_yu\partial_x\theta]_\Phi|e^\Psi\Delta_k\bar{\theta}_\Phi)\\
&+2(e^\Psi[(\theta+\theta^E)\partial_yu\partial_y^2u]_\Phi|e^\Psi\Delta_k\bar{\theta}_\Phi).
\end{aligned}
\end{equation}
Denote the three terms on the left side by $\mathcal{L}_1^k,\mathcal{L}_2^k,\mathcal{L}_3^k$ and the seven terms on the right side by $\mathcal{R}_1^k,\cdots,\mathcal{R}_7^k$.

Due to $\bar{\theta}|_{y=0}=0$, by using integration by parts, it follows
$$\begin{aligned}
    \mathcal{L}_3^k=&(e^\Psi\Delta_k[(\theta+\theta^E)\partial_y \bar{\theta}]_\Phi|e^\Psi\Delta_k \partial_y \bar{\theta}_\Phi)\\
    &+(2e^\Psi\partial_y\Psi\Delta_k[(\theta+\theta^E)\partial_y \bar{\theta}]_\Phi|e^\Psi\Delta_k \bar{\theta}_\Phi)+(e^\Psi\Delta_k[\partial_y\theta\partial_y \bar{\theta}]_\Phi|e^\Psi\Delta_k \bar{\theta}_\Phi).
\end{aligned}$$

The following two lemmas can be obtained in a way similar to that given in Lemmas \ref{6}-\ref{13}.
\begin{lemma}\label{40-1}
    For any $ t<T^*$ and $\eta>0$, it holds that  
    $$\begin{aligned}    \sum\limits_{k\in\mathbb{Z}}&2^k\left(\sqrt{|\mathcal{L}_1^k+\mathcal{L}_2^k+\mathcal{L}_{3}^k|}\right)\gtrsim\|\bar{\theta}_\Phi(t)\|_{B^{1,0}}-\|e^{\delta|D|}\partial_y \theta_0\|_{B^{1,0}_{\Psi_0}}+(\sqrt{\lambda}-C_\eta)\|\bar{\theta}_\Phi\|_{\widetilde{L}_{t,\dot{\mu}}^2(B^{\frac{3}{2},0})}\\
    &+\|\sqrt{-(\partial_t\Psi+4\theta^E(\partial_y\Psi)^2)}\bar{\theta}_\Phi\|_{\widetilde{L}_t^2(B^{1,0})}+(\frac{\sqrt{\theta^E}}{2}-\epsilon^{\frac{1}{2}}\langle t\rangle^{\frac{1}{8}}-\eta)\|\partial_y\bar{\theta}_\Phi\|_{\widetilde{L}_t^2(B^{1,0})}-\eta\|\partial_y\theta_\Phi\|_{\widetilde{L}_t^2(B^{1,0})}.
    \end{aligned}$$
\end{lemma}
\begin{lemma}\label{7.2}
    For any $ t<T^*$ and $\eta>0$, the terms $\{{\cal R}_j^k\}_{j\le 6}$ given in \eqref{39} satisfy the estimate  $$\sum\limits_{j=1}^6\sum\limits_{k\in\mathbb{Z}}2^k\sqrt{|\mathcal{R}_j^k|}\lesssim \eta\|\partial_y\bar{\theta}_\Phi\|_{\widetilde{L}_t^2(B^{1,0})}+C_\eta\|\bar{\theta}_\Phi\|_{\widetilde{L}_{t,\dot{\mu}}^2(B^{\frac{3}{2},0})}+\|u_\Phi\|_{\widetilde{L}_{t,\dot{\mu}}^2(B^{\frac{3}{2},1})}.$$
\end{lemma}

The estimate of the remaining term $\mathcal{R}_7^k$ relies on the following lemma.

\begin{lemma}\label{37}
For any $ t<T^*$ and $\eta>0$, it holds that
 \begin{equation}\label{36} \sum\limits_{k\in\mathbb{Z}}2^k\sqrt{|\mathcal{R}_7^k|}\lesssim \eta\|\partial_y\bar{u}_\Phi\|_{{\widetilde{L}}^2(B^{1,0})}+C_\eta\|\bar{\theta}_\Phi\|_{\widetilde{L}_{t,\dot{\mu}}^2(B^{\frac{3}{2},0})}+\|\partial_yu_\Phi\|_{\widetilde{L}_{t,\dot{\mu}}^2(B^{\frac{3}{2},0})}+\|\theta_\Phi\|_{\widetilde{L}_{t,\dot{\mu}}^2(B^{\frac{3}{2},0})}.
 \end{equation}    
\end{lemma}

\begin{proof}
We write
$$\begin{aligned}
\mathcal{R}_7^k=&2(e^\Psi\Delta_k[(\theta+\theta^E)\partial_yu\partial_y^2u]_\Phi|e^\Psi\Delta_k\bar{\theta}_\Phi)\\
=&2(e^\Psi\Delta_k[\theta\partial_yu\partial_y^2u]_\Phi|e^\Psi\Delta_k\bar{\theta}_\Phi)+2\theta^E(e^\Psi\Delta_k[\partial_yu\partial_y^2u]_\Phi|e^\Psi\Delta_k\bar{\theta}_\Phi)\\
=&:2I_1^k+2I_2^k.
\end{aligned}$$
By applying Bony's decomposition to $\theta\partial_yu\partial_y^2u$ it gives
$$\begin{aligned}
I_1^k=&(e^\Psi\Delta_k[T_{\theta\partial_yu}\partial_y^2u]_\Phi|e^\Psi\Delta_k\bar{\theta}_\Phi)+(e^\Psi\Delta_k[T_{\partial_y^2u}\theta\partial_yu]_\Phi|e^\Psi\Delta_k\bar{\theta}_\Phi)\\
&+(e^\Psi\Delta_k[R(\theta\partial_yu,\partial_y^2u)]_\Phi|e^\Psi\Delta_k\bar{\theta}_\Phi).
\end{aligned}$$
In a way similar to that of (\ref{27}), (\ref{28}), we get
$$\begin{aligned}
    \|\Delta_{k}[\theta\partial_yu]_\Phi\|_{L_x^2}=&\|\int_y^\infty\Delta_{k}[\bar{\theta}\partial_yu+\theta\partial_y^2u]_\Phi dy'\|_{L_x^2}\\
    \lesssim&\sum\limits_{|l-k|\leq 4}(\|\bar{\theta}_\Phi\|_{B^{\frac{1}{2},0}}\|\Delta_l\partial_yu_\Phi\|_{L_+^2}+\|\partial_yu_\Phi\|_{B^{\frac{1}{2},0}}\|\Delta_l\bar{\theta}_\Phi\|_{L_+^2}\\
    &+\|\theta_\Phi\|_{B^{\frac{1}{2},0}}\|\Delta_l\partial_y^2u_\Phi\|_{L_+^2}+\|\partial_y^2u_\Phi\|_{B^{\frac{1}{2},0}}\|\Delta_l\theta_\Phi\|_{L_+^2})\\
    &+\sum\limits_{l\geq k-3}(\|\bar{\theta}_\Phi\|_{B^{\frac{1}{2},0}}\|\Delta_l\partial_yu_\Phi\|_{L_+^2}+\|\partial_y^2u_\Phi\|_{B^{\frac{1}{2},0}}\|\Delta_l\theta_\Phi\|_{L_+^2}).
\end{aligned}$$
It follows in a way similar to the proof of Lemma \ref{19} that
$$\sum\limits_{k\in\mathbb{Z}}2^k\left(\sqrt{|I_1^k|}+\sqrt{|I_2^k|}\right)\lesssim\eta\|\partial_y\bar{u}_\Phi\|_{{\widetilde{L}}^2(B^{1,0})}+C_\eta\|\bar{\theta}_\Phi\|_{\widetilde{L}_{t,\dot{\mu}}^2(B^{\frac{3}{2},0})}+\|\partial_yu_\Phi\|_{\widetilde{L}_{t,\dot{\mu}}^2(B^{\frac{3}{2},0})}+\|\theta_\Phi\|_{\widetilde{L}_{t,\dot{\mu}}^2(B^{\frac{3}{2},0})}.$$
This completes the proof of Lemma \ref{37}.
\end{proof}

\begin{proof}[Proof of Proposition \ref{38}]
We deduce from (\ref{39}) that $$\sum\limits_{k\in\mathbb{Z}}2^k\left(\sqrt{|\mathcal{L}_1^k+\mathcal{L}_2^k+\mathcal{L}_{3}^k|}\right)\lesssim\sum\limits_{j=1}^7\sum\limits_{k\in\mathbb{Z}}2^k\sqrt{|\mathcal{R}_j^k|}.$$
Plugging the estimates in Lemmas \ref{40-1}-\ref{37} to the above inequality, we conclude the estimate (\ref{41}). 
\end{proof}

\subsection{Estimate of \texorpdfstring{$\partial_y^2 u$}{partial y2 u}}
From (\ref{1}), we know that $W=\partial_y^2 u$ satisfies the following equation 
$$\begin{aligned}
&\partial_tW+u\partial_xW+(\partial_x\partial_yu+\partial_y^3\theta)\partial_yu+(-\partial_xu+4(\partial_yu)^2+\partial_y^2\theta) W\\
&\qquad+(-\int_0^y\partial_xudy'+\int_0^y(\partial_yu)^2dy'-\partial_y\theta)\partial_yW=(\theta+\theta^E)\partial_y^2W,\end{aligned}$$
which implies that 
 $W_\Phi=
\mathcal{F}^{-1}_{\xi\rightarrow x}(e^{\Phi(t,\xi)}{\hat W}(t,\xi,y))$ satisfies 
$$\begin{aligned}
\partial_tW_\Phi+\lambda\dot{\mu}|D|W_\Phi+[u\partial_xW]_\Phi+[(\partial_x\partial_yu+\partial_y^3\theta)\partial_yu]_\Phi+[(-\partial_xu+4(\partial_yu)^2+\partial_y^2\theta) W]_\Phi&\\
+[(-\int_0^y\partial_xudy'+\int_0^y(\partial_yu)^2dy'-\partial_y\theta)\partial_yW]_\Phi=[(\theta+\theta^E)\partial_y^2W]_\Phi.&
\end{aligned}$$
By applying the dyadic operator $\Delta_k$ to the above equation and then take the $L_{t,+}^2$ inner product of the resulting equation with $e^{2\Psi}\Delta_kW_\Phi$, it follows
\begin{equation}\label{81}
\begin{aligned}
&(e^\Psi\Delta_k\partial_tW_\Phi|e^\Psi\Delta_kW_\Phi)+\lambda(\dot{\mu}|D|e^\Psi\Delta_kW_\Phi|e^\Psi\Delta_kW_\Phi)-(e^\Psi\Delta_k[(\theta+\theta^E)\partial_y^2W]_\Phi|e^\Psi\Delta_kW_\Phi)\\
=&-(e^\Psi\Delta_k[u\partial_xW]_\Phi|e^\Psi\Delta_kW_\Phi)-(e^\Psi\Delta_k[(\partial_x\partial_yu+\partial_y^3\theta)\partial_yu]_\Phi|e^\Psi\Delta_kW_\Phi)\\
&-(e^\Psi\Delta_k[(-\partial_xu+4(\partial_yu)^2+\partial_y^2\theta) W]_\Phi|e^\Psi\Delta_kW_\Phi)\\
&-(e^\Psi\Delta_k[(-\int_0^y\partial_xudy'+\int_0^y(\partial_yu)^2dy'-\partial_y\theta)\partial_yW]_\Phi|e^\Psi\Delta_kW_\Phi).
\end{aligned}\end{equation}
Denote the three terms on the left side by $\mathcal{L}_1^k,\mathcal{L}_2^k,\mathcal{L}_3^k$ and the four terms on the right side by $\mathcal{R}_1^k,\cdots,\mathcal{R}_4^k$.

By using the fact that $W|_{y=0}=0$ from (\ref{1}) and integrating by parts, we obtain that
$$\begin{aligned}
    \mathcal{L}_3^k=&(e^\Psi\Delta_k[(\theta+\theta^E)\partial_y W]_\Phi|e^\Psi\Delta_k \partial_y W_\Phi)\\
    &+(2e^\Psi\partial_y\Psi\Delta_k[(\theta+\theta^E)\partial_y W]_\Phi|e^\Psi\Delta_k W_\Phi)+(e^\Psi\Delta_k[\partial_y\theta\partial_y W]_\Phi|e^\Psi\Delta_k W_\Phi).
\end{aligned}$$

\begin{lemma}\label{82}
    For any $ t<T^*$ and $\eta>0$, we have the following estimate  
    $$\begin{aligned}    \sum\limits_{k\in\mathbb{Z}}&2^\frac{k}2\left(\sqrt{|\mathcal{L}_1^k+\mathcal{L}_2^k+\mathcal{L}_{3}^k|}\right)\gtrsim\|W_\Phi(t)\|_{B^{\frac{1}{2},0}}-\|e^{\delta|D|}\partial_y^2u_0\|_{B^{\frac{1}{2},0}_{\Psi_0}}+(\sqrt{\lambda}-C_\eta)\|W_\Phi\|_{\widetilde{L}_{t,\dot{\mu}}^2(B^{1,0})}\\
    &+\|\sqrt{-(\partial_t\Psi+4\theta^E(\partial_y\Psi)^2)}W_\Phi\|_{\widetilde{L}_t^2(B^{\frac{1}{2},0})}+(\frac{\sqrt{\theta^E}}{2}-\epsilon^{\frac{1}{2}}\langle t\rangle^{\frac{1}{8}}-\eta)\|\partial_yW_\Phi\|_{\widetilde{L}_t^2(B^{\frac{1}{2},0})}-\|\theta_\Phi\|_{\widetilde{L}_{t,\dot{\mu}}^2(B^{1,2})}.
    \end{aligned}$$
\end{lemma}

Indeed, for the term  $(e^\Psi\partial_y\Psi\Delta_k[T_{\partial_yW}\theta]_\Phi|e^\Psi\Delta_kW_\Phi)$ appeared in ${\cal L}_3^k$, we have
    $$\begin{aligned}
    &|(e^\Psi\partial_y\Psi\Delta_k[T_{\partial_yW}\theta]_\Phi|e^\Psi\Delta_kW_\Phi)|\\
    \lesssim&\int_0^t\int_0^\infty e^{2\Psi}\sum\limits_{|k'-k|\leq 4}\sum\limits_{k''\leq k'-2}2^{\frac{k''}{2}}\|\Delta_{k''}\partial_yW_\Phi\|_{L_x^2}\|\partial_y\Psi\Delta_{k'}\theta_\Phi\|_{L_x^2}\|\Delta_kW_\Phi\|_{L_x^2}dydt'\\
    \lesssim&\int_0^t\sum\limits_{|k'-k|\leq 4}\sum\limits_{k''\leq k'-2}2^{\frac{k''}{2}}\langle t\rangle^{\frac{1}{4}}\|e^\Psi\Delta_{k''}\partial_yW_\Phi\|_{L_+^2}(\|\partial_y\Psi e^\Psi\Delta_{k'}\partial_y\theta_\Phi\|_{L_+^2}+\langle t\rangle^{-1}\|e^\Psi\Delta_{k'}\theta_\Phi\|_{L_+^2})\\
    &\times\|e^\Psi\Delta_kW_\Phi\|_{L_+^2}dt'\\
    \lesssim&\int_0^t\sum\limits_{|k'-k|\leq 4}\sum\limits_{k''\leq k'-2}2^{\frac{k''}{2}}\langle t\rangle^{\frac{1}{4}}\|e^\Psi\Delta_{k''}\partial_yW_\Phi\|_{L_+^2}(\|e^\Psi\Delta_{k'}\partial_y^2\theta_\Phi\|_{L_+^2}+\langle t\rangle^{-1}\|e^\Psi\Delta_{k'}\theta_\Phi\|_{L_+^2})\\
    &\times\|e^\Psi\Delta_kW_\Phi\|_{L_+^2}dt'\\
    \lesssim&\sum\limits_{|k'-k|\leq 4}\left(\int_0^t\dot{\mu}\|e^\Psi\Delta_kW_\Phi\|_{L_+^2}^2dt'\right)^{\frac{1}{2}}\left(\int_0^t\dot{\mu}\|e^\Psi\Delta_{k'}\partial_y^2\theta_\Phi\|_{L_+^2}^2dt'\right)^{\frac{1}{2}}\\
    &+\sum\limits_{|k'-k|\leq 4}\left(\int_0^t\dot{\mu}\|e^\Psi\Delta_kW_\Phi\|_{L_+^2}^2dt'\right)^{\frac{1}{2}}\left(\int_0^t\dot{\mu}\|e^\Psi\Delta_{k'}\theta_\Phi\|_{L_+^2}^2dt'\right)^{\frac{1}{2}},
    \end{aligned}$$
    which implies $$\sum\limits_{k\in\mathbb{Z}}2^{\frac{k}{2}}\sqrt{|(e^\Psi\partial_y\Psi\Delta_k[T_{\partial_yW}\theta]_\Phi|e^\Psi\Delta_kW_\Phi)|}\lesssim \|W_\Phi\|_{\widetilde{L}_{t,\dot{\mu}}^2(B^{1,0})}+\|\partial_y^2\theta_\Phi\|_{\widetilde{L}_{t,\dot{\mu}}^2(B^{1,0})}+\|\theta_\Phi\|_{\widetilde{L}_{t,\dot{\mu}}^2(B^{1,0})}.$$
Other terms in $\{{\cal L}_j^k\}_{j\le 3}$ can be studied in a way similar to the proofs of Lemmas \ref{6}-\ref{13}, and conclude the estimate given in Lemma \ref{82}. Similarly, one can have

\begin{lemma}\label{83}
    For any $ t<T^*$ and $\eta>0$, it holds that  
    $$\begin{aligned}\sum\limits_{j=1}^4\sum\limits_{k\in\mathbb{Z}}2^{\frac{k}{2}}\sqrt{|\mathcal{R}_j^k|}\lesssim & \eta(\|\partial_yW_\Phi\|_{\widetilde{L}_t^2(B^{\frac{1}{2},0})}+\|\partial_yS_\Phi\|_{\widetilde{L}_t^2(B^{\frac{1}{2},0})})+C_\eta\|W_\Phi\|_{\widetilde{L}_{t,\dot{\mu}}^2(B^{1,0})}\\
    &+\|S_\Phi\|_{\widetilde{L}_{t,\dot{\mu}}^2(B^{1,0})}+\|u_\Phi\|_{\widetilde{L}_{t,\dot{\mu}}^2(B^{1,1})}.\end{aligned}$$  
\end{lemma}

\begin{proof}[Proof of Proposition \ref{79}]
We deduce from (\ref{81}) that 
$$\sum\limits_{k\in\mathbb{Z}}2^k\left(\sqrt{|\mathcal{L}_1^k+\mathcal{L}_2^k+\mathcal{L}_{3}^k|}\right)\lesssim\sum\limits_{j=1}^4\sum\limits_{k\in\mathbb{Z}}2^k\sqrt{|\mathcal{R}_j^k|}.$$
Combining Lemmas \ref{82}-\ref{83}, it leads to estimate (\ref{80}) given in Proposition \ref{79}.
\end{proof}

\subsection{Estimate of \texorpdfstring{$\partial_y^2 \theta$}{partial y2 theta}}\label{77}

Denote by $S=\partial_y^2\theta$. At this time, we do not have the Dirichlet boundary condition of $S$, but  from (\ref{85}), we get 
\begin{equation}\label{9.1}
[\partial_yu\partial_x\theta-(\theta+\theta^E)\partial_y S]|_{y=0}=0.
\end{equation}
Thus in order to use this boundary condition, we write the equation of $S$ from (\ref{85}) in the form 
$$\begin{aligned}
&\partial_tS+\partial_yu\partial_x\partial_y\theta-\partial_y\theta\partial_x\partial_yu+u\partial_xS+(-2\partial_xu+\partial_y^2\theta+(\partial_yu)^2)S\\
&+(-\int_0^y\partial_xudy'+\partial_y\theta+\int_0^y(\partial_yu)^2dy')\partial_yS\\
&-2\partial_y\theta\partial_yu\partial_y^2u-2(\theta+\theta^E)(\partial_y^2u)^2-2(\theta+\theta^E)\partial_yu\partial_y^3u\\
&=\partial_y((\theta+\theta^E)\partial_yS-\partial_yu\partial_x\theta),
\end{aligned}$$
which implies that 
 $S_\Phi=
\mathcal{F}^{-1}_{\xi\rightarrow x}(e^{\Phi(t,\xi)}{\hat S}(t,\xi,y))$ satisfies the following equation
$$\begin{aligned}
&\partial_tS_\Phi+\lambda\dot{\mu}|D|S_\Phi+[\partial_yu\partial_x\partial_y\theta-\partial_y\theta\partial_x\partial_yu]_\Phi+[u\partial_xS]_\Phi\\
&+[(-2\partial_xu+\partial_y^2\theta+(\partial_yu)^2)S]_\Phi+[(-\int_0^y\partial_xudy'+\partial_y\theta+\int_0^y(\partial_yu)^2dy')\partial_yS]_\Phi\\
&-[2\partial_y\theta\partial_yu\partial_y^2u]_\Phi-[2(\theta+\theta^E)(\partial_y^2u)^2]_\Phi-[2(\theta+\theta^E)\partial_yu\partial_y^3u]_\Phi\\
=&[\partial_y((\theta+\theta^E)\partial_yS-\partial_yu\partial_x\theta)]_\Phi.
\end{aligned}$$

By applying the dyadic operator $\Delta_k$ to the above equation and then take the $L_{t,+}^2$ inner product of the resulting equation with $e^{2\Psi}\Delta_kS_\Phi$, we obtain
\begin{equation}\label{95}\begin{aligned}
&(e^\Psi\Delta_k\partial_tS_\Phi|e^\Psi\Delta_kS_\Phi)+(\lambda\dot{\mu}|D|e^\Psi\Delta_kS_\Phi|e^\Psi\Delta_kS_\Phi)-(e^\Psi\Delta_k[\partial_y((\theta+\theta^E)\partial_yS-\partial_yu\partial_x\theta)]_\Phi|e^\Psi\Delta_kS_\Phi)\\
=&-(e^\Psi\Delta_k[\partial_yu\partial_x\partial_y\theta-\partial_y\theta\partial_x\partial_yu]_\Phi|e^\Psi\Delta_kS_\Phi)-(e^\Psi\Delta_k[u\partial_xS]_\Phi|e^\Psi\Delta_kS_\Phi)\\
&-(e^\Psi\Delta_k[(-2\partial_xu+\partial_y^2\theta+(\partial_yu)^2)S]_\Phi|e^\Psi\Delta_kS_\Phi)\\
&-(e^\Psi\Delta_k[(-\int_0^y\partial_xudy'+\partial_y\theta+\int_0^y(\partial_yu)^2dy')\partial_yS]_\Phi|e^\Psi\Delta_kS_\Phi)\\
&+(e^\Psi\Delta_k[2\partial_y\theta\partial_yu\partial_y^2u]_\Phi|e^\Psi\Delta_kS_\Phi)+(e^\Psi\Delta_k[2(\theta+\theta^E)(\partial_y^2u)^2]_\Phi|e^\Psi\Delta_kS_\Phi)\\
&+(e^\Psi\Delta_k[2(\theta+\theta^E)\partial_yu\partial_y^3u]_\Phi|e^\Psi\Delta_kS_\Phi).
\end{aligned}\end{equation}
Denote the three terms on the left side by $\mathcal{L}_1^k,\mathcal{L}_2^k,\mathcal{L}_3^k$ and the seven terms on the right side by $\mathcal{R}_1^k,\cdots,\mathcal{R}_7^k$.

By integrating by parts, and using \eqref{9.1} we obtain that
\begin{equation}\label{84}
\begin{aligned}
    \mathcal{L}_3^k=&(e^\Psi\Delta_k[(\theta+\theta^E)\partial_y S]_\Phi|e^\Psi\Delta_k \partial_y S_\Phi)+(2e^\Psi\partial_y\Psi\Delta_k[(\theta+\theta^E)\partial_y S]_\Phi|e^\Psi\Delta_k S_\Phi)\\
    &-(e^\Psi\Delta_k[\partial_yu\partial_x\theta]_\Phi|e^\Psi\Delta_k \partial_y S_\Phi)-(2e^\Psi\partial_y\Psi\Delta_k[\partial_yu\partial_x\theta]_\Phi|e^\Psi\Delta_k S_\Phi)\\
    =:&  \mathcal{L}_{3,1}^k+\mathcal{L}_{3,2}^k+\mathcal{L}_{3,3}^k+\mathcal{L}_{3,4}^k,
\end{aligned}
\end{equation}
with obvious notations $\{\mathcal{L}_{3,j}^k\}_{j\le 4}$. 

\begin{lemma}\label{125}
For any $t<T^*$ and $\eta>0$, it holds that
\begin{equation}\label{124}
    \sum\limits_{k\in\mathbb{Z}}2^\frac{k}{2}\sqrt{|\mathcal{L}_{3,3}^k|}\lesssim\eta\|\partial_yS_\Phi\|_{\widetilde{L}_{t}^2(B^{\frac{1}{2},0})}+C_\eta(\|\partial_yu_\Phi\|_{\widetilde{L}_{t,\dot{\mu}}^2(B^{\frac{3}{2},0})}+\|\partial_y\theta_\Phi\|_{\widetilde{L}_{t,\dot{\mu}}^2(B^{\frac{3}{2},0})}).
\end{equation}
\end{lemma}
\begin{proof}
By Young's inequality, we have
\begin{equation}\label{119}
\begin{aligned}|\mathcal{L}_{3,3}^k|\lesssim\eta\|e^\Psi\Delta_k \partial_y S_\Phi\|^2_{L_{t,+}^2}+C_\eta\|e^\Psi\Delta_k[\partial_yu\partial_x\theta]_\Phi\|^2_{L^2_{t,+}}
=:\eta I_1^k+C_\eta I_2^k.
\end{aligned}
\end{equation}
Obviously, one has
\begin{equation}\label{120}
\sum\limits_{k\in\mathbb{Z}}2^{\frac{k}{2}}\sqrt{|I_1^k|}=\|\partial_yS_\Phi\|_{\widetilde{L}_{t}^2(B^{\frac{1}{2},0})}.\end{equation}
On the other hand, by applying Bony's decomposition to $\partial_yu\partial_x\theta$, it follows
\begin{equation}\label{121}
\begin{aligned}
I_2^k=&(e^\Psi\Delta_k[T_{\partial_yu}\partial_x\theta]_\Phi|e^\Psi\Delta_k[T_{\partial_yu}\partial_x\theta]_\Phi)+(e^\Psi\Delta_k[T_{\partial_yu}\partial_x\theta]_\Phi|e^\Psi\Delta_k[T_{\partial_x\theta}\partial_yu]_\Phi)\\
&+(e^\Psi\Delta_k[T_{\partial_yu}\partial_x\theta]_\Phi|e^\Psi\Delta_k[R(\partial_yu,\partial_x\theta)]_\Phi)+(e^\Psi\Delta_k[T_{\partial_x\theta}\partial_yu]_\Phi|e^\Psi\Delta_k[T_{\partial_yu}\partial_x\theta]_\Phi)\\
&+(e^\Psi\Delta_k[T_{\partial_x\theta}\partial_yu]_\Phi|e^\Psi\Delta_k[T_{\partial_x\theta}\partial_yu]_\Phi)+(e^\Psi\Delta_k[T_{\partial_x\theta}\partial_yu]_\Phi|e^\Psi\Delta_k[R(\partial_yu,\partial_x\theta)]_\Phi)\\
&+(e^\Psi\Delta_k[R(\partial_yu,\partial_x\theta)]_\Phi|e^\Psi\Delta_k[T_{\partial_yu}\partial_x\theta]_\Phi)+(e^\Psi\Delta_k[R(\partial_yu,\partial_x\theta)]_\Phi|e^\Psi\Delta_k[T_{\partial_x\theta}\partial_yu]_\Phi)\\
&+(e^\Psi\Delta_k[R(\partial_yu,\partial_x\theta)]_\Phi|e^\Psi\Delta_k[R(\partial_yu,\partial_x\theta)]_\Phi)\\
=&:J_1^k+\cdots+J_9^k.
\end{aligned}
\end{equation}
By a similar derivation of Lemma \ref{12}, we first get that
$$\begin{aligned}
|J_1^k|\lesssim&\int_0^t\sum\limits_{|k'-k|\leq 4}2^{k'}\langle t\rangle^{\frac{1}{4}}\|\partial_yu_\Phi\|_{B^{\frac{1}{2},0}}\|e^\Psi\Delta_{k'}\partial_y\theta_\Phi\|_{L_+^2}\sum\limits_{|l-k|\leq 4}2^{l}\langle t\rangle^{\frac{1}{4}}\|\partial_yu_\Phi\|_{B^{\frac{1}{2},0}}\|e^\Psi\Delta_{l}\partial_y\theta_\Phi\|_{L_+^2}dt'\\
\leq&\sum\limits_{|k'-k|\leq 4}\sum\limits_{|l-k|\leq 4}2^{k'+l}\left(\int_0^t\dot{\mu}\|e^\Psi\Delta_{k'}\partial_y\theta_\Phi\|^2_{L_+^2}dt'\right)^\frac{1}{2}\left(\int_0^t\dot{\mu}\|e^\Psi\Delta_l\partial_y\theta_\Phi\|^2_{L_+^2}dt'\right)^\frac{1}{2}.
\end{aligned}$$
By taking square root of the resulting inequality, muliplying by $2^{\frac{k}{2}}$ and summing over $k\in\mathbb{Z}$, we find that
\begin{equation}\label{122}
\sum\limits_{k\in\mathbb{Z}}2^{\frac{k}{2}}\sqrt{|J_1^k|}\lesssim\|\partial_y\theta_\Phi\|_{\widetilde{L}_{t,\dot{\mu}}^2(B^{\frac{3}{2},0})}.
\end{equation}
In the  same way as above, we have
\begin{equation}\label{123}
\sum_{j=2}^9\sum\limits_{k\in\mathbb{Z}}2^{\frac{k}{2}}\sqrt{|J_j^k|}\lesssim\|\partial_yu_\Phi\|_{\widetilde{L}_{t,\dot{\mu}}^2(B^{\frac{3}{2},0})}+\|\partial_y\theta_\Phi\|_{\widetilde{L}_{t,\dot{\mu}}^2(B^{\frac{3}{2},0})}.
\end{equation}
Inserting the estimates (\ref{122}) and (\ref{123}) into (\ref{121}) and combining (\ref{120}), we obtain the estimate (\ref{124}) from (\ref{119}). 
\end{proof}

By using similar arguments for the proofs of Lemmas \ref{6}-\ref{13}, and using Lemma \ref{125}, we can get the following two estimates.
\begin{lemma}\label{91}
    For any $ t<T^*$ and $\eta>0$, it holds that  
    $$\begin{aligned}    \sum\limits_{k\in\mathbb{Z}}&2^{\frac{k}{2}}\left(\sqrt{|\mathcal{L}_1^k+\mathcal{L}_2^k+\mathcal{L}_{3}^k|}\right)\gtrsim\|S_\Phi(t)\|_{B^{\frac{1}{2},0}}-\|e^{\delta|D|}\partial_y^2\theta_0\|_{B^{\frac{1}{2},0}_{\Psi_0}}+(\sqrt{\lambda}-C_\eta)\|S_\Phi\|_{\widetilde{L}_{t,\dot{\mu}}^2(B^{1,0})}\\
    &+\|\sqrt{-(\partial_t\Psi+4\theta^E(\partial_y\Psi)^2)}S_\Phi\|_{\widetilde{L}_t^2(B^{\frac{1}{2},0})}+(\frac{\sqrt{\theta^E}}{2}-\epsilon^{\frac{1}{2}}\langle t\rangle^{\frac{1}{8}}-\eta)\|\partial_yS_\Phi\|_{\widetilde{L}_t^2(B^{\frac{1}{2},0})}\\
    &-C_\eta(\|\partial_yu_\Phi\|_{\widetilde{L}_{t,\dot{\mu}}^2(B^{\frac{3}{2},0})}+\|\partial_y\theta_\Phi\|_{\widetilde{L}_{t,\dot{\mu}}^2(B^{\frac{3}{2},0})})-\|\theta_\Phi\|_{\widetilde{L}_{t,\dot{\mu}}^2(B^{1,0})}.
    \end{aligned}$$
and
$$\begin{aligned}\sum\limits_{j=1}^7\sum\limits_{k\in\mathbb{Z}}&2^\frac{k}{2}\sqrt{|\mathcal{R}_j^k|}\lesssim \eta(\|\partial_yS_\Phi\|_{\widetilde{L}_t^2(B^{\frac{1}{2},0})}+\|\partial_yW_\Phi\|_{\widetilde{L}_t^2(B^{\frac{1}{2},0})})+C_\eta\|S_\Phi\|_{\widetilde{L}_{t,\dot{\mu}}^2(B^{1,0})}\\
    &+\|u_\Phi\|_{\widetilde{L}_{t,\dot{\mu}}^2(B^{1,2})}+\|\theta_\Phi\|_{\widetilde{L}_{t,\dot{\mu}}^2(B^{1,1})}.\end{aligned}$$  
\end{lemma}

\begin{proof}[Proof of Proposition \ref{93}]
We deduce from (\ref{95}) that 
$$\sum\limits_{k\in\mathbb{Z}}2^{\frac{k}{2}}\left(\sqrt{|\mathcal{L}_1^k+\mathcal{L}_2^k+\mathcal{L}_{3}^k|}\right)\lesssim\sum\limits_{j=1}^7\sum\limits_{k\in\mathbb{Z}}2^{\frac{k}{2}}\sqrt{|\mathcal{R}_j^k|}.$$
Using the estimates given in Lemmas \ref{91}, it follows the estimate  given in Proposition \ref{93}.
\end{proof}

\section{Proof of the main result}

The goal of this section is to prove the existence and uniqueness of the solution to the boundary layer problem \eqref{1} given in Theorem \ref{50}.

\subsection{Existence of solutions to the compressible Prandtl equations}

In this subsection, we will prove the existence part stated in Theorem \ref{50}. For this, ﬁrst we construct a regularized  boundary layer equation with adding artificial viscosity to obtain an appropriate approximate solutions, then deduce uniform estimates for such approximate solution sequence, and ﬁnally pass to the limit of the approximate problems.

\noindent\emph{Step 1}. Associating with the problem \eqref{1}, we introduce the following regularized boundary layer equation by adding artificial viscosity $\nu>0$, 
\begin{equation}\label{56}\left\{
    \begin{aligned}
        &\partial_t u^\nu+u^\nu\partial_x u^\nu+v^\nu\partial_y u^\nu=\nu\partial_x^2 u^\nu+(\theta^\nu+\theta^E)\partial_y^2 u^\nu,\\
        &\partial_t \theta^\nu+u^\nu\partial_x \theta^\nu+v^\nu\partial_y \theta^\nu=\nu\partial_x^2\theta^\nu+(\theta^\nu+\theta^E)\partial_y^2 \theta^\nu+(\theta^\nu+\theta^E)(\partial_y u^\nu)^2,\\
        &\partial_x u^\nu+\partial_y v^\nu=\partial_y^2\theta^\nu+(\partial_y u^\nu)^2,\\
        &u^\nu|_{y=0}=v^\nu|_{y=0}=\partial_y\theta^\nu|_{y=0}=0,\quad \lim\limits_{y\to\infty}u^\nu=\lim\limits_{y\to\infty}\theta^\nu=0,\\
        &(u^\nu,\theta^\nu)|_{t=0}=(u_0,\theta_0)(x,y).
    \end{aligned}
    \right.
\end{equation}
\emph{Step 2}. We introduce a new phase function 
$$\widetilde{\Phi}(t,\xi)=(\frac{\delta}{2}-\lambda\mu(t))|\xi|,$$
which yields
$$\begin{aligned}&\|u^\nu_{\widetilde{\Phi}}\|_{\widetilde{L}_t^\infty(B^{3,1})}\lesssim\|u^\nu_\Phi\|_{\widetilde{L}_t^\infty(B^{1,1})},\|\partial_y u^\nu_{\widetilde{\Phi}}\|_{\widetilde{L}_t^2(B^{3,1})}\lesssim\|\partial_y u^\nu_{\Phi}\|_{\widetilde{L}_t^2(B^{1,1})},\\
&\|\theta^\nu_{\widetilde{\Phi}}\|_{\widetilde{L}_t^\infty(B^{3,1})}\lesssim\|\theta^\nu_\Phi\|_{\widetilde{L}_t^\infty(B^{1,1})},\|\partial_y \theta^\nu_{\widetilde{\Phi}}\|_{\widetilde{L}_t^2(B^{3,1})}\lesssim\|\partial_y \theta^\nu_{\Phi}\|_{\widetilde{L}_t^2(B^{1,1})}.
\end{aligned}$$
One can establish the same estimate as given in \eqref{45} uniformly in $\nu$ for the solution of \eqref{56}, which implies that for any $ t<T,~\nu>0$, $$\|u^\nu_{\widetilde{\Phi}}\|_{\widetilde{L}_t^\infty(B^{3,1})}+\|\partial_y u^\nu_{\widetilde{\Phi}}\|_{\widetilde{L}_t^2(B^{3,1})}+\|\theta^\nu_{\widetilde{\Phi}}\|_{\widetilde{L}_t^\infty(B^{3,1})}+\|\partial_y \theta^\nu_{\widetilde{\Phi}}\|_{\widetilde{L}_t^2(B^{3,1})}\leq C,$$
where $C$ is a positive constant only depending on the initial data $(u_0,\theta_0)$.
Therefore, from Definition \ref{48}, we infer that
\begin{equation}\label{72}\|u^\nu\|_{L^\infty(0,T;H^{3,1}(\mathbb{R}^2_+))}+\|u^\nu\|_{L^2(0,T;H^{3,2}(\mathbb{R}^2_+))}+\|\theta^\nu\|_{L^\infty(0,T;H^{3,1}(\mathbb{R}^2_+))}+\|\theta^\nu\|_{L^2(0,T;H^{3,2}(\mathbb{R}^2_+))}\leq C,~\forall\nu>0.\end{equation}
\emph{Step 3}. By using the Lions-Aubin lemma and the compact embedding of $H^2$ in $H^1_{\mathrm{loc}}$, to prove that  $(u^\nu,\theta^\nu)$ converges in $L^2(0,T;H^1_{\mathrm{loc}}(\mathbb{R}^2_+))$, it suffices to justify that $(\partial_tu^\nu,\partial_t\theta^\nu)$ are bounded in $L^2(0,T;L^2(\mathbb{R}^2_+))$ uniformly in $\nu$. In view of (\ref{56}), we write 
$$\begin{aligned}\partial_t u^\nu=&\nu\partial_x^2 u^\nu+(\theta^\nu+\theta^E)\partial_y^2 u^\nu-u^\nu\partial_x u^\nu-v^\nu\partial_y u^\nu\\
=&I_1+I_2+I_3+I_4,
\end{aligned}$$
and
$$\begin{aligned}\partial_t \theta^\nu=&\nu\partial_x^2 \theta^\nu+(\theta^\nu+\theta^E)\partial_y^2 \theta^\nu-u^\nu\partial_x \theta^\nu-v^\nu\partial_y \theta^\nu+(\theta^\nu+\theta^E)(\partial_yu^\nu)^2\\
=&J_1+J_2+J_3+J_4+J_5.
\end{aligned}$$
First, it is obvious that  $$\|I_1\|_{L^2(0,T;L^2(\mathbb{R}^2_+))}\leq\nu\|\partial_x^2 u^\nu\|_{L^2(0,T;L^2(\mathbb{R}^2_+))}\lesssim C.$$
From (\ref{5}), one has
$$\|I_2\|_{L^2(0,T;L^2(\mathbb{R}^2_+))}\leq\|\theta^\nu+\theta^E\|_{L^\infty(0,T;L^\infty(\mathbb{R}^2_+))}\|\partial_y^2 u^\nu\|_{L^2(0,T;L^2(\mathbb{R}^2_+))}\lesssim C.$$
By using the Sobolev embedding  $$\|u^\nu\|_{L^\infty(\mathbb{R}^2_+)}\lesssim\|u^\nu\|_{H^{1,1}(\mathbb{R}^2_+)},$$
we obtain 
$$\|I_3\|_{L^2(0,T;L^2(\mathbb{R}^2_+))}\leq\|u^\nu\|_{L^\infty(0,T;L^\infty(\mathbb{R}^2_+))}\|\partial_xu^\nu\|_{L^2(0,T;L^2(\mathbb{R}^2_+))}\lesssim C^2.$$
Finally, again from (\ref{56}), we get
$$I_4=\int_0^y\partial_xu^\nu dy'\partial_yu^\nu-\int_0^y(\partial_yu^\nu)^2dy'\partial_yu^\nu-\partial_y\theta^\nu\partial_yu^\nu.$$
It follows from Hölder's inequality and Definition \ref{58} that
$$\|\int_0^y\partial_xu^\nu dy'\|_{L^\infty(\mathbb{R}^2_+)}\leq\int_0^\infty\|\partial_xu^\nu\|_{L_x^\infty}dy
{
\lesssim
}
\langle t\rangle^\frac{1}{4}\|u^\nu\|_{B^{2,0}}.$$
Similarly, 
$$\|\int_0^y(\partial_yu^\nu)^2dy'\|_{L^\infty(\mathbb{R}^2_+)}\leq\int_0^\infty\|\partial_y u^\nu\|_{L_x^\infty}^2dy\lesssim\|u^\nu\|_{H^{1,1}(\mathbb{R}^2_+)}^2.$$
As a result, by Sobolev embedding, we get
$$\begin{aligned}
&\|I_4\|_{L^2(0,T;L^2(\mathbb{R}^2_+))}\\
\lesssim&\left(\|\int_0^y\partial_xu^\nu dy'\|_{L^\infty(0,T;L^\infty(\mathbb{R}^2_+))}+\|\int_0^y(\partial_yu^\nu)^2dy'\|_{L^\infty(0,T;L^\infty(\mathbb{R}^2_+))}\right)\|\partial_yu^\nu\|_{L^2(0,T;L^2(\mathbb{R}^2_+))}\\
&+\|\partial_y\theta^\nu\|_{L^\infty(0,T;L^2(\mathbb{R}^2_+))}\|\partial_yu^\nu\|_{L^2(0,T;L^\infty(\mathbb{R}^2_+))}\\
\lesssim&(\langle t\rangle^\frac{1}{4}\|u^\nu\|_{\widetilde{L}_t^\infty(B^{2,0})}+\|u^\nu\|_{L^\infty(0,T;H^{1,1}(\mathbb{R}^2_+))}^2)\|\partial_yu^\nu\|_{L^2(0,T;L^2(\mathbb{R}^2_+))}\\
&+\|\partial_y\theta^\nu\|_{L^\infty(0,T;L^2(\mathbb{R}^2_+))}\|\partial_yu^\nu\|_{L^2(0,T;H^{1,1}(\mathbb{R}^2_+))}\\
\lesssim&(\langle t\rangle^\frac{1}{4}C+C^2)C+C^2.
\end{aligned}$$
Similar to the estimates above, $J_1,\cdots,J_4$ are also uniformly bounded in $L^2(0,T;L^2(\mathbb{R}^2_+))$, while
$$\|J_5\|_{L^2(0,T;L^2(\mathbb{R}^2_+))}\leq\|\theta^E+\theta^\nu\|_{L^\infty(0,T;L^\infty(\mathbb{R}^2_+))}\|\partial_yu^\nu\|_{L^\infty(0,T;L^2(\mathbb{R}^2_+))}\|\partial_yu^\nu\|_{L^2(0,T;L^\infty(\mathbb{R}^2_+))}\lesssim C^2.$$
Summarizing the above estimates, we get that
$(\partial_tu^\nu,\partial_t\theta^\nu)$ are bounded in $L^2(0,T;L^2(\mathbb{R}^2_+))$ uniformly in $\nu$. Thus, there exists $(u,\theta)$ such that as $\nu\to 0$, by extracting subsequence if necessary of $(u^\nu,\theta^\nu)$, which we still denote it by itself,
\begin{equation}\label{59}
    (u^\nu,\theta^\nu)\to (u,\theta)\quad\text{in~}L^2(0,T;H^1_{\mathrm{loc}}(\mathbb{R}^2_+)).
\end{equation}
\emph{Step 4}. 
Let $K\subset [-a,a]\times[0,b]$ be an arbitrary compact set of $\mathbb{R}_+^2$. Due to (\ref{56}), we write $v^\nu=-\int_0^y\partial_xu^\nu dy'+\int_0^y(\partial_yu^\nu)^2dy'+\partial_y\theta^\nu$. Using Hölder's inequality and the strong convergence (\ref{59}), we have
$$\begin{aligned}
&\int_0^T\int_K|\int_0^y\partial_xu^\nu dy'-\int_0^y\partial_xu dy'|dxdydt
\lesssim\|\partial_xu^\nu-\partial_xu\|_{L^2(0,T;L^2([-a,a]\times[0,b]))}\longrightarrow 0,
\end{aligned}$$
as $\nu\to0$. Similarly, one has
$$\begin{aligned}
&\int_0^T\int_K|\int_0^y(\partial_yu^\nu)^2 dy'-\int_0^y(\partial_yu)^2 dy'|dxdydt+\int_0^T\int_K|(\partial_yu^\nu)^2-(\partial_yu)^2|dxdydt\\
&\qquad\qquad
\lesssim\|\partial_yu^\nu-\partial_yu\|_{L^2(0,T;L^2([-a,a]\times[0,b]))}\|\partial_yu^\nu+\partial_yu\|_{L^2(0,T;L^2([-a,a]\times[0,b]))}\longrightarrow 0.
\end{aligned}$$

Therefore, we deduce the strong convergence of $v^\nu$ and $(\partial_y u^\nu)^2$:
\begin{equation}v^\nu\to v,\quad (\partial_y u^\nu)^2\to(\partial_y u)^2\quad\text{in~}L^1(0,T;L^1(K)).\end{equation}
Using the strong convergence (\ref{59}), uniform bound (\ref{72}) and uniform bound of $v^\nu$ and $(\partial_y u^\nu)^2$ stated in \emph{Step 3}, we have, after taking a subsequence, as $\nu_k\to0^+$,
\begin{equation}\left\{
\begin{array}{ll}
(u^{\nu_k},\theta^{\nu_k})\to (u,\theta)&\text{in~}L^2(0,T;H_{loc}^1(\R_+^2)),\\[2mm]
(\partial_y^2 u^{\nu_k},\partial_y^2 \theta^{\nu_k}){\rightharpoonup}(\partial_y^2 u,\partial_y^2\theta)&\text{in~}L^2(0,T;L^2_{loc}(\R_+^2)),\\[2mm]
v^{\nu_k}\to v&\text{in~}L^2(0,T;L^2_{loc}(\R_+^2)),\\[2mm]
(\partial_y u^{\nu_k})^2\to (\partial_y u)^2&\text{in~}L^2(0,T;L^2_{loc}(\R_+^2)).
\end{array}
\right.\end{equation}
Thus, passing to the limit $\nu^k\to0^+$ in (\ref{56}), we know that the limit $(u,\theta)$ solves the problem (\ref{1}).

\subsection{Uniqueness of the solution to the compressible Prandtl equations}

This subsection is devoted to the proof of the uniqueness part of Theorem \ref{50}. Let $(u_1,\theta_1)$ and $(u_2,\theta_2)$ be two solutions of (\ref{1}) obtained in Theorem \ref{50}. Denote by $U:=u_1-u_2$, $V:=v_1-v_2$ and $\Theta:=\theta_1-\theta_2$. Then from (\ref{1}), we know that $(U,\Theta)$ satisfy the following problem:
\begin{equation}\label{55}\left\{
\begin{aligned}
    &\partial_t U+U\partial_x u_1+u_2\partial_x U+V\partial_y u_1+v_2\partial_y U=\Theta\partial_y^2 u_1+(\theta_2+\theta^E)\partial_y^2 U,\\
    &\partial_t \Theta+U\partial_x \theta_1+u_2\partial_x \Theta+V\partial_y \theta_1+v_2\partial_y \Theta=\Theta\partial_y^2 \theta_1+(\theta_2+\theta^E)\partial_y^2 \Theta\\
    &\hspace{1.7in}
    +\Theta(\partial_y u_1)^2+(\theta_2+\theta^E)\partial_yU\partial_y(u_1+u_2),\\
    &\partial_x U+\partial_y V=\partial_y^2\Theta+\partial_y U\partial_y(u_1+u_2),\\
    &U|_{y=0}=V|_{y=0}=\partial_y\Theta|_{y=0}=0,\quad \lim\limits_{y\to\infty}U=\lim\limits_{y\to\infty}\Theta=0,\\
    &(U,\Theta)|_{t=0}=0.
\end{aligned}\right.
\end{equation}

Let $\mu_i(t;t_0),i=1,2$, be determined respectively by
\begin{equation}
\left\{
    \begin{aligned}
        \dot{\mu}_i(t;t_0) &= \langle t\rangle^{\frac{1}{4}} 
        \big(\|\partial_y u_{i\Phi_i}\|_{B^{\frac{1}{2},1}} 
        + \|\partial_y \theta_{i\Phi_i}\|_{B^{\frac{1}{2},1}} 
        + \|u_{i\Phi_i}\|_{B^{1,0}} 
        + \|\theta_{i\Phi_i}\|_{B^{1,0}}\big) \\
        &\quad + \langle t\rangle^{\frac{1}{2}} 
        \big(\|u_{i\Phi_i}\|_{B^{1,1}}^2 
        + \|\theta_{i\Phi_i}\|_{B^{1,1}}^2\big) \\
        &\quad + \|\partial_y u_{i\Phi_i}\|_{B^{1,0}}^4 
        + \|\partial_y \theta_{i\Phi_i}\|_{B^{1,0}}^4 
        + \|\partial_y u_{i\Phi_i}\|_{B^{1,0}} 
        \|\partial_y^2 u_{i\Phi_i}\|_{B^{1,0}}, \\
        \mu_i|_{t=t_0} &= 0,
    \end{aligned}
\right.
\end{equation}
with $
{
\Phi_i(t,\xi;t_0)
}
=(\delta-\lambda \mu_i(t;t_0))|\xi| ~ (i=1,2)$, 
and define the phase function $\widehat{\Phi}$ by
$$\widehat{\Phi}(t,\xi;t_0)=(\frac{\delta}{2}-\lambda M(t;t_0))|\xi|$$
with
$$M(t;t_0)=\mu_1(t;t_0)+\mu_2(t;t_0)+\int_{t_0}^t\langle t'\rangle^{\frac{1}{2}}(\|\partial_y^2 u_{1\Phi_1}\|^2_{B^{1,0}}+\|\partial_y^2 u_{2\Phi_2}\|^2_{B^{1,0}}+\|\partial_y^2 \theta_{1\Phi_1}\|^2_{B^{1,0}}+\|\partial_y^2 \theta_{2\Phi_2}\|^2_{B^{1,0}})dt'.$$
For simplicity, we will omit the hat of $\widehat{\Phi}$ in this subsection. 

\begin{remark}
(1) From the a priori estimate (\ref{45}), one has 
$$\int_0^T\langle t\rangle^{\frac{1}{2}}(\|\partial_y^2u_{i\Phi_i}\|^2_{B^{1,0}}+\|\partial_y^2\theta_{i\Phi_i}\|^2_{B^{1,0}})dt\leq\langle T\rangle^{\frac{1}{2}}(\|\partial_y^2u_{i\Phi_i}\|^2_{\widetilde{L}_t^2(B^{1,0})}+\|\partial_y^2\theta_{i\Phi_i}\|^2_{\widetilde{L}_t^2(B^{1,0})})<+\infty,$$
which implies that for fixed $\delta>0, \lambda>0$, there is $\tilde \delta>0$ such that whenever $0<t-t_0<\tilde \delta$,
$$\frac{\delta}{2}-\lambda M(t;t_0)>0,$$
so that there holds
$$\Phi(t,\xi;t_0)\leq\Phi(t,\xi-\eta;t_0)+\Phi(t,\eta;t_0),\quad\forall\xi,\eta\in\mathbb{R}.$$
In what follows, we shall always take $t_0=0$ since the uniqueness for the whole time period of existence can be deduced by a continuous argument.

\vspace{.1in}
(2) Since $\Phi(t)\leq\min(\Phi_1(t),\Phi_2(t))$, we have for $i=1,2$,
$$\begin{aligned}
&\|u_{i\Phi}\|_{\widetilde{L}_t^\infty(B^{1,1})}\leq\|u_{i\Phi_i}\|_{\widetilde{L}_t^\infty(B^{1,1})},\|\theta_{i\Phi}\|_{\widetilde{L}_t^\infty(B^{1,1})}\leq\|\theta_{i\Phi_i}\|_{\widetilde{L}_t^\infty(B^{1,1})},\\
&\|\partial_yu_{i\Phi}\|_{\widetilde{L}_t^2(B^{1,1})}\leq\|\partial_yu_{i\Phi_i}\|_{\widetilde{L}_t^2(B^{1,1})},\|\partial_y\theta_{i\Phi}\|_{\widetilde{L}_t^2(B^{1,1})}\leq\|\partial_y\theta_{i\Phi_i}\|_{\widetilde{L}_t^2(B^{1,1})},\\
&\|u_{i\Phi}\|_{\widetilde{L}_t^\infty(B^{\frac{3}{2},0})}\lesssim\|u_{i\Phi_i}\|_{\widetilde{L}_t^\infty(B^{1,0})},\|\theta_{i\Phi}\|_{\widetilde{L}_t^\infty(B^{\frac{3}{2},0})}\lesssim\|\theta_{i\Phi_i}\|_{\widetilde{L}_t^\infty(B^{1,0})}.\end{aligned}$$
These estimates have been presented as well in Section 4 of \cite{zhang2016long}.
\end{remark}

In this subsection, we will present the proof of estimates given in the following proposition, which are the key ingredients used in the proof of the uniqueness part of Theorem \ref{50}.

\begin{proposition}\label{62}
Suppose that $(U,\Theta)$ is a solution of (\ref{55}), for any $0<\eta<1$, it holds that
\begin{equation}\label{65}
\begin{aligned}
    &\|U_\Phi(t)\|_{B^{1,0}}+(\frac{\sqrt{\theta^E}}{2}-\epsilon^{\frac{1}{2}}\langle t\rangle^{\frac{1}{8}})\|\partial_y U_\Phi\|_{\widetilde{L}_t^2(B^{1,0})}+(\sqrt{\lambda}-C_\eta)\|U_\Phi\|_{\widetilde{L}_{t,\dot{M}}^2(B^{\frac{3}{2},0})}\\
    \lesssim&C_\eta\|U_\Phi\|_{\widetilde{L}_{t,\dot{M}}^2(B^{\frac{3}{2},0})}+\eta\|\partial_y U_\Phi\|_{\widetilde{L}_t^2(B^{1,0})}+C_\eta t^{\frac{1}{2}}(\|u_{1\Phi_1}\|_{\widetilde{L}_t^\infty(B^{1,0})}+\|u_{2\Phi_2}\|_{\widetilde{L}_t^\infty(B^{1,0})})\|U_\Phi\|_{\widetilde{L}_t^\infty(B^{1,0})}\\
    &+\eta(\|\partial_y(u_1+u_2)_\Phi\|_{\widetilde{L}_t^\infty(B^{1,0})}\|\partial_yU_\Phi\|_{\widetilde{L}_t^2(B^{1,0})}+\|\partial_y\Theta_\Phi\|_{\widetilde{L}_t^2(B^{1,0})}).
\end{aligned}
\end{equation}
and
\begin{equation}\label{69}
\begin{aligned}
    &\|\Theta_\Phi(t)\|_{B^{1,0}}+(\frac{\sqrt{\theta^E}}{2}-\epsilon^{\frac{1}{2}}\langle t\rangle^{\frac{1}{8}})\|\partial_y \Theta_\Phi\|_{\widetilde{L}_t^2(B^{1,0})}+(\sqrt{\lambda}-C_\eta)\|\Theta_\Phi\|_{\widetilde{L}_{t,\dot{M}}^2(B^{\frac{3}{2},0})}\\
    \lesssim&C_\eta(\|U_\Phi\|_{\widetilde{L}_{t,\dot{M}}^2(B^{\frac{3}{2},0})}+\|\Theta_\Phi\|_{\widetilde{L}_{t,\dot{M}}^2(B^{\frac{3}{2},0})})+\eta(\|\partial_y \Theta_\Phi\|_{\widetilde{L}_t^2(B^{1,0})}+\|\partial_y U_\Phi\|_{\widetilde{L}_t^2(B^{1,0})})\\
    &+C_\eta t^{\frac{1}{2}}(\|\theta_{1\Phi_1}\|_{\widetilde{L}_t^\infty(B^{1,0})}+\|u_{2\Phi_2}\|_{\widetilde{L}_t^\infty(B^{1,0})})\|\Theta_\Phi\|_{\widetilde{L}_t^\infty(B^{1,0})}+\eta\|\partial_y(u_1+u_2)_\Phi\|_{\widetilde{L}_t^\infty(B^{1,0})}\|\partial_yU_\Phi\|_{\widetilde{L}_t^2(B^{1,0})}\\
    &+(\sqrt{\epsilon}+\sqrt{\theta^E})(\eta\|\partial_yU_\Phi\|_{\widetilde{L}_t^2(B^{1,0})}
    +C_\eta\|\Theta_\Phi\|_{\widetilde{L}_{t,\dot{M}}^2(B^{\frac{3}{2},0})})
\end{aligned}
\end{equation}
\end{proposition}

We are now in a position to complete the proof of the uniqueness part of Theorem \ref{50}.

\begin{proof}[Proof of the uniqueness part of Theorem \ref{50}]
Combining estimates (\ref{65}) and (\ref{69}),  we obtain that
\begin{equation}
\begin{aligned}
    &\|U_\Phi(t)\|_{\widetilde{L}_t^\infty(B^{1,0})}+\|\Theta_\Phi(t)\|_{\widetilde{L}_t^\infty(B^{1,0})}+(\frac{\sqrt{\theta^E}}{2}-\epsilon^{\frac{1}{2}}\langle t\rangle^{\frac{1}{8}})(\|\partial_y U_\Phi\|_{\widetilde{L}_t^2(B^{1,0})}+\|\partial_y \Theta_\Phi\|_{\widetilde{L}_t^2(B^{1,0})})\\
    &+\sqrt{\lambda}(\|U_\Phi\|_{\widetilde{L}_{t,\dot{M}}^2(B^{\frac{3}{2},0})}
    +\|\Theta_\Phi\|_{\widetilde{L}_{t,\dot{M}}^2(B^{\frac{3}{2},0})})\\
    &\lesssim C_\eta(\|U_\Phi\|_{\widetilde{L}_{t,\dot{M}}^2(B^{\frac{3}{2},0})}+\|\Theta_\Phi\|_{\widetilde{L}_{t,\dot{M}}^2(B^{\frac{3}{2},0})})+\eta(\|\partial_y U_\Phi\|_{\widetilde{L}_t^2(B^{1,0})}+\|\partial_y \Theta_\Phi\|_{\widetilde{L}_t^2(B^{1,0})})\\
    &+C_\eta t^{\frac{1}{2}}\left(\|(u_{1\Phi_1}, u_{2\Phi_2})\|_{\widetilde{L}_t^\infty(B^{1,0})}\|U_\Phi\|_{\widetilde{L}_t^\infty(B^{1,0})} +(\|\theta_{1\Phi_1}\|_{\widetilde{L}_t^\infty(B^{1,0})}+\|u_{2\Phi_2}\|_{\widetilde{L}_t^\infty(B^{1,0})})\|\Theta_\Phi\|_{\widetilde{L}_t^\infty(B^{1,0})}\right)\\
&+\eta\|\partial_y(u_1+u_2)_\Phi\|_{\widetilde{L}_t^\infty(B^{1,0})}\|\partial_yU_\Phi\|_{\widetilde{L}_t^2(B^{1,0})}.
\end{aligned}
\end{equation}
Hence, by choosing $\eta$ sufficiently small and $\lambda$ large enough such that $\sqrt{\lambda}\gtrsim C_\eta$, in the above inequality, and using the boundedness of $u_i$ and $\theta_i$ ($i=1,2$) given by the a priori estimate (\ref{45}), we deduce that $\|U_\Phi(t)\|_{\widetilde{L}_t^\infty(B^{1,0})}+\|\Theta_\Phi(t)\|_{\widetilde{L}_t^\infty(B^{1,0})}=0$ for some small time $t$. Therefore, we obtain the uniqueness of solutions obtained in Theorem \ref{50}.
\end{proof}

\subsubsection{Estimate of $U$}
By acting the Fourier multiplier $e^{\Phi(t,D)}$ and applying $\Delta_k$ to the first equation given in \eqref{55}, then taking $L^2_{t,+}$ inner product of the resulting equation with $e^{2\Psi}\Delta_k U_\Phi$, we obtain 
\begin{equation}\label{53}
\begin{aligned}
    &(e^\Psi\partial_t\Delta_k U_\Phi|e^\Psi\Delta_k U_\Phi)+\lambda(\dot{M}\langle D\rangle e^\Psi\Delta_k U_\Phi|e^\Psi\Delta_k U_\Phi)-(e^\Psi\Delta_k[(\theta_2+\theta^E)\partial_y^2 U]_\Phi|e^\Psi\Delta_k U_\Phi)\\
    =&-(e^\Psi\Delta_k[U\partial_x u_1]_\Phi|e^\Psi\Delta_k U_\Phi)-(e^\Psi\Delta_k[u_2\partial_x U]_\Phi|e^\Psi\Delta_k U_\Phi)-(e^\Psi\Delta_k[V\partial_y u_1]_\Phi|e^\Psi\Delta_k U_\Phi)\\
    &-(e^\Psi\Delta_k[v_2\partial_y U]_\Phi|e^\Psi\Delta_k U_\Phi)+(e^\Psi\Delta_k[\Theta \partial_y^2 u_1]_\Phi|e^\Psi\Delta_k U_\Phi).
\end{aligned}
\end{equation}
Denote the three terms on the left side by $\mathcal{L}_1^k,\mathcal{L}_2^k,\mathcal{L}_3^k$ and the five terms on the right side by $\mathcal{R}_1^k,\cdots,\mathcal{R}_5^k$ respectively.

One can establish the following two lemmas by an argument similar to that given in Lemmas \ref{6}-\ref{13} and the proof of uniqueness part in \cite{zhang2016long}.

\begin{lemma}\label{63}
    For any $\eta>0$, it holds that
    $$\begin{aligned}
    \sum\limits_{k\in\mathbb{Z}}&2^k\left(\sqrt{|\mathcal{L}_1^k+\mathcal{L}_2^k+\mathcal{L}_{3}^k|}\right)\gtrsim\|U_\Phi(t)\|_{B^{1,0}}+(\sqrt{\lambda}-C_\eta)\|U_\Phi\|_{\widetilde{L}_{t,\dot{M}}^2(B^{\frac{3}{2},0})}\\
    &+\|\sqrt{-(\partial_t\Psi+4\theta^E(\partial_y\Psi)^2)}U_\Phi\|_{\widetilde{L}_t^2(B^{1,0})}
    +(\frac{\sqrt{\theta^E}}{2}-\epsilon^{\frac{1}{2}}\langle t\rangle^{\frac{1}{8}}-\eta)\|\partial_y U_\Phi\|_{\widetilde{L}_t^2(B^{1,0})}.\end{aligned}$$
\end{lemma}

\begin{lemma}\label{54}
    For any $\eta>0$, it holds that
    $$\begin{aligned}\sum\limits_{k\in\mathbb{Z}}&2^k\left(\sqrt{|\mathcal{R}_1^k|}+\sqrt{|\mathcal{R}_2^k|}+\sqrt{|\mathcal{R}_4^k|}\right)\lesssim C_\eta\|U_\Phi\|_{\widetilde{L}_{t,\dot{M}}^2(B^{\frac{3}{2},0})}+\eta\|\partial_y U_\Phi\|_{\widetilde{L}_t^2(B^{1,0})}\\
    &+C_\eta t^{\frac{1}{2}}(\|u_{1\Phi_1}\|_{\widetilde{L}_t^\infty(B^{1,0})}+\|u_{2\Phi_2}\|_{\widetilde{L}_t^\infty(B^{1,0})})\|U_\Phi\|_{\widetilde{L}_t^\infty(B^{1,0})}. \end{aligned}$$
\end{lemma}

\begin{lemma}\label{61}
    For any $\eta>0$, it holds that
    \begin{equation}\label{60}\sum\limits_{k\in\mathbb{Z}}2^k\sqrt{|\mathcal{R}_3^k|}\lesssim\eta(\|\partial_y(u_1+u_2)_\Phi\|_{\widetilde{L}_t^\infty(B^{1,0})}\|\partial_yU_\Phi\|_{\widetilde{L}_t^2(B^{1,0})}+\|\partial_y\Theta_\Phi\|_{\widetilde{L}_t^2(B^{1,0})})+C_\eta\|U_\Phi\|_{\widetilde{L}_{t,\dot{M}}^2(B^{\frac{3}{2},0})}\end{equation}

\end{lemma}

\begin{proof}
By using the third equation given in (\ref{55}), we have
$$
-(e^\Psi\Delta_k[V\partial_y u_1]_\Phi|e^\Psi\Delta_k U_\Phi)=\mathcal{R}_{3,1}^k+\mathcal{R}_{3,2}^k+\mathcal{R}_{3,3}^k$$
with
$$\begin{cases}
\mathcal{R}_{3,1}^k=(e^\Psi\Delta_k[\int_0^y \partial_x Udy'
\partial_y u_1]_\Phi|e^\Psi\Delta_k U_\Phi),\\
\mathcal{R}_{3,2}^k=-(e^\Psi\Delta_k[\int_0^y\partial_yU\partial_y(u_1+u_2)dy'\partial_y u_1]_\Phi|e^\Psi\Delta_k U_\Phi),\\
\mathcal{R}_{3,3}^k=-(e^\Psi\Delta_k[(
\partial_y\Theta)\partial_y u_1]_\Phi|e^\Psi\Delta_k U_\Phi)    
\end{cases}
$$
From a similar derivation of Lemma \ref{54}, we get
$$\sum\limits_{k\in\mathbb{Z}}2^k\sqrt{|\mathcal{R}_{3,1}^k|}\lesssim\|U_\Phi\|_{\widetilde{L}_{t,\dot{M}}^2(B^{\frac{3}{2},0})}.$$
By applying Bony's decomposition to $\displaystyle\int_0^y\partial_yU\partial_y(u_1+u_2)dy'\partial_y u_1$, we obtain
$$\begin{aligned}
    \mathcal{R}_{3,2}^k
    =&-(e^\Psi\Delta_k[T_{\int_0^y\partial_yU\partial_y(u_1+u_2)dy'}\partial_y u_1]_\Phi|e^\Psi\Delta_k U_\Phi)-(e^\Psi\Delta_k[T_{\partial_y u_1}\int_0^y\partial_yU\partial_y(u_1+u_2)dy']_\Phi|e^\Psi\Delta_k U_\Phi)\\
    &-(e^\Psi\Delta_k[R(\int_0^y\partial_yU\partial_y(u_1+u_2)dy',\partial_y u_1)]_\Phi|e^\Psi\Delta_k U_\Phi)\\
    =&:I_1^k+I_2^k+I_3^k.
\end{aligned}$$
Considering the support properties to the Fourier transform of the terms in $T_{\int_0^y\partial_yU\partial_y(u_1+u_2)dy'}\partial_y u_1$, we deduce that
$$\begin{aligned}
    |I_{1}^k|\lesssim&\int_0^t\int_0^\infty e^{2\Psi}\sum\limits_{|k'-k|\leq 4}\sum\limits_{k''\leq k'-2}2^{\frac{k''}{2}}\|\Delta_{k''}[\int_0^y\partial_yU\partial_y(u_1+u_2)dy']_\Phi\|_{L_x^2}\|\Delta_{k'}\partial_y u_{1\Phi}\|_{L_x^2}\|\Delta_kU_\Phi\|_{L_x^2}dydt'\\
    \lesssim&\sum\limits_{|k'-k|\leq 4}\int_0^t\|\partial_yU_\Phi\|_{B^{1,0}}\|\partial_y(u_1+u_2)_\Phi\|_{B^{1,0}}\|e^\Psi\Delta_{k'}\partial_y u_{1\Phi}\|_{L_+^2}\|e^\Psi\Delta_kU_\Phi\|_{L_+^2}dt'\\
    \leq&\|\partial_yu_{1\Phi}\|_{\widetilde{L}_t^\infty(B^{1,0})}2^{-\frac{3}{2}k}d_k\sum\limits_{|k'-k|\leq 4}2^{-k'}\widetilde{d}_{k'}\|\partial_yU_\Phi\|_{\widetilde{L}_t^2(B^{1,0})}\|U_\Phi\|_{\widetilde{L}_{t,\dot{M}}^2(B^{\frac{3}{2},0})},
\end{aligned}$$
with $${d_k}=\frac{2^{\frac{3}{2}k}\left(\int_0^t\|e^\Psi\Delta_k U_\Phi\|_{L_+^2}^2dt'\right)^{\frac{1}{2}}}{\|U_\Phi\|_{\widetilde{L}_t^2(B^{\frac{3}{2},0})}}$$ and $$\widetilde{d}_{k}=\frac{2^k\|e^\Psi\Delta_k\partial_yu_{1\Phi}\|_{L_t^\infty(L_+^2)}}{\|\partial_yu_{1\Phi}\|_{\widetilde{L}_t^\infty(B^{1,0})}}.$$
In the same way as above, we have that
$$\begin{aligned}
    |I_{2}^k|\lesssim&\int_0^t\int_0^\infty e^{2\Psi}\sum\limits_{|k'-k|\leq 4}\sum\limits_{k''\leq k'-2}2^{\frac{k''}{2}}\|\Delta_{k''}\partial_y u_{1\Phi}\|_{L_x^2}\|\Delta_{k'}[\int_0^y\partial_yU\partial_y(u_1+u_2)dy']_\Phi\|_{L_x^2}\|\Delta_kU_\Phi\|_{L_x^2}dydt'\\
    \lesssim&\left(\int_0^t\|\partial_yu_{1\Phi}\|_{B^{1,0}}^2\|e^\Psi\Delta_kU_\Phi\|_{L_+^2}^2\right)^\frac{1}{2}\sum\limits_{|k'-k|\leq 4}2^{-k'}\bar{d}_{k'}\|\partial_y(u_1+u_2)_\Phi\|_{\widetilde{L}_t^\infty(B^{1,0})}\|\partial_yU_\Phi\|_{\widetilde{L}_t^2(B^{1,0})} \\
    \lesssim&\|\partial_y(u_1+u_2)_\Phi\|_{\widetilde{L}_t^\infty(B^{1,0})}2^{-\frac{3}{2}k}d_k\sum\limits_{|k'-k|\leq 4}2^{-k'}\bar{d}_{k'}\|\partial_yU_\Phi\|_{\widetilde{L}_t^2(B^{1,0})}\|U_\Phi\|_{\widetilde{L}_{t,\dot{M}}^2(B^{\frac{3}{2},0})},
\end{aligned}$$
where $\bar{d}_k$ is properly defined such that $\sum_{k=1}^{+\infty}|\bar{d}_k|=1$ by using the following inequality, which can be obtained in a way similar to the proof of Lemma \ref{13},
$$\sum\limits_{k'\in\Z}2^{k'}\|\Delta_{k'}[\int_0^y\partial_yU\partial_y(u_1+u_2)dy']_\Phi\|_{L_x^2}\lesssim\|\partial_y(u_1+u_2)_\Phi\|_{B^{1,0}}\|\partial_yU_\Phi\|_{B^{1,0}}.
$$
Finally, considering the support properties to the Fourier transform of the terms in $R(\theta\partial_yu,\partial_y^2u)$, we deduce that
$$\begin{aligned}
    |I_{3}^k|\lesssim&\int_0^t\int_0^\infty e^{2\Psi}\sum\limits_{k'\geq k-3}\sum\limits_{k''= k'-1}^{k'+1}2^{\frac{k''}{2}}\|\Delta_{k''}[\int_0^y\partial_yU\partial_y(u_1+u_2)dy']_\Phi\|_{L_x^2}\|\Delta_{k'}\partial_y u_{1\Phi}\|_{L_x^2}\|\Delta_kU_\Phi\|_{L_x^2}dydt'\\
    \lesssim&\left(\int_0^t\|\partial_yu_{1\Phi}\|_{B^{1,0}}^2\|e^\Psi\Delta_kU_\Phi\|_{L_+^2}^2\right)^\frac{1}{2}\sum\limits_{k'\geq k-3}2^{-k'}\bar{d}_{k'}\|\partial_y(u_1+u_2)_\Phi\|_{\widetilde{L}_t^\infty(B^{1,0})}\|\partial_yU_\Phi\|_{\widetilde{L}_t^2(B^{1,0})}\\
    \lesssim&\|\partial_y(u_1+u_2)_\Phi\|_{\widetilde{L}_t^\infty(B^{1,0})}2^{-\frac{3}{2}k}d_k\sum\limits_{k'\geq k-3}2^{-k'}\bar{d}_{k'}\|\partial_yU_\Phi\|_{\widetilde{L}_t^2(B^{1,0})}\|U_\Phi\|_{\widetilde{L}_{t,\dot{M}}^2(B^{\frac{3}{2},0})}.
\end{aligned}$$
By applying Young's inequality and summing up the above three terms, we get
$$\sum\limits_{k\in\mathbb{Z}}2^k\sqrt{|\mathcal{R}_{3,2}^k|}\lesssim\eta\|\partial_y(u_1+u_2)_\Phi\|_{\widetilde{L}_t^\infty(B^{1,0})}\|\partial_yU_\Phi\|_{\widetilde{L}_t^2(B^{1,0})}+C_\eta\|U_\Phi\|_{\widetilde{L}_{t,\dot{M}}^2(B^{\frac{3}{2},0})}.$$
By applying Bony's decomposition to $\partial_y\Theta\partial_y u_1$, we obtain 
$$\begin{aligned}
    \mathcal{R}_{3,3}^k
    =&-(e^\Psi\Delta_k[T_{\partial_y u_1}\partial_y\Theta]_\Phi|e^\Psi\Delta_k U_\Phi)- (e^\Psi\Delta_k[T_{\partial_y\Theta}\partial_y u_1]_\Phi|e^\Psi\Delta_k U_\Phi)\\
    &-(e^\Psi\Delta_k[R(\partial_y\Theta,\partial_y u_1)]_\Phi|e^\Psi\Delta_k U_\Phi)\\
    =&:I_1^k+I_2^k+I_3^k.
\end{aligned}$$
By a similar derivation of Lemma \ref{54}, we find
$$\begin{aligned}
|I_1^k|\lesssim&\sum\limits_{|k'-k|\leq4}\int_0^t\langle t\rangle^{\frac{1}{4}}\|\partial_y^2 u_{1\Phi}\|_{B^{1,0}}\|e^\Psi\Delta_{k'}\partial_y\Theta_\Phi\|_{L_+^2}\|e^\Psi\Delta_k U_\Phi\|_{L_+^2}dt'\\
\lesssim&\sum\limits_{|k'-k|\leq4}\bigg(\int_0^t \dot{M}\|e^\Psi\Delta_k U_\Phi\|_{L_+^2}^2dt'\bigg)^{\frac{1}{2}}\bigg(\int_0^t\|e^\Psi\Delta_{k'}\partial_y\Theta_\Phi\|_{L_+^2}^2\bigg)^{\frac{1}{2}}.
\end{aligned}$$
Along the same line, we can deduce that 
$$\begin{aligned}
|I_2^k|\lesssim&\sum\limits_{k''\leq k+2}2^{\frac{k''-2k}{2}}\int_0^t\langle t\rangle^{\frac{1}{4}}\|\partial_y^2 u_{1\Phi}\|_{B^{1,0}}\|e^\Psi\Delta_{k''}\partial_y\Theta_\Phi\|_{L_+^2}\|e^\Psi\Delta_k U_\Phi\|_{L_+^2}dt'\\
\lesssim&\sum\limits_{k''\leq k+2}2^{\frac{k''-2k}{2}}\bigg(\int_0^t \dot{M}\|e^\Psi\Delta_k U_\Phi\|_{L_+^2}^2dt'\bigg)^{\frac{1}{2}}\bigg(\int_0^t\|e^\Psi\Delta_{k''}\partial_y\Theta_\Phi\|_{L_+^2}^2\bigg)^{\frac{1}{2}},
\end{aligned}$$
and
$$\begin{aligned}
|I_3^k|\lesssim&\sum\limits_{k'\geq k-3}\int_0^t\langle t\rangle^{\frac{1}{4}}\|\partial_y^2 u_{1\Phi}\|_{B^{1,0}}\|e^\Psi\Delta_{k'}\partial_y\Theta_\Phi\|_{L_+^2}\|e^\Psi\Delta_k U_\Phi\|_{L_+^2}dt'\\
\lesssim&\sum\limits_{k'\geq k-3}\bigg(\int_0^t \dot{M}\|e^\Psi\Delta_k U_\Phi\|_{L_+^2}^2dt'\bigg)^{\frac{1}{2}}\bigg(\int_0^t\|e^\Psi\Delta_{k'}\partial_y\Theta_\Phi\|_{L_+^2}^2\bigg)^{\frac{1}{2}}.
\end{aligned}$$
By applying Young's inequality and summing up the above three terms, we have
\begin{equation}\label{57}\sum\limits_{k\in\mathbb{Z}}2^k\sqrt{|\mathcal{R}_{3,3}^k|}\lesssim\eta\|\partial_y\Theta_\Phi\|_{\widetilde{L}_t^2(B^{1,0})}+C_\eta\|U_\Phi\|_{\widetilde{L}_{t,\dot{M}}^2(B^{\frac{3}{2},0})}.\end{equation}
Summing up the above estimates gives rise to (\ref{60}). This finishes the proof of Lemma \ref{61}.
\end{proof}

By a similar derivation as to that of (\ref{57}), we obtain the following lemma:
\begin{lemma}\label{64}
    For any $\eta>0$, it holds that    $$\sum\limits_{k\in\mathbb{Z}}2^k\sqrt{|\mathcal{R}_5^k|}\lesssim\eta\|\partial_y\Theta_\Phi\|_{\widetilde{L}_t^2(B^{1,0})}+C_\eta\|U_\Phi\|_{\widetilde{L}_{t,\dot{M}}^2(B^{\frac{3}{2},0})}.$$
\end{lemma}

\begin{proof}[Proof of The Estimate \eqref{65}]
From (\ref{53}), we have that
$$\sum\limits_{k\in\mathbb{Z}}2^k\sqrt{|\mathcal{L}_1^k+\mathcal{L}_2^k+\mathcal{L}_{3}^k|}\lesssim\sum\limits_{j=1}^5\sum\limits_{k\in\mathbb{Z}}2^k\sqrt{|\mathcal{R}_j^k|}.$$
By combining Lemmas \ref{63}-\ref{64}, we obtain the estimate (\ref{65}) immediately.
\end{proof}

\subsubsection{Estimate of $\Theta$}

By acting the Fourier multiplier $e^{\Phi(t,D)}$ and applying $\Delta_k$ to the second equation given in \eqref{55}, then taking $L^2_{t,+}$ inner product of the resulting equation with $e^{2\Psi}\Delta_k \Theta_\Phi$, we obtain 
\begin{equation}\label{70}
\begin{aligned}
    &(e^\Psi\partial_t\Delta_k \Theta_\Phi|e^\Psi\Delta_k \Theta_\Phi)+\lambda(\dot{M}\langle D\rangle e^\Psi\Delta_k \Theta_\Phi|e^\Psi\Delta_k \Theta_\Phi)-(e^\Psi\Delta_k[(\theta_2+\theta^E)\partial_y^2 \Theta]_\Phi|e^\Psi\Delta_k \Theta_\Phi)\\
    =&-(e^\Psi\Delta_k[U\partial_x \theta_1]_\Phi|e^\Psi\Delta_k \Theta_\Phi)-(e^\Psi\Delta_k[u_2\partial_x \Theta]_\Phi|e^\Psi\Delta_k \Theta_\Phi)-(e^\Psi\Delta_k[V\partial_y \theta_1]_\Phi|e^\Psi\Delta_k \Theta_\Phi)\\
    &-(e^\Psi\Delta_k[v_2\partial_y \Theta]_\Phi|e^\Psi\Delta_k \Theta_\Phi)+(e^\Psi\Delta_k[\Theta \partial_y^2 \theta_1]_\Phi|e^\Psi\Delta_k \Theta_\Phi)+(e^\Psi\Delta_k[\Theta (\partial_y u_1)^2]_\Phi|e^\Psi\Delta_k \Theta_\Phi)\\
    &+(e^\Psi\Delta_k[(\theta_2+\theta^E)\partial_y U\partial_y(u_1+u_2)]_\Phi|e^\Psi\Delta_k \Theta_\Phi).
\end{aligned}
\end{equation}
Denote the three terms on the left side by $\mathcal{L}_1^k,\mathcal{L}_2^k,\mathcal{L}_3^k$ and the seven terms on the right side by $\mathcal{R}_1^k,\cdots,\mathcal{R}_7^k$.

First, we can obtain the following two lemmas by a similar derivation of Lemmas \ref{63}-\ref{64} and Lemma \ref{19}.

\begin{lemma}\label{71}
    For any $\eta>0$, it holds that
    $$\begin{aligned}
    \sum\limits_{k\in\mathbb{Z}}&2^k\left(\sqrt{|\mathcal{L}_1^k+\mathcal{L}_2^k+\mathcal{L}_{3}^k|}\right)\gtrsim\|\Theta_\Phi(t)\|_{B^{1,0}}+(\sqrt{\lambda}-C_\eta)\|\Theta_\Phi\|_{\widetilde{L}_{t,\dot{M}}^2(B^{\frac{3}{2},0})}\\
    &+\|\sqrt{-(\partial_t\Psi+4\theta^E(\partial_y\Psi)^2)}\Theta_\Phi\|_{\widetilde{L}_t^2(B^{1,0})}
    +(\frac{\sqrt{\theta^E}}{2}-\epsilon^{\frac{1}{2}}\langle t\rangle^{\frac{1}{8}}-\eta)\|\partial_y \Theta_\Phi\|_{\widetilde{L}_t^2(B^{1,0})}.\end{aligned}$$
\end{lemma}

\begin{lemma}
    For any $\eta>0$, it holds that
    $$\begin{aligned}\sum\limits_{j=1}^6&\sum\limits_{k\in\mathbb{Z}}2^k\sqrt{|\mathcal{R}_j^k|}\lesssim C_\eta(\|U_\Phi\|_{\widetilde{L}_{t,\dot{M}}^2(B^{\frac{3}{2},0})}+\|\Theta_\Phi\|_{\widetilde{L}_{t,\dot{M}}^2(B^{\frac{3}{2},0})})+\eta(\|\partial_y \Theta_\Phi\|_{\widetilde{L}_t^2(B^{1,0})}+\|\partial_y U_\Phi\|_{\widetilde{L}_t^2(B^{1,0})})\\
    &+C_\eta t^{\frac{1}{2}}(\|\theta_{1\Phi_1}\|_{\widetilde{L}_t^\infty(B^{1,0})}+\|u_{2\Phi_2}\|_{\widetilde{L}_t^\infty(B^{1,0})})\|\Theta_\Phi\|_{\widetilde{L}_t^\infty(B^{1,0})}+\eta(\|\partial_y(u_1+u_2)_\Phi\|_{\widetilde{L}_t^\infty(B^{1,0})}\|\partial_yU_\Phi\|_{\widetilde{L}_t^2(B^{1,0})})
    \end{aligned}$$

\end{lemma}

\begin{lemma}\label{67}
    For any $\eta>0$, the term $\mathcal{R}_7^k$ given in \eqref{70} satisfies the following estimate
    \begin{equation}\label{66}
    \sum\limits_{k\in\mathbb{Z}}2^k\sqrt{|\mathcal{R}_7^k|}\lesssim(\sqrt{\epsilon}+\sqrt{\theta^E})(\eta\|\partial_yU_\Phi\|_{\widetilde{L}_t^2(B^{1,0})}+C_\eta\|\Theta_\Phi\|_{\widetilde{L}_{t,\dot{M}}^2(B^{\frac{3}{2},0})})
    \end{equation}
\end{lemma}

\begin{proof}
Let
$$\begin{aligned}
    \mathcal{R}_7^k=&(e^\Psi\Delta_k[(\theta_2+\theta^E)\partial_y U\partial_y(u_1+u_2)]_\Phi|e^\Psi\Delta_k \Theta_\Phi)\\
    =&(e^\Psi\Delta_k[\theta_2\partial_y(u_1+u_2)\partial_y U]_\Phi|e^\Psi\Delta_k \Theta_\Phi)+\theta^E(e^\Psi\Delta_k[\partial_y U\partial_y(u_1+u_2)]_\Phi|e^\Psi\Delta_k \Theta_\Phi)\\
    =&:I_1^k+I_2^k.
\end{aligned}$$
Applying Bony's decomposition to $\theta_2\partial_y(u_1+u_2)\partial_y U$, it follows
$$\begin{aligned}
I_1^k=&(e^\Psi\Delta_k[T_{\theta_2\partial_y(u_1+u_2)}\partial_y U]_\Phi|e^\Psi\Delta_k \Theta_\Phi)+(e^\Psi\Delta_k[T_{\partial_y U}\theta_2\partial_y(u_1+u_2)]_\Phi|e^\Psi\Delta_k \Theta_\Phi)\\
&+(e^\Psi\Delta_k[R(\theta_2\partial_y(u_1+u_2),\partial_y U)]_\Phi|e^\Psi\Delta_k \Theta_\Phi)\\
=&:I_{1,1}^k+I_{1,2}^k+I_{1,3}^k.
\end{aligned}$$
By a similar derivation as to that of (\ref{27}), (\ref{28}), we get
$$\sum\limits_{k\in\mathbb{Z}}2^k\|\Delta_k[\theta_2\partial_y(u_1+u_2)]_\Phi\|_{L_x^2}\lesssim\|\theta_{2\Phi_2}\|_{B^{1,1}}\|\partial_y(u_1+u_2)_\Phi\|_{B^{1,1}}.$$
Considering the support properties to the Fourier transform of the terms in $T_{\theta_2\partial_y(u_1+u_2)}\partial_y U$, we deduce that
$$\begin{aligned}
    |I_{1,1}^k|\lesssim&\int_0^t\int_0^\infty e^{2\Psi}\sum\limits_{|k'-k|\leq 4}\sum\limits_{k''\leq k'-2}2^{\frac{k''}{2}}\|\Delta_{k''}[\theta_2\partial_y(u_1+u_2)]_\Phi\|_{L_x^2}\|\Delta_{k'}\partial_y U_\Phi\|_{L_x^2}\|\Delta_k\Theta_\Phi\|_{L_x^2}dydt'\\
    \lesssim&\sum\limits_{|k'-k|\leq 4}\int_0^t\|\theta_{2\Phi_2}\|_{B^{1,1}}\|\partial_y(u_1+u_2)_\Phi\|_{B^{1,1}}\|e^\Psi\Delta_{k'}\partial_yU_\Phi\|_{L_+^2}\|e^\Psi\Delta_k\Theta_\Phi\|_{L_+^2}dt'\\
    \leq&\epsilon\sum\limits_{|k'-k|\leq 4}\left(\int_0^t\|e^\Psi\Delta_{k'}\partial_yU_\Phi\|_{L_+^2}^2dt'\right)^{\frac{1}{2}}\left(\int_0^t\dot{M}\|e^\Psi\Delta_k\Theta_\Phi\|_{L_+^2}^2dt'\right)^{\frac{1}{2}},
\end{aligned}$$
by using the smallness of $\theta$ given in \eqref{51}.
In the same manner, we infer that
$$\begin{aligned}
    |I_{1,2}^k|\lesssim&\int_0^t\int_0^\infty e^{2\Psi}\sum\limits_{|k'-k|\leq 4}\sum\limits_{k''\leq k'-2}2^{\frac{k''}{2}}\|\Delta_{k''}\partial_yU_\Phi\|_{L_x^2}\|\Delta_{k'}[\theta_2\partial_y(u_1+u_2)]_\Phi\|_{L_x^2}\|\Delta_k\Theta_\Phi\|_{L_x^2}dydt'\\
    \lesssim&\sum\limits_{k''\leq k+2}2^{\frac{k''}{2}-k}\int_0^t\|\theta_{2\Phi_2}\|_{B^{1,1}}\|\partial_y(u_1+u_2)_\Phi\|_{B^{1,1}}\|\Delta_{k''}\partial_yU_\Phi\|_{L_x^2}\|\Delta_k\Theta_\Phi\|_{L_x^2}dydt'\\
    \leq&\epsilon\sum\limits_{k''\leq k+2}2^{\frac{k''}{2}-k}\left(\int_0^t\|e^\Psi\Delta_{k''}\partial_yU_\Phi\|_{L_+^2}^2dt'\right)^{\frac{1}{2}}\left(\int_0^t\dot{M}\|e^\Psi\Delta_k\Theta_\Phi\|_{L_+^2}^2dt'\right)^{\frac{1}{2}}.
\end{aligned}$$
Finally, we deduce from the support properties to the Fourier transform of the terms in $R(\theta_2\partial_y(u_1+u_2),\partial_y U)$ that
$$\begin{aligned}
    |I_{1,3}^k|\lesssim&\int_0^t\int_0^\infty e^{2\Psi}\sum\limits_{k'\geq k-3}\sum\limits_{k''= k'-1}^{k'+1}2^{\frac{k''}{2}}\|\Delta_{k''}[\theta_2\partial_y(u_1+u_2)]_\Phi\|_{L_x^2}\|\Delta_{k'}\partial_y U_\Phi\|_{L_x^2}\|\Delta_k\Theta_\Phi\|_{L_x^2}dydt'\\
    \lesssim&\sum\limits_{k'\geq k-3}\int_0^t\|\theta_{2\Phi_2}\|_{B^{1,1}}\|\partial_y(u_1+u_2)_\Phi\|_{B^{1,1}}\|e^\Psi\Delta_{k'}\partial_yU_\Phi\|_{L_+^2}\|e^\Psi\Delta_k\Theta_\Phi\|_{L_+^2}dt'\\
    \leq&\epsilon\sum\limits_{k'\geq k-3}\left(\int_0^t\|e^\Psi\Delta_{k'}\partial_yU_\Phi\|_{L_+^2}^2dt'\right)^{\frac{1}{2}}\left(\int_0^t\dot{M}\|e^\Psi\Delta_k\Theta_\Phi\|_{L_+^2}^2dt'\right)^{\frac{1}{2}}.
\end{aligned}$$
On the other hand, in a similar proof of Lemma \ref{19} we can get that
    $$\sum\limits_{k\in\mathbb{Z}}2^k\sqrt{|I_2^k|}\lesssim\sqrt{\theta^E}(\eta\|\partial_yU_\Phi\|_{\widetilde{L}_t^2(B^{1,0})}+C_\eta\|\Theta_\Phi\|_{\widetilde{L}_{t,\dot{M}}^2(B^{\frac{3}{2},0})}).$$
Summing up the above estimates gives rise to the estimate (\ref{66}). We thus complete the proof of Lemma \ref{67}.
\end{proof}

\begin{proof}[Proof of The Estimate \eqref{69}]
Due to (\ref{70}), we infer that
$$\sum\limits_{k\in\mathbb{Z}}2^k\sqrt{|\mathcal{L}_1^k+\mathcal{L}_2^k+\mathcal{L}_{3}^k|}\lesssim\sum\limits_{j=1}^7\sum\limits_{k\in\mathbb{Z}}2^k\sqrt{|\mathcal{R}_j^k|}.$$
By combining Lemmas \ref{71}-\ref{67}, we obtain the estimate (\ref{69}) immediately.
\end{proof}

 \noindent
    {\bf Acknowledgments:}
    This research was supported by National Key R\&D Program of China under grant 2024YFA1013302, and National Natural Science Foundation of China under Grant Nos. 12331008, 12171317 and 12250710674. 

\bibliographystyle{plain}
\bibliography{ref.bib}

\end{document}